\documentclass[preprint,12pt]{elsarticle}
\usepackage{multirow}%
\usepackage{mathrsfs}%
\usepackage[title]{appendix}%
\usepackage{textcomp}%
\usepackage{manyfoot}%
\usepackage{booktabs}%
\usepackage[linesnumbered,ruled,vlined]{algorithm2e}
\usepackage{listings}%
\usepackage{lmodern}
\usepackage[section]{placeins}

\usepackage{amsmath,amssymb,amsfonts}%
\usepackage{amsthm} 
\usepackage{dsfont} 

\usepackage{mathtools}
\usepackage{enumitem} 

\usepackage{ragged2e}
\usepackage{soul}
\usepackage[dvipsnames]{xcolor}

\definecolor{cb1}{RGB}{228,26,28}
\definecolor{cb2}{RGB}{55,126,184}
\definecolor{cb3}{RGB}{77,175,74}
\definecolor{cb4}{RGB}{152,78,163}
\definecolor{cb5}{RGB}{255,127,0}

 \usepackage{subfig} 
 \usepackage{float}
 \usepackage{graphicx} 
 \usepackage{caption}
 \usepackage{import}
 \usepackage{url}
 \usepackage{pgfplots}

\usepgfplotslibrary{groupplots, colorbrewer}
\pgfplotsset{cycle list/Set1-5}
\usepackage{tikz}
\usetikzlibrary{matrix}
\usetikzlibrary{shapes.geometric, arrows.meta, positioning, calc, shadows.blur, fit, backgrounds, colorbrewer, plotmarks}
\usetikzlibrary{external}
\usetikzlibrary{spy} 
\tikzsetexternalprefix{figures/}
\newtheoremstyle{iremark}
  {\topsep}   
  {\topsep}   
  {\upshape}  
  {0pt}       
  {\bfseries}  
  {.}         
  {5pt plus 1pt minus 1pt} 
  {\thmname{#1}\thmnumber{ \itshape#2}\thmnote{ (#3)}} 
\theoremstyle{iremark}
\newtheorem{remark}{Remark}
\usepackage{svg}

\newcommand*\diff{\mathop{}\!\mathrm{d}}

\colorlet{shadecolor}{gray!15}
\usepackage{amsthm}
\usepackage{framed}
\usepackage{mdframed}

\newtheorem{theorem}{Theorem}
\newtheorem{lemma}{Lemma}
\newtheorem{definition}{Definition}
\newtheorem{problem}{Problem}
\newtheorem{proposition}{Proposition}

\newcommand{\fx}{\boldsymbol{x}}
\newcommand{\fu}{\boldsymbol{u}}
\newcommand{\fv}{\boldsymbol{v}}
\newcommand{\fg}{\boldsymbol{g}}

\newcommand{\ff}{\boldsymbol{f}}

\newcommand{\ffC}{\mathbb{C}}
\newcommand{\rC}{\mathrm{C}}
\newcommand{\rG}{\mathrm{G}}

\usepackage{silence}
\newcommand{\primed}[1]{#1'}
\usepackage{hyperref}
\pgfplotsset{compat=1.17}
\pgfplotsset{select coords between index/.style 2 args={
    x filter/.code={
        \ifnum\coordindex<#1\fi
        \ifnum\coordindex>#2\fi
    }
}}
\newcounter{robplotcount}
\newcommand{\makeRobOptimizationPlots}[3]{%
    \stepcounter{robplotcount}
    \begin{tikzpicture}
        \def\legendrows{}
        \foreach \thefile [count=\i] in {#1} {
            \def\thislabel{}
            \foreach \mylabel [count=\j] in {#3} {
                \ifnum\i=\j \xdef\thislabel{\mylabel}\fi
            }
            \ifnum\i=1
                \xdef\legendrows{\protect\ref*{plotline:\arabic{robplotcount}:\thefile} \& \thislabel}
            \else
                \xdef\legendrows{\legendrows \& $\quad$ \& \protect\ref*{plotline:\arabic{robplotcount}:\thefile} \& \thislabel}
            \fi
        }

        \begin{groupplot}[
            group style={
                group size=2 by 2,
                horizontal sep=2cm,
                vertical sep=1cm
            },
            width=0.45\textwidth,
            height=5.5cm,
            xlabel={Iteration},
            grid=major,
            every axis plot/.append style={thick, no marks, smooth},
            label style={font=\footnotesize},
            tick label style={font=\footnotesize},
            legend image code/.code={
                \draw[-] (0cm,0.0ex) -- (0.7cm,0.0ex);
            }            
        ]

        \nextgroupplot[ylabel={Objective Function}]
        \foreach \thefile in {#1} {
            \addplot+ [no marks, smooth, select coords between index={0}{#2}] 
            table[x=ITER, y=COST, col sep=comma]{\thefile};
            \label{plotline:\arabic{robplotcount}:\thefile}
        }

        \nextgroupplot[ylabel={Volume Frac.}]
        \foreach \thefile in {#1} {
            \addplot+ [no marks, smooth, select coords between index={0}{#2}] 
            table[x=ITER, y=VOLUME, col sep=comma]{\thefile};
        }

        \nextgroupplot[ylabel={Expected Compliance}]
        \foreach \thefile in {#1} {
            \addplot+ [no marks, smooth, select coords between index={0}{#2}] 
            table[x=ITER, y=COMPLIANCE, col sep=comma]{\thefile};
        }

        \nextgroupplot[ylabel={Std. Deviation}]
        \foreach \thefile in {#1} {
            \addplot+ [no marks, smooth, select coords between index={0}{#2}] 
            table[x=ITER, y=SDEV, col sep=comma]{\thefile};
        }
        \end{groupplot}

        \matrix [
            draw, 
            matrix of nodes, 
            anchor=north, 
            at={($(group c1r2.south)!0.5!(group c2r2.south)$)}, 
            yshift=-1cm, 
            column sep=1ex, 
            nodes={anchor=base,           
            font=\footnotesize,    
            text height=1.2ex,     
            text depth=0.3ex,      
            inner xsep=2pt,         
            inner ysep=0pt
            },
            ampersand replacement=\& 
        ] (legend) {
            \legendrows \\
        };
    \end{tikzpicture}
}

\newcounter{example}[section]
\newenvironment{example}[1][]{%
   \refstepcounter{example}%
   \setcounter{subexample}{0}
   \par\medskip\noindent
   \textbf{Example~\theexample. #1}%
   \par\nopagebreak\smallskip
}{%
   \par\medskip
}

\newcounter{subexample}[example]
\newcommand{\subex}[1][]{%
  \refstepcounter{subexample}%
  \par\medskip\noindent
  \textit{#1}%
  \par\nopagebreak\smallskip
}
\newbool{showimages}
\setbool{showimages}{true} 

\newcommand{\placeholder}[3][width=\textwidth]{%
    \ifbool{showimages}{%
        \begin{figure}[h]
            \centering
            \includegraphics[#1]{#2}
            \caption{#3}
        \end{figure}
    }{%
    }%
}

\begin{document}
\begin{frontmatter}
\title{Robust topology optimization of isogeometric Kirchhoff--Love shells by the Speed Method}
\author[poly]{Philipp Eisenhardt \fnref{fn1}}
\author[poly]{Fernando Hübner Scherer \fnref{fn1}}

\fntext[1]{Both authors contributed equally to this work.}

\author[tum]{Barbara Wohlmuth }

\author[poly]{Andrei Constantinescu \corref{cor1}}
\cortext[cor1]{andrei.constantinescu@polytechnique.edu}

\address[poly]{
    Laboratoire de Mécanique des Solides, CNRS, École Polytechnique, IP Paris,
    91128,
    Palaiseau,
    France
}
\address[tum]{
    Institute for Numerical Mathematics, 
    TU Munich, 
    Boltzmannstr. 3,
85748 Garching, Germany
}
\begin{abstract}
This paper presents a framework for robust topology optimization of Kirchhoff–Love (KL) multi-patch shells using isogeometric analysis and non-matching meshes and provides a CAD-compatible method. 

The material distribution is optimized via a NURBS-based level set method, enabling accurate boundary representation and seamless export to standard CAD/CAM workflows. The shape derivative is computed using the speed method, yielding a volume-based formulation that avoids boundary approximations and aligns with the isogeometric principle of coarse, high-order discretizations. This eliminates the need for local refinement, adaptive quadrature, or trimming as the descent direction for the compliance and volume objectives are obtained through a biharmonic equation.

The framework is extended to robust optimization under uncertainty, accounting for stochastic loading and material properties via Monte Carlo sampling and Gaussian Random Fields. Numerical examples demonstrate the method’s efficiency and robustness across single- and multi-patch shells, reproducing literature results with fewer elements than traditional FEM or refined IGA approaches. 
The developed methods and tools can be directly extended to problems involving stress-based optimization and risk-averse metrics. 

The numerical implementations are open-source and available on GitHub.
\end{abstract}
\end{frontmatter}

\section{Introduction}

Thin Kirchhoff--Love (KL) shells are widely used in engineering applications, including civil structures \cite{bathe_geometric_1980,cirak_subdivision_2000,kiendl_isogeometric_2009,bouclier_efficient_2013}, biological tissues \cite{tepole_isogeometric_2015}, and complex material behaviors such as elasto-plasticity and phase-field fracture \cite{ambati_isogeometric_2018,kiendl_phase-field_2016}. The $C^1$-continuity requirement for KL shells is naturally satisfied by Isogeometric Analysis (IGA) \cite{hughes_isogeometric_2005,bazilevs_isogeometric_2006,cottrell_isogeometric_2009}, which integrates with finite element frameworks via Bézier extraction \cite{borden_isogeometric_2011,scott_isogeometric_2011,hennig_bezier_2016,kamensky_tigar_2019,yu_three-dimensional_2019,zhao_automated_2024}.

This paper presents a CAD-compatible topology optimization framework for KL multi-patch shells using a patch-wise isogeometric paradigm with non-matching meshes. We employ a NURBS level set method to implicitly represent evolving shapes, ensuring full CAD compatibility throughout the optimization as the shape at every iteration is encoded by the zero level set. The shape derivative is computed using the speed method, yielding a volume-based formulation that avoids boundary approximations and aligns with IGA’s philosophy of coarse, high-order discretizations. This eliminates the need for trimming, local refinement, or adaptive quadrature, while requiring the solution of a biharmonic equation to obtain a descent direction. We extend the framework to robust optimization under loading and material uncertainty using Monte Carlo sampling and Gaussian Random Fields. Numerical examples demonstrate that our method reproduces literature results with fewer elements than comparable FEM or refined IGA techniques.

Previous work on shell optimization includes parametric methods \cite{wall_isogeometric_2008,qian_full_2010,hsu_high-performance_2011,kiendl_isogeometric_2014}, Free Form Deformation (FFD) \cite{sederberg_free-form_1986,hirschler_embedded_2019,zhao_automated_2024}, and the Adaptive Bubble Method (ABM) \cite{cai_integrated_2023}. Level set topology optimization \cite{allaire_level-set_2005,de_gournay_velocity_2006,de_gournay_shape_2008,amstutz_topological_2010,laurain_level_2017,nika_design_2019,agnelli_design_2022} has been applied to 2D elasticity \cite{wang_isogeometric_2016,jahangiry_isogeometrical_2017,xu_level_2019,teixeira_isogeometric_2025} and Reissner--Mindlin shells \cite{hubner_scherer_topology_2024,hubner_scherer_stress-based_2024}. For KL shells, \cite{xing_isogeometric_2026} used boundary-based shape derivatives requiring fine refinement. Robust optimization approaches include \cite{kharmanda_reliability-based_2004,guo_confidence_2009,allaire_linearized_2014,de_topology_2020,dapogny_distributionally_2025,antil_risk-averse_2026,kodakkal_risk-averse_2022,sigmund_manufacturing_2009,chen_new_2011, baumgarten2026risk}.

The paper is structured as follows: Section~\ref{sec:speedMeth} recalls the Kirchhoff--Love shell model on surfaces and derives the shape derivative for shell compliance using the speed method for both deterministic and uncertainty-based optimization, along with its isogeometric implementation. The optimized geometries for a wide variety of single and multi-patch problems are shown in Section~\ref{sec:numRes}.

\section{Derivation of the Shape Derivatives for Kirchhoff--Love shells by the speed method}
\label{sec:speedMeth}
The Kirchhoff--Love (KL) model discussed next is structured along the presentation in \cite{stein_models_2004, ciarlet_introduction_2006, chapelle_finite_2010} and relies on concepts of differential geometry from \cite{frankel_geometry_2004, spivak_comprehensive_1979, abraham_manifolds_1988, ciarlet_introduction_2006}.
In our description, the mid-surface KL shell occupies a domain $\Omega$, embedded in a surface $\mathcal{S}$. 
\paragraph{Notation} Common differential geometry notation is adopted throughout this paper. Greek indices and exponents $\alpha, \beta$ have values in $\{1, 2\}$, Latin lowercase indices $i,j,k,l$ range in $\{1, 2, 3\}$. Einstein notation over repeated indices is adopted. Upper indices denote the contra-variant components of a vector or tensor, while lower indices its covariant components. The partial derivatives of a field $\mathbf{X}(\xi^1, \xi^2)$ with respect to the parametric coordinates are denoted by $\mathbf{X}_{,\alpha} = \frac{\partial \mathbf{X}}{\partial \xi^\alpha}$.\\
\paragraph{Kirchhoff--Love shell formulation on parametrized surfaces}

The modeled shell occupies a volume of thickness $h$ around its mid-surface $\Omega$ and is embedded in  a two-dimensional surface $\mathcal{S} \subset \mathbb{R}^3$. $\mathcal{S}$ is described by the $C^2$-smooth vector field $\mathbf{X}$ mapping each point of  the two-dimensional parametric domain $\boldsymbol{\xi} \in \hat{\mathcal{S}}\subset \mathbb{R}^2$ into a material point $\mathbf{X}(\boldsymbol{\xi}) \in \mathcal{S}$:
\begin{equation}
\begin{aligned}
\mathbf{X} : \quad \hat{\mathcal{S}}  \to \mathcal{S} 
\qquad 
\boldsymbol{\xi} \mapsto \mathbf{X}(\boldsymbol{\xi}), 
\end{aligned}
\label{eqn:def_parameterization}
\end{equation}
 Let us further notice that the material mid-surface $\Omega$ is the image of the domain $\hat{\Omega}\subset\hat{\mathcal{S}}$ through the mapping $\mathbf{X}$ and embedded within $\mathcal{S}$. The thickness of the shell $d$ is measured at each point on the mid-surface along the normal vector.  
\\
The geometry of the surface $\mathcal{S}$ is defined by the parametric mapping $\mathbf{X}$. $\boldsymbol{X}$ and its inverse $\boldsymbol{X}^{-1}$ will permit to transport all quantities of interest between the parametric and the material surfaces, $\hat{\mathcal{S}}$  and $\mathcal{S}$, respectively. Moreover, its gradient permits defining the Jacobian, the associated metric tensor and the metric determinant as:
\begin{align}
\mathbf{\hat{F}} = \nabla_{\mathbf{\xi}} \mathbf{X}
&& 
\mathbf{\hat{G}} = \mathbf{\hat{F}}^T \mathbf{\hat{F}}
&& 
\mathbf{\hat{g}} = \sqrt{\operatorname{det} \mathbf{\hat{G}}}.
\label{eqn:def_jac_met_det}
\end{align}
The covariant vector basis $\mathbf{A}$ of the surface defines the tangent and normal vectors on the surface:
\begin{equation}
    \mathbf{A}_\alpha = \mathbf{X}_{,\alpha} \qquad \text{ and } \qquad
    \mathbf{A}_3 = \frac{\mathbf{A}_1 \times \mathbf{A}_2}{|\mathbf{A}_1 \times \mathbf{A}_2|}.
\end{equation}
Under the Kirchhoff--Love kinematic assumption, material points situated on the normal to the initial mid-surface remain on the normal of the deformed mid-surface of the shell, without deformation along the normal. Hence, the shell deformation is completely characterized by the displacement field $\boldsymbol{u}$ of the mid-surface. As a consequence, the material points of the mid-surface of deformed shell $\mathbf{x}$ and the  associated covariant basis vectors $\mathbf{a}_\alpha$ on this new surface are obtained by adding the displacement field $\boldsymbol{u}$ of the mid-surface to initial mid-surface $\mathbf{X}$:
\begin{align}
    \mathbf{x} = \mathbf{X} + \boldsymbol{u}
    && 
    \mathbf{a}_\alpha =  \mathbf{x}_{,\alpha} = \mathbf{X}_{,\alpha} + \boldsymbol{u}_{,\alpha} 
\end{align}
Furthermore, the strain tensor can be linearized at small displacements and split into \emph{membrane strain} $\varepsilon_{\alpha \beta}$ and \emph{bending strain} $\kappa_{\alpha \beta}$:
\begin{align}
    \varepsilon_{\alpha \beta}(\boldsymbol{u}) &= \frac{1}{2} \left(\mathbf{A}_{\alpha} \cdot \boldsymbol{u}_{,\beta} + \mathbf{A}_{\beta} \cdot \boldsymbol{u}_{,\alpha} \right),
    \label{eqn:membraneStrainSmallDispl}\\
    \kappa_{\alpha \beta}(\boldsymbol{u}) &= -\mathbf{A}_3\cdot \left(\boldsymbol{u}_{,\alpha\beta} - \Gamma^\lambda_{\alpha\beta} \boldsymbol{u}_{,\lambda}\right),
    \label{eqn:bendingStrainSmallDispl}
\end{align}
where $\Gamma^\lambda_{\alpha\beta}$ corresponds to the Christoffel symbol of second kind associated with the mid-surface (see \cite{ciarlet_introduction_2008} for the derivation and more details). 

For isotropic linear elastic materials, the Hooke tensor of elastic moduli $\mathbb{C}$ is defined by the Young's modulus,  $E$, the Poisson coefficient $\nu$ and the metric tensor $\mathbf{G}$ \cite[equation 4.34]{chapelle_finite_2010}. Its contravariant coefficients are defined as:   
\begin{equation}
\rC^{\alpha\beta\lambda\mu} = \frac{E}{2(1+ \nu)}\left(\rG^{\alpha\lambda}\rG^{\beta\mu} + \rG^{\alpha\mu}\rG^{\beta\lambda} +\frac{2\nu}{1-\nu} \rG^{\alpha\beta}\rG^{\lambda\mu}  \right),
\label{eqn:consitutiveRelationWithBasis}
\end{equation}

\medskip

The elasticity problem can be defined as follows: 
\begin{problem}[Linear Elasticity on KL shells: $\mathcal{P}_{KL}(\boldsymbol{g}, \boldsymbol{u}_D)$]
    Find a displacement field $\boldsymbol{u}$ such that:
    \begin{equation}
        a_\Omega(\boldsymbol{u}, \boldsymbol{v}) = L(\boldsymbol{v}), \qquad \forall \boldsymbol{v}
    \end{equation}
    where the operators $a_\Omega$ and $L$ denote the virtual work of internal and external forces and are expressed as:
    \begin{equation}
    \begin{aligned}
        a_\Omega(\boldsymbol{u}, \boldsymbol{v}) &:= \int_{\Omega} \left( d \varepsilon(\boldsymbol{u}) :\ffC  : \varepsilon(\boldsymbol{v}) + \frac{d^3}{12}\kappa(\boldsymbol{u}):\ffC : \kappa(\boldsymbol{v}) \right) \quad \text{ and } \\
        L_{\Gamma_N}(\boldsymbol{v})&:= \int_{\partial \Omega_N} \boldsymbol{g} \cdot \boldsymbol{v}, 
    \end{aligned}
    \label{eqn:linElasKL}
    \end{equation}
  where $\boldsymbol{g}$ corresponds to a given traction field on the boundary $\partial \Omega_N$. 
  Moreover, the displacement field $\boldsymbol{u}$ has to satisfy the boundary condition $ \boldsymbol{u} = \boldsymbol{u}^D \text{ on }\partial\Omega_D$
  and  $\partial \Omega_D$ and $\partial \Omega_N$ are disjoint and complementary parts of   $\partial \Omega$, denoted as the Dirichlet and Neumann region, respectively.
    \label{prob:linElasKLOmega}
\end{problem}

\subsection{Shape optimization problem}
The shape optimization problem of the Kirchhoff--Love shells, seeks to determine a domain $\Omega$ embedded in  $\mathcal{S}$ optimizing the performance of the shell as measured by the compliance $\mathcal{J}_C$ and/or the volume function $\mathcal{J}_V$ under given loading conditions of the shell. These cost functionals are defined as:
\begin{align}
    \mathcal{J}_C(\Omega, \boldsymbol{u}) &=  a_\Omega(\boldsymbol{u}, \boldsymbol{u}) , \label{eqn:defCompCost}\\
    \mathcal{J}_V(\Omega) &=\int_\Omega 1, \label{eqn:defVolCost}
\end{align}
where $\boldsymbol{u}$ is the solution of the corresponding elasticity problem $\mathcal{P}_{KL}(\boldsymbol{g}, \boldsymbol{u}_D)$ \ref{prob:linElasKLOmega}.
\medskip
\begin{problem}[Optimization Problem: $\mathcal{P}_{O}(\mathcal{J})$]
\label{eq:optProb}
Find the optimal domain $\Omega^* \subset \mathcal{S}$ such that 
\begin{align}
    \Omega^* &= \underset{\Omega}{\operatorname{argmin}}\ \mathcal{J}(\Omega,\boldsymbol{u} ) 
    \\
    &=  \underset{\Omega}{\operatorname{argmin}} \ \left(\mathcal{J}_C(\Omega,   \boldsymbol{u})
    + \Lambda \mathcal{J}_V(\Omega) \right )
\end{align}
with  $\Lambda >0$ a fixed parameter and $\boldsymbol{u}\in \mathcal{P}(\boldsymbol{g}, \boldsymbol{u}_D)$. \\
In addition to solving the optimization with a fixed Lagrange multiplier $\Lambda$, a specific target volume fraction can be obtained through an augmented Lagrangian approach. Here, the objective function is augmented with a quadratic penalty term, as 
$ \Lambda_k \left(\mathcal{J}_V(\Omega)-V^* \right)+ \frac{\mu_k}{2} \left(\mathcal{J}_V(\Omega)-V^* \right)^2$, where $V^*$ corresponds to the target volume fraction and $\mu$ to an additional penalty parameter. 
The multipliers are updated dynamically across iterations to strictly enforce the volume constraint (see \cite{nocedal_numerical_2006, birgin_practical_2014} for the derivation, which directly extends to the current setting).
\end{problem}

\paragraph{Shape Derivative by the Speed Method}

The shape derivative is computed using the speed method, following \cite{sokolowski_introduction_1992, delfour_shapes_2011, laurain_level_2017} and the previous definition and notations for the Kirchhoff--Love shell. 

Let us consider  an evolution of the domain $\Omega$ as a function of a fictitious time $t$ designated by $\Omega_t$. Naturally, this will induce in the parametric space  a  map  $\hat{T}_t$ from  $\hat{\Omega}=\hat{\Omega}_0$ to $\hat{\Omega}_t$: $\hat{T}_t(\boldsymbol{\xi}) := \boldsymbol{\xi}_t$ and define an associated velocity field $\mathbf{\hat{\vartheta}}$ with Jacobian determinant $\psi$:
\begin{align*}
    \mathbf{\hat{\vartheta}}(t, \xi) = \frac{\diff}{\diff t}\xi_t && \hat{\psi}(t) = \text{det}\left(\left[\frac{\partial (\xi_t)^\lambda}{\partial \xi^\mu}\right]\right)
\end{align*}
Through the surface mapping $\mathbf{X}$, this parametric velocity is pushed forward to the material domain, defining a 3D tangent velocity field $\boldsymbol{\theta}$ on the surface $\mathcal{S}$:
\begin{align*}
    \hat{\boldsymbol{\theta}} := \hat{\vartheta}^\alpha\hat{\mathbf{A}}_\alpha, &&
    \boldsymbol{\theta} := (\hat{\vartheta}^\alpha\hat{\mathbf{A}}_\alpha) \circ \boldsymbol{X}^{-1}.
\end{align*}

The shape derivative of $\mathcal{J}$ for the shape $\Omega$  in direction $\boldsymbol{\theta}$ is defined by (see \cite{sokolowski_introduction_1992,delfour_shapes_2011}):
\begin{equation}
d\mathcal{J}(\Omega) [\boldsymbol{\theta}] := \lim_{t \searrow 0} \frac{\mathcal{J}(\Omega_t)-\mathcal{J}(\Omega)}{t}
\label{eqn:eulShapeDer}
\end{equation}

\begin{figure}[h]
  \centering
\includegraphics[width=\linewidth]{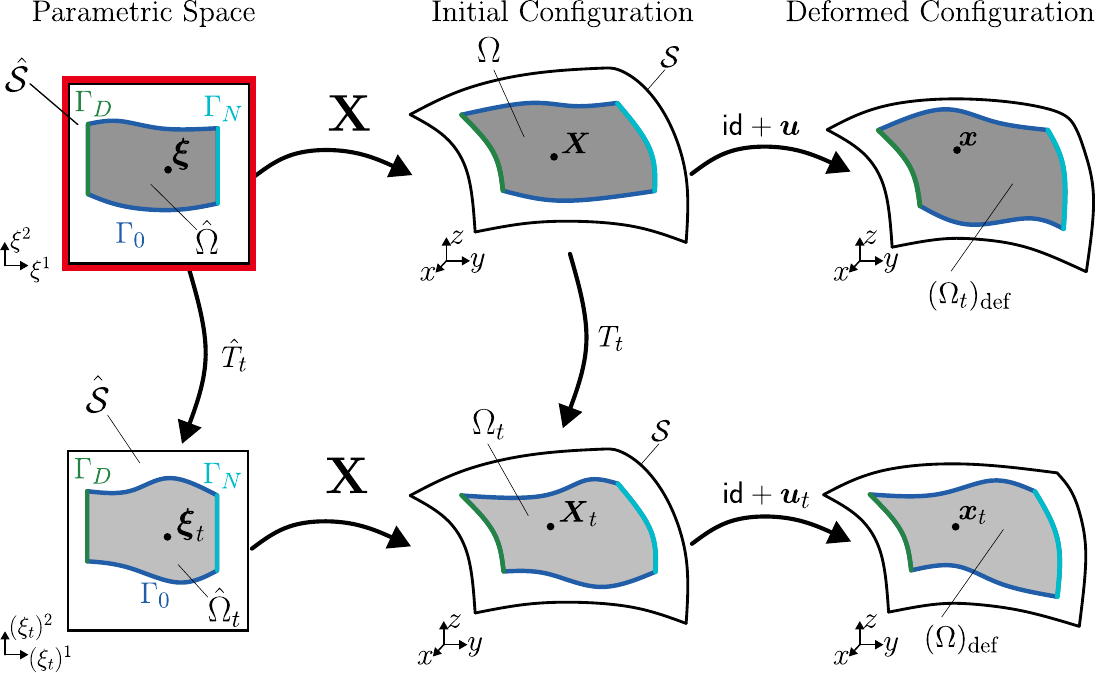}
  \caption{
Schematic view of the parametric-to-material-to-deformed domain mapping for shape perturbation in Kirchhoff-–Love shells. On the upper row:  the parametric domain  $\hat{\Omega}$  (left) the initial and deformed material configurations of  $\Omega$. On the lower row the corresponding variations of $\Omega_t$  (blue) under a velocity field $\boldsymbol{\vartheta}$.
}
  \label{fig:schemaShapes}
\end{figure}
To evaluate this shape derivative, the state equation evaluated on the perturbed domain $\Omega_t$ must be differentiated. Because the integration domain $\Omega_t$ depends on $t$, the integrals and the perturbed state variables $\boldsymbol{u}_t$ are pulled back to the fixed reference configuration $\Omega$ using the transformation $T_t$ as $\hat{\boldsymbol{u}}^t = \hat{\boldsymbol{u}}_t \circ T_t$. This allows for the computation of the material derivatives of the state variables and the subsequent application of the chain rule to the cost functional.\\
In general, we define the time derivative of a variable $\ff$ as:
\begin{equation}
    \ff^{t\prime} = \frac{\partial}{\partial t}\ff^t.
\end{equation}

\paragraph{Mechanical problem on perturbed shape \texorpdfstring{$\Omega_t$}{Ω}} 
Let us further assume that the boundary for the mechanical problem $\mathcal{P}(\boldsymbol{g}, \boldsymbol{u}^D)$ can be partitioned into three complementary regions: $\Gamma_D = \partial \Omega_D$ with imposed displacements, $\Gamma_N$ with prescribed traction $\boldsymbol{g}$, and a free surface $\Gamma_0$ where $\boldsymbol{g} = \boldsymbol{0}$ to give $\partial \Omega_N=\Gamma_N \cup \Gamma_0$.
The Dirichlet boundary $\Gamma_D$ and traction boundary $\Gamma_N$ will remain fixed by enforcing $\boldsymbol{\theta} =\boldsymbol{0}$.

Under this assumption, the mechanical problem is mapped to the perturbed domain $\Omega_t$, as illustrated in Figure~\ref{fig:schemaShapes} by the transformation $T_t$.
The mechanical problem~\ref{prob:linElasKLOmega} $\mathcal{P}(\fg, \fu_D)$ is now defined by:
\begin{equation}
        \underbrace{\int_{\Omega_t} d \, \varepsilon(\boldsymbol{u}_t) :  \mathbb{C}_t : \varepsilon (\boldsymbol{v}_t) + \frac{d^3}{12} \kappa(\boldsymbol{u}_t) :  \mathbb{C}_t : \kappa(\boldsymbol{v}_t)}_{a_{\Omega_t}} = \underbrace{\int_{\Gamma_N} \boldsymbol{g} \cdot \boldsymbol{v}_t}_{{\Gamma_N}}, 
    \label{prob:linelasticity}
\end{equation}
where $\hat{(\bullet)}_t$ corresponds to functions defined on the parametric domain $\hat{\Omega}_t$. 

\goodbreak

The problem can subsequently be rewritten over the original shape $\Omega$ and its parametric equivalent $\hat{\Omega}$ by using the transformation:
\begin{equation}
    \hat{\boldsymbol{u}}^t = \hat{\boldsymbol{u}}_t \circ T_t,
\end{equation}
where $\hat{(\cdot)}^t$ corresponds to properties in the perturbed configuration, rewritten over the original shape. 
The transformed expressions are derived in \autoref{sec:technicalLemma}.
With analogous formulae, the mechanical problem can hence be rewritten for the solution over the deformed shape as $\boldsymbol{u}^t$:
    \begin{equation}
        \int_{\Omega} \left( d \varepsilon^t(\boldsymbol{u}^t):\ffC^t  : \varepsilon^t(\boldsymbol{v}^t) + \frac{d^3}{12} \kappa^t(\boldsymbol{u}^t):\ffC^t : \kappa^t(\boldsymbol{v}^t) \right) \psi(t)= \int_{\partial\Omega_N} \boldsymbol{g}^t \cdot \boldsymbol{v}^t\psi(t)
        \label{eq:linElasPerturb}
    \end{equation}
\goodbreak
As the different problems have been rewritten with respect to a perturbation of the material phase $\Omega$, Equation~\ref{eqn:eulShapeDer} can be applied to compute the shape derivative of the different terms. 
Based on this, in the following expressions for the sensitivity of the \emph{volume function} and the \emph{compliance function} are developed.
\goodbreak
\noindent\textbf{Shape Derivative of the Volume function}\\
Applying these findings, the shape derivative of the volume of the domain can be computed:
The cost function
\begin{equation}
    \mathcal{J}_V(\Omega) = \int_{\Omega} 1
    \label{eqn:costVol}
\end{equation}
admits the following shape derivative
\begin{equation}
   d\mathcal{J}_V(\Omega)[\boldsymbol{\theta}]=\int_\Omega \operatorname{div} (\boldsymbol{\theta})
   \label{eqn:shapederivvol}
\end{equation}
This result is demonstrated in Appendix~\ref{sec:demSDvol} using the definition of the shape derivative from Definition~(\ref{eqn:eulShapeDer}), and the definitions of differential operators in the surface. 
\smallskip

\noindent\textbf{Shape Derivative of the Compliance Function} \\
To compute the shape derivative of the compliance function, we introduce a Lagrangian to account for the change of the displacement variable under a change of geometry as it is common in level set topology optimization (see for example \cite{allaire_conception_2007, allaire_level-set_2005, laurain_distributed_2015} and the references therein).
It thus consists of the cost function and the variational formulation and depends on the displacement $\boldsymbol{u}$ and the Lagrange multiplier of the state equation $\boldsymbol{v}$. 
We define the Lagrangian $\mathcal{L}$, which is written as:
\begin{equation}
\begin{aligned}
    \mathcal{L}(t, \boldsymbol{w}^t, \boldsymbol{v}^t) =&  \mathcal{J}_C(\boldsymbol{w}^t, \Omega_t) &&\\ 
    +&a_{\Omega_t}(\boldsymbol{w}^t, \boldsymbol{v}^t)-L_{\Gamma_N}(\boldsymbol{v}^t). 
    \label{eqn:lagrangeParametric}
\end{aligned}
\end{equation}
Notice, that $a_{\Omega_t}(\boldsymbol{w}^t, \boldsymbol{v}^t)$ can  be computed on $\Omega$ using a change of variable and the Jacobian $\psi(t)$. 
Analogously to previous results in level set topology optimization, the direct and adjoint equation are obtained through the optimality conditions of the Lagrangian (see \cite{allaire_level-set_2005, allaire_conception_2007, laurain_distributed_2015, laurain_level_2017} for detailed derivations of the optimality conditions.). 
Differentiation with respect to $\boldsymbol{v}^t$ yields the \emph{direct equation}, here the linear elasticity in equation~\ref{eq:linElasPerturb}. The \emph{adjoint equation} is obtained by derivation with respect to $\boldsymbol{w}^t$.\\
By evaluating said derivative at time zero and $\boldsymbol{w}^t=\boldsymbol{u}^t$,  the solution of $\boldsymbol{p}^0$ can be rewritten in terms of the direct problem's solution as $\boldsymbol{p}^0 = -2 \boldsymbol{u}^0$. \\

The shape derivative is obtained through a time derivative after these simplifications:

\begin{equation}
    d\mathcal{J}_C(\Omega)[{\boldsymbol{\theta}}]  = \frac{\diff}{\diff t} \mathcal{L}(t, \boldsymbol{w}^t, \boldsymbol{v}^t)\bigg|_{t=0} = \frac{\partial}{\partial t} \mathcal{L}(0, \boldsymbol{u}^0, \boldsymbol{p}^0) 
\end{equation}

To differentiate the Lagrangian (\ref{eqn:lagrangeParametric}) partially with respect to time it is rewritten in parametric space accounting for the appropriate mappings and metrics. 
With the derivation from Appendix~\ref{sec:technicalLemma} and for $\fu\in \mathcal{P}_{KL}$, the time derivative can be expressed as:
{\footnotesize
\begin{equation}
\begin{aligned}
    d\mathcal{J}_C(\Omega)[{\boldsymbol{\theta}}] = \quad  d &\int_{\Omega} \left( -\varepsilon(\boldsymbol{u}) :\primed{\ffC}: \varepsilon(\boldsymbol{u})- 2 \varepsilon(\boldsymbol{u}) :\ffC:\primed{\varepsilon}(\boldsymbol{u})  -  \varepsilon(\boldsymbol{u}) :\ffC :\varepsilon(\boldsymbol{u}) \operatorname{div}(\boldsymbol{\theta})   \right) \\  +  
     \frac{d^3}{12} & \int_\Omega  \left( -\kappa(\boldsymbol{u}) : \primed{\ffC}:\kappa(\boldsymbol{u}) -2 \kappa(\boldsymbol{u}) :\ffC:\primed{\kappa}(\boldsymbol{u})  - \kappa(\boldsymbol{u}) :\ffC:\kappa(\boldsymbol{u}) \operatorname{div}(\boldsymbol{\theta}) \right),
     \label{eqn:totShapeDeriv}
\end{aligned}
\end{equation}
}
with $\primed{\varepsilon}$, $\primed{\kappa}$ and $\primed{\ffC}$ given by equations (\ref{eqn:dotepsilon}), (\ref{eqn:dotkappa}) and (\ref{eqn:cdot}).

The time derivatives of (i) $\varepsilon^t$, (ii) $\kappa^t$ and (iii) $C^t$  are expressed as follows:
\begin{enumerate}[label=(\roman*)]
    \item Membrane strain time derivative $(\primed{{\varepsilon}})_{\alpha\beta}(\boldsymbol{u}) := \left(\frac{\partial {\varepsilon}^t({\boldsymbol{u}}^t)}{\partial t}\right)\bigg|_{t=0}$:
\begin{equation}
\begin{aligned}
    (\primed{\varepsilon})_{\alpha\beta}(\boldsymbol{u}) &=&& \frac{1}{2}\bigg( \primed{{\mathbf{a}}}_\alpha \cdot {\boldsymbol{u}}_{,\beta} + {\mathbf{a}}_\alpha \cdot \primed{{\boldsymbol{u}}}_{,\beta}+\primed{{\mathbf{a}}}_\beta \cdot {\boldsymbol{u}}_{,\alpha} + {\mathbf{a}}_\beta \cdot \primed{{\boldsymbol{u}}}_{,\alpha} \bigg).
\end{aligned}
\label{eqn:dotepsilon}
\end{equation}    
    \item Bending strain time derivative $(\primed{{\kappa}})_{\alpha\beta}(\boldsymbol{u}) := \left(\frac{\partial {\kappa}^t({\boldsymbol{u}}^t)}{\partial t}\right)\bigg|_{t=0}$:
\begin{equation}
\begin{aligned}
    (\primed{\kappa})_{\alpha\beta}(\boldsymbol{u}) &= && -\primed{{\mathbf{a}}}_{3} \cdot\left({\boldsymbol{u}}_{,\alpha,\beta} - \Gamma_{\alpha\beta}^\lambda {\boldsymbol{u}}_{,\lambda} \right) \\
    & && -{\mathbf{a}}_{3} \cdot \left(\primed{\boldsymbol{u}}_{,\alpha,\beta} - \primed{\Gamma}^\lambda_{\alpha\beta} {\boldsymbol{u}}_{,\lambda} -  \Gamma_{\alpha\beta}^\lambda \primed{\boldsymbol{u}}_{,\lambda} \right),
\end{aligned}
\label{eqn:dotkappa}
\end{equation}

where

\begin{equation}
\begin{aligned}
    \primed{\mathbf{a}}_\alpha &= \mathbf{a}_{\alpha, \mu} \vartheta^\mu \\
    \primed{\boldsymbol{u}}_{,\alpha} &= -\boldsymbol{u}_{, \lambda} \vartheta^\lambda_{,\alpha} \\
    \primed{\mathbf{a}}_3 &= \frac{\mathbf{a}_{1,\lambda } \vartheta^{\lambda } \times \mathbf{a}_2   +  \mathbf{a}_1 \times \mathbf{a}_{2,\lambda } \vartheta^{\lambda } }{|\mathbf{a}_1 \times \mathbf{a}_2|} - \frac{1}{2}\left(G^{\mu\gamma} G_{\gamma\mu,\lambda} \vartheta^\lambda\right)\mathbf{a}_3 \\
    \primed{\boldsymbol{u}}_{,\alpha,\beta} &= - \boldsymbol{u}_{,\gamma\alpha} \vartheta^\gamma_{,\beta}- \boldsymbol{u}_{,\gamma\beta} \vartheta^\gamma_{,\alpha} - \boldsymbol{u}_{,\gamma} \vartheta^\gamma_{,\alpha\beta}\\
    \primed{\Gamma}^\lambda_{\alpha\beta} &= {\Gamma}^\lambda_{\alpha\beta, \mu} \vartheta^\mu
\end{aligned}
\end{equation}
    \item Elasticity tensor time derivative $(\primed{{\rC}})^{\alpha\beta\lambda\mu} := \left(\frac{\partial {\rC}^t}{\partial t}\right)\Big|_{t=0}$:
{\footnotesize
\begin{equation}
\begin{aligned}
(\primed{{\rC}})^{\alpha\beta\lambda\mu} = \frac{E}{2(1+ \nu)} {\vartheta}^\gamma \bigg( &\rG^{\alpha\lambda}{}_{,\gamma} G^{\beta\mu} +\rG^{\alpha\lambda} G^{\beta\mu}{}_{,\gamma} +  \rG^{\alpha\mu}{}_{,\gamma}G^{\beta\lambda} +
\rG^{\alpha\mu}G^{\beta\lambda}{}_{,\gamma} 
\\ +& \frac{2\nu}{1-\nu} \left( \rG^{\alpha\beta}{}_{,\gamma} G^{\lambda\mu} +\rG^{\alpha\beta} G^{\lambda\mu}{}_{,\gamma} \right)\bigg).
\end{aligned}
\label{eqn:cdot}
\end{equation}
}
\end{enumerate}

\begingroup
\paragraph{Optimization with a level set approach}
To solve the optimization problem posed in Problem \ref{eq:optProb}, a level set based method is applied using isogeometric analysis. 
The material domain $\Omega$ is represented implicitly by the level set, the shape derivative is constructed via the speed method within a Lagrangian formulation of the constrained cost functional, and a regularized velocity field enables the integration of the Hamilton--Jacobi equation to evolve the level set and implicitly the domain shape.
\goodbreak
The material domain $\Omega$ of the mid-surface of the shell is defined by a level set function $\phi$:  
\begin{align}
    \phi:\quad {\mathcal{S}} \rightarrow \mathbb{R}, 
&&  
    \phi(\boldsymbol{x}) = 
    \begin{cases}
            <0: \quad \fx \in \Omega\\
            >0: \quad \fx \in \mathcal{S}/\Omega\\
            =0: \quad \fx \in \partial \Omega
    \end{cases}
\end{align}
For the void regions $\mathcal{S} \setminus \Omega$, a weak ersatz material is introduced to prevent singularities in the mechanical fields. This leads to the following distribution of the elastic moduli:
\begin{align}
\mathbb{C}(\phi) =
\left\{
\begin{array}{lllll}
\mathbb{C} & \text{in } \Omega & \qquad & \text{material}, & \phi < 0, \\[6pt]
\epsilon  \mathbb{C} & \text{in } \mathcal{S} \setminus \Omega & \qquad & \text{ersatz-material, } & \phi > 0, 
\end{array}
\right.
\label{eqn:ymNoInterp}
\end{align}
where $0<\epsilon \ll 1$ is a small ersatz parameter, chosen as $\epsilon = 0.001$ in the numerical examples.

\paragraph{Regularized velocity field}
The preceding computations of $d\mathcal{J}_C(\Omega)[\mathbf{\theta}]$ provide the variation of the cost for a given velocity field. To compute a descent direction, we conduct a biharmonic regularization of the shape derivative on $\Omega$ in analogy to \cite{laurain_distributed_2015, laurain_level_2017}.

As the kinematic definition of the bending strains $\kappa(\boldsymbol{u})$ involves second-order derivatives of the velocity field $\boldsymbol{\theta}$ (see Equation~\ref{eqn:dotkappa}), standard $H^1$-regularization is insufficient. Given this higher regularity requirement, we apply a biharmonic regularization to find the descent direction $\boldsymbol{\theta}$ which leads to the following regularization problem:
\begin{problem}[Regularization Problem: $\mathcal{P}_{R}(\mathcal{J})$]
\label{prob:regProb}
Find a descent direction $\boldsymbol{\theta}$ such that:
\begin{align}
    \underbrace{\int_{\mathcal{S}} 
    \left(
    \alpha_B \Delta \boldsymbol{\upsilon} \cdot \Delta \boldsymbol{\theta} + \alpha_H \nabla \boldsymbol{\upsilon} \cdot \nabla \boldsymbol{\theta} + \boldsymbol{\upsilon} \cdot \boldsymbol{\theta} 
     \right) }_{b(\boldsymbol{\theta}, \ \boldsymbol{\upsilon} )}
    &= -d\mathcal{J}(\Omega)[{\boldsymbol{\upsilon}}] \qquad && \forall \boldsymbol{\upsilon}, 
        \label{eqn:regProb}
    \\
    \boldsymbol{\theta}\cdot \boldsymbol{n}&=0 \qquad &&\text{on }\partial \mathcal{S} \nonumber
\end{align}
where $\alpha_B=h_{\text{elem}}^4>0$ and $\alpha_H=2h_{\text{elem}}^2>0$ are discretization-dependent parameters, chosen numerically based on a factor of the element size and $\Delta$ corresponds to the Laplace-Beltrami operator.

\end{problem}

\paragraph{Hamilton--Jacobi equation to evolve the level set}
The material domain is updated at each optimization iteration by evolving the level set function $\phi$ over a pseudo-time variable $t$. Following standard approaches \cite{allaire_level-set_2005, hubner_scherer_topology_2024}, this is governed by the Hamilton--Jacobi advection equation which can be written as the following problem:
\begin{problem}[Advection Problem on level set: $\mathcal{P}_{A}(\phi_0, \boldsymbol{\theta})$]
    Find the evolution of the level set function as the evolution of the level set variable over a pseudo-time t $\phi(\boldsymbol{x}, t)$ such that:
    \begin{align}
        \frac{\partial \phi}{\partial t} + \boldsymbol{\theta} \cdot \nabla\phi = 0, \qquad &\text{ in }  \mathcal{S} 
        \label{eqn:HamiltonJacobi}
        \\
        \text{ with } \phi(\boldsymbol{x}, 0) = \phi_0(\boldsymbol{x}) \qquad &\text{ on }  \mathcal{S}, \nonumber
    \end{align}
    where $\boldsymbol{\theta}$ denotes the regularized velocity field solution of $\mathcal{P}_R$ and $\phi_0$. 
\end{problem}
During the advection governed by Equation~\eqref{eqn:HamiltonJacobi}, the level set function $\phi$ typically ceases to remain a signed distance function ($|\nabla \phi| = 1$), leading to steep or flat gradients that can cause numerical instabilities. 
To preserve the signed distance property and ensure accurate interface identification, the level set is reinitialized every five iterations. 
\endgroup

\subsection{Numerical Optimization Algorithm}
\label{sec:numericalOptim}
\paragraph{Implementation with Isogeometric Analysis}
The developed variational problems are solved with Isogeometric Analysis (IGA), where both the geometry and the shape functions are encoded by NURBS functions. 
As more complex surfaces $\mathcal{S}$ often require the division into multiple patches, continuity needs to be ensured. 
As it is standard in multi-patch IGA, this division can be written as:
\begin{align}
\overline{\mathcal{S}} &= \bigcup_{k=1}^K \overline{\mathcal{D}}_k \qquad \text{ with } \mathcal{D}_i \cap \mathcal{D}_j = \emptyset \text{ when } i \ne j, \\
    \gamma_{kl} &= \partial \mathcal{D}_k \cap \partial \mathcal{D}_l \quad k<l; \quad k, l \in \mathcal{I}.
\end{align}
Hence, the surface $\mathcal{S}$ is decomposed into the different open patches $\mathcal{D}_k$ and neighboring patches share a common interface $\gamma_{kl}$. \\
In the multi-patch case, the continuity of the fields $\fu, \ \boldsymbol{\theta}$ has to be enforced at the interfaces. 
The preceding problems $\mathcal{P}_{KL}$, $\mathcal{P}_{R}$ are thus modified with additional penalty terms to ensure said continuity. 
For the \textit{linear elasticity problem} ($\mathcal{P}_{KL}$), the $C^1$-continuity of the displacement field $\boldsymbol{u}$ is enforced by a penalty approach based on \cite{herrema_penalty_2019, proserpio_penalty_2022} with the implementation from \cite{zhao_open-source_2022}. \\ 
 The \textbf{regularization problem} ($\mathcal{P}_{R}$)  enforces $C^1$-continuity of the descent direction $\boldsymbol{\theta}$ by a symmetric interior penalty Galerkin method, based on \cite{epshteyn_estimation_2007, moore_discontinuous_2018, moore_multipatch_2020}.
The modified problems for $N$ patches can be written as:
\begin{align}
    \sum_{k=1}^N a_{\Omega, k} (\fu_k,\ \fv_k) + \sum_{\gamma_{kl}} p^{KL}(\fu_k,\ \fv_k) &= \sum_{k=1}^NL_{\Gamma_N}(\fv) \\
    \sum_{k=1}^N b_k (\boldsymbol{\theta}_k, \ \boldsymbol{\upsilon}_k) + \sum_{\gamma_{kl}} p^R(\boldsymbol{\theta}_k, \ \boldsymbol{\upsilon}_k) &= -\sum_{k=1}^N dJ(\Omega_k), 
\end{align}
where $a_{\Omega, k}$ and $b_k$ correspond to the patch wise bilinear forms from Equations \ref{eqn:linElasKL} and \ref{eqn:regProb} and $p^{KL}, \ p^{R}$ correspond to the penalty coupling terms, whose expression is given in Equations~\ref{eq:penKL}~\&~\ref{eq:penReg}.
{\footnotesize
\begin{align}
        p^{KL} &= \sum_{(k,l)\in \mathcal{I}}\int_{\gamma_{kl}} \alpha_D \left(\boldsymbol{u}^{(k)}-\boldsymbol{u}^{(l)}\right)\cdot\left(\boldsymbol{v}^{(k)}-\boldsymbol{v}^{(l)}\right) +\sum_{(k,l)\in \mathcal{I}} \int_{\gamma_{kl}} \alpha_R \biggl[  \nonumber         \\  
    &\left( \boldsymbol{a}_3^{(k)}(\boldsymbol{u}^{(k)})\cdot \boldsymbol{a}_3^{(l)}(\boldsymbol{u}^{(l)}) -\boldsymbol{A}_3^{(k)}\cdot \boldsymbol{A}_3^{(l)} \right)\left( \boldsymbol{a}_3^{(k)}(\boldsymbol{v}^{(k)})\cdot \boldsymbol{a}_3^{(l)}(\boldsymbol{u}^{(l)})+\boldsymbol{a}_3^{(k)}(\boldsymbol{u}^{(k)})\cdot \boldsymbol{a}_3^{(l)}(\boldsymbol{v}^{(l)})\right) \nonumber \\
+&\left( \boldsymbol{a}_N^{(k)}(\boldsymbol{u}^{(k)})\cdot \boldsymbol{a}_3^{(l)}(\boldsymbol{u}^{(l)}) -\boldsymbol{A}_N^{(k)}\cdot \boldsymbol{A}_3^{(l)} \right)\left( \boldsymbol{a}_N^{(k)}(\boldsymbol{v}^{(k)})\cdot \boldsymbol{a}_3^{(l)}(\boldsymbol{u}^{(l)}) +\boldsymbol{a}_N^{(k)}(\boldsymbol{u}^{(k)})\cdot \boldsymbol{a}_3^{(l)}(\boldsymbol{v}^{(l)}) \right) \biggr], 
\label{eq:penKL}
\end{align}
}
where $\boldsymbol{A}_N$ and $\boldsymbol{a}_N$ corresponds to the interface normal. 
The corresponding penalty parameters $\alpha_D, \ \alpha_R$ are chosen in accordance to \cite{herrema_penalty_2019} as 
\begin{align*}
    \alpha_D = \alpha \frac{E d}{h(1-\nu^2)} && \alpha_R = \alpha \frac{E d^3}{12 h(1-\nu^2)}.
\end{align*}

\begin{equation}
\begin{split}
p^R = \beta \Bigg[ 
            - & \sum_{(k,l)\in\mathcal{I}} \int_{\gamma_{kl}}[\![ \nabla \boldsymbol{\theta} ]\!] \{\Delta \boldsymbol{\upsilon} \}-\sum_{(k,l)\in\mathcal{I}} \int_{\gamma_{kl}}[\![ \nabla \boldsymbol{\upsilon} ]\!] \{\Delta \boldsymbol{\theta} \}\\
            + & \sum_{(k,l)\in\mathcal{I}} \int_{\gamma_{kl}}[\![ \boldsymbol{\theta} ]\!] \{\nabla \Delta \boldsymbol{\upsilon} \}+\sum_{(k,l)\in\mathcal{I}} \int_{\gamma_{kl}}[\![ \boldsymbol{\upsilon} ]\!] \{\nabla\Delta \boldsymbol{\theta} \} \\
+ &\sum_{(k,l)\in\mathcal{I}} \int_{\gamma_{kl}}\frac{\delta_1}{h_{kl}^3}[\![ \boldsymbol{\upsilon} ]\!] [\![ \boldsymbol{\theta} ]\!]+\sum_{(k,l)\in\mathcal{I}} \int_{\gamma_{kl}} \frac{\delta_0}{h_{kl}} [\![ \nabla \boldsymbol{\upsilon} ]\!] \{\nabla \boldsymbol{\theta} \}
\Bigg]\\
            +&\alpha \Bigg[  \sum_{(k,l)\in\mathcal{I}} \int_{\gamma_{kl}}\frac{\delta_1}{h_{kl}}[\![ \boldsymbol{\upsilon} ]\!] [\![ \boldsymbol{\theta} ]\!] \\
            -& \sum_{(k,l)\in\mathcal{I}} \int_{\gamma_{kl}}[\![ \boldsymbol{\theta} ]\!] \{\nabla \boldsymbol{\upsilon} \}-\sum_{(k,l)\in\mathcal{I}} \int_{\gamma_{kl}}[\![ \boldsymbol{\upsilon} ]\!] \{\nabla \boldsymbol{\theta} \}
\Bigg]\\
    +& \sum_{k=1}^K \int_{\mathcal{D}_k}  \boldsymbol{\upsilon} \cdot \boldsymbol{\theta},
\end{split}
\label{eq:penReg}
\end{equation}
where $\beta, \ \alpha$ correspond to the prefactors of the biharmonic and harmonic components of the regularization respectively and $[\![ \cdot ]\!] $ and $\{\cdot \}$ correspond to the jump and average operator along the interface.
The overall optimization is conducted as a gradient descent method, the structure of which is detailed in \autoref{alg:robOptimAlg}. 
\goodbreak
The isogeometric implementation is based on the \textsc{FEniCS} project \cite{logg_automated_2012} with the IGA extension via Bezier extraction using the \textsc{tIGAr} package \cite{kamensky_tigar_2019} and the non-matching KL shells coupling from the \textsc{PENGoLINS} package \cite{zhao_open-source_2022}. 
Furthermore, the numerical results and the associated codes are made publicly available in a  repository.

The optimization is conducted as a gradient method. 
Its structure is shown in \autoref{alg:robOptimAlg}.

\begin{minipage}{\linewidth}
\begin{algorithm}[H]
    \small
    \KwData{Multi-patch geometry $\{ \mathcal{D}_k \}$, Initial level set $\{\phi_k\}^{(0)}$, Volume penalization $\Lambda$}
    \KwIn{Initial timestep $\tau > 0$, Increase factor $\beta > 1$, Decrease factor $\gamma \in (0,1)$, lsMax $l>0$, min timestep $\tau^{min}$, Max Iterations $N^{\mathrm{max}}$}
    \BlankLine
    \textcolor{blue}{Solve elasticity problem $\mathcal{P}_{KL}$ using $\phi^{(0)}$ as $\fu$ \;}
    \textcolor{Plum}{Solve samples $\omega_i$ from linear elasticity problem $\mathcal{P}_{KL}$ using $\phi^{(0)}$ as ${\fu_\omega}$ \;}
    Compute initial cost $\mathcal{J}^{(0)}$\;
    $m\gets 0$, $convergence \gets \textbf{false}$\;
    \While{$m < N^{\mathrm{max}}$\textbf{ and not } convergence}{
        \textcolor{blue}{Compute shape derivative $d\mathcal{J}(\Omega)$ using (\ref{eqn:totShapeDeriv}) and current state $\fu$ \;}
        \textcolor{Plum}{Compute $d\mathcal{J}_R(\Omega)$ using (\ref{eqn:totRobustShapeDeriv}) and samples from current state ${\fu_\omega}$  \;}
        Find $\theta^m$ by regularization problem (\ref{prob:regProb}) $\mathcal{P}_{R}$ with $d\mathcal{J}(\Omega)$ or $d\mathcal{J}_R(\Omega)$\;   
        $accepted \gets \textbf{false}$, $lsIt \gets 0$\;
        \While{\upshape \textbf{not} $accepted$ \textbf{and} $lsIt < l$ }{
            Compute trial $\phi^{trial}$ via evolution eq. \ref{eqn:HamiltonJacobi} using $\boldsymbol{\theta}^{(m)}$, $\phi^{(m)}$ and $\tau$\;
            \textcolor{blue}{Solve linear elasticity $\mathcal{P}_{KL}$ using $\phi^{trial}$ as  $\fu$ \;}
            \textcolor{Plum}{Solve samples $\omega_i$ from $\mathcal{P}_{KL}$ using $\phi^{trial}$ as ${\fu_\omega}$ \;}
            Compute cost $\mathcal{J}^{trial}$\;
            \eIf{$\mathcal{J}^{trial} < \mathcal{J}^{(m)}$}{
                \textbf{Accept step}\;
                $\phi^{(m+1)} \gets \phi^{trial}$, \ $\mathcal{J}^{(m+1)} \gets \mathcal{J}^{trial}$, $accepted \gets \textbf{true}$\;
                \If{$lsIt=0$}{
                    $\tau \gets \tau \cdot \beta$ \tcc*[r]{Accepted: Increase time step}
                }
            }{
                \textbf{Reject step}\;
                $\tau \gets \tau \cdot \gamma$ \tcc*[r]{Cost increased: Red. time step}
            }
            $lsIt \gets lsIt + 1$\;
            \If{$\tau < \tau^{min}$}{
                $convergence \gets \textbf{true}$\;
            }
        }
        $m\gets m+1$\;
    }
    \caption{Gradient based topology optimization with line search}
    \label{alg:robOptimAlg}
\end{algorithm}
\small Alg. \ref{alg:robOptimAlg} describes the optimization method. 
Lines in blue refer to problems that need to be solved for the deterministic problem while the purple lines refer to the extension under uncertainty that will be explained in the following section. 
\end{minipage}
Let us further specify that $\boldsymbol{u}, \ \boldsymbol{\theta}$ are written as NURBS  functions of degree 3 while $\phi$ is written as a B-Spline function with the same discretization and degree. 

\subsection{Extension of the Shape Derivative to Loading and Material Uncertainties}
\paragraph*{Robust Compliance with Standard Deviation $\sigma$}
The preceding shape optimization framework will be extended next for problems exhibiting material or loading uncertainties. Let us denote by $\omega$ a parameter sampling a distribution of the elastic moduli $\mathbb{C}(\omega)$ or the external loadings $g(\omega)$. Moreover, these spatial distributions of uncertainties are modeled by Gaussian random fields (GRFs) with Matérn covariance, which are widely popular for various engineering applications \cite{chu_stochastic_2021,koh_stochastic_2023} and are generated through the variational problem described in Appendix~\ref{sec:grfWithSPDE}. For further details, the reader is referred to \cite{lindgren_explicit_2011,lindgren_spde_2022,khristenko_analysis_2019,croci_efficient_2018}. \\
In such a case, each sample $\omega$ generates a mechanical problem and a given compliance measured by the cost $\mathcal{J}_C(\omega)$. The use of the \emph{robust compliance} is motivated by previous work that can for example be found in \cite{chen_level_2010,conti_risk_2011,conti_stochastic_2018,martinez-frutos_robust_2016,martinez-frutos_risk-averse_2021}.

For a given distribution of $\omega$, one can define the expectancy $\mathbb{E}_\omega[\mathcal{J}_C(\omega)]$ and the variance $\mathbb{V}_\omega[\mathcal{J}_C(\omega)]$ of the compliance. The robust compliance cost functional is then defined as follows:
    \begin{equation}
        \mathcal{J}_R = \mathbb{E}_\omega [\mathcal{J}_C(\omega)]+\alpha \sqrt{\mathbb{V}_\omega[\mathcal{J}_C(\omega)]}, 
    \end{equation}
where $\mathcal{J}_C(\omega)$ corresponds to the compliance function of an individual material sample $\omega$ and $\alpha$ denotes a robustness parameter $\alpha\geq 0$. 
Analogous with the preceding deterministic shape optimization problem, the derivative of the cost functional will be constructed with the help of a Lagrangian function. Since the standard deviation of a random variable $X$ can be expanded as $\sqrt{\mathbb{V}_\omega[X]} = \sqrt{\mathbb{E}_\omega[X^2] - \mathbb{E}_\omega[X]^2}$, the Lagrangian $\mathcal{L}_R$ for the robust problem is written in analogy to the deterministic Lagrangian from Equation~\ref{eqn:lagrangeParametric} as:
\begin{equation}
    \begin{aligned}
        \mathcal{L}_R(t, \{\boldsymbol{w}_\omega^t\}, \{\boldsymbol{v}_\omega^t\}) =&  \mathbb{E}_\omega [ \mathcal{J}_C(\boldsymbol{w}_\omega^t) ] + \alpha \sqrt{\mathbb{V}_\omega [ \mathcal{J}_C(\boldsymbol{w}_\omega^t) ]} \\ 
        +& \mathbb{E}_\omega \left[ a_{\Omega_t}(\boldsymbol{w}_\omega^t,\ \boldsymbol{v}_\omega^t)-L_{\Gamma_N}(\boldsymbol{v}_\omega^t)\right], 
        \label{eqn:lagrangeRobust}
    \end{aligned}
\end{equation}

where $\mathcal{J}_C$ corresponds to the compliance as defined in Equation \ref{eqn:defCompCost}. 

Analogous to the deterministic framework, the direct and adjoint equations are obtained through the optimality conditions of the Lagrangian. 

To simplify the resulting expressions, we define a realization-dependent scalar weight $r(\omega)$:
\begin{equation}
    r(\omega) := 1 + \frac{\alpha}{\sqrt{\mathbb{V}_\omega[J_C]}}(J_C(\omega) - \mathbb{E}_\omega[J_C]).
    \label{eqn:weightH}
\end{equation}

By evaluating said derivative at time zero and substituting $\boldsymbol{w}_\omega^t=\boldsymbol{u}_\omega^t$, the solution of $\boldsymbol{v}_\omega^0$ can through the self-adjointness of the problem be rewritten in terms of the direct problem's solution:
\begin{equation}
    \boldsymbol{v}_\omega^0 = -2 r(\omega) \boldsymbol{u}_\omega^0.
    \label{eqn:adjointRobust}
\end{equation}
The shape derivative is obtained through a time derivative of the Lagrangian at $t=0$ after substituting the adjoint variable as shown in Equation \ref{eqn:adjointRobust}. 
By differentiating the Lagrangian partially with respect to time, factoring out $r(\omega)$ and further simplifications, the time derivatives of the strain and elasticity tensors ($\primed{\varepsilon}$, $\primed{\kappa}$, and $\primed{\ffC}$) are substituted analogously to the deterministic expressions. 

Assembling the different terms, we obtain the shape derivative of the robust optimization problem as the expectation of the scaled deterministic derivatives:
\begin{align}
    d\mathcal{J}_R(\Omega)[{\boldsymbol{\theta}}]  = \quad  &\mathbb{E}_\omega \Bigg[ r(\omega) \Bigg( d \int_{\Omega} \Big( -\varepsilon(\boldsymbol{u}_\omega) :\primed{\ffC}_\omega : \varepsilon(\boldsymbol{u}_\omega) - 2 \varepsilon(\boldsymbol{u}_\omega) : \ffC_\omega:\primed{\varepsilon}(\boldsymbol{u}_\omega) \nonumber \\
     & \qquad -  \varepsilon(\boldsymbol{u}_\omega) : \ffC_\omega:\varepsilon(\boldsymbol{u}_\omega) \operatorname{div}(\boldsymbol{\theta})   \Big) \nonumber\\  
    + &\frac{d^3}{12} \int_\Omega \Big( -\kappa(\boldsymbol{u}_\omega) :\primed{\ffC}_\omega : \kappa(\boldsymbol{u}_\omega)-2 \kappa(\boldsymbol{u}_\omega) : \ffC_\omega:\primed{\kappa}(\boldsymbol{u}_\omega) \nonumber\\
    & \qquad - \kappa(\boldsymbol{u}_\omega) : \ffC_\omega:\kappa(\boldsymbol{u}_\omega) \operatorname{div}(\boldsymbol{\theta}) \Big) \Bigg) \Bigg].
    \label{eqn:totRobustShapeDeriv}
\end{align}

With this developed shape derivative, the same approach as in Section~\ref{sec:numericalOptim} can be applied with the regularization problem $\mathcal{P}_R$ and the following advection problem $\mathcal{P}_A$, however, with the shape derivative being replaced by Equation~\ref{eqn:totRobustShapeDeriv}.\\
The process is illustrated in Algorithm~\ref{alg:robOptimAlg}. 

\section{Results and Discussion}
\label{sec:numRes}
The methodology is applied to a variety of benchmark problems to demonstrate the effectiveness of the developed approach and the use of IGA KL shells for topology optimization. 

In addition to examples with pure deterministic loading and materials, examples including stochastic deviations are included. 
The examples are compared to common benchmark examples, where available. 
For further examples, the reader is referred to the GitHub repository. 
\\ 
Unless explicitly stated, all physical dimensions are given in coherent SI units.

\begin{example}[Multi-Patch Hemisphere]
    \label{ex:hemisphere_det}

    To test the approach on curved multi-patch domains, with non-uniform element size in the physical domain, the optimization is tested on a hemisphere with a central cutout. 
    The problem setting is adapted from \cite{macneal_proposed_1985}, as part of the shell obstacle course. 
    Half a sphere of radius 10 with a hole of $18^\circ$ is modeled with four shell patches of thickness $0.02$, as depicted in \autoref{fig:hemisphere_det_ini}, along with the boundary conditions of the mechanical problem and the distributed load of magnitude $6\cdot 10^{4}$ in negative $e_z$ applied to the structure. 
    The Young's modulus and Poisson ratio are chosen as $200\cdot 10^9$ and $0.3$ respectively with a volume Lagrange multiplier $\Lambda = 20$. 

    \begin{figure}[H]
        \centering
\includegraphics[width=\linewidth]{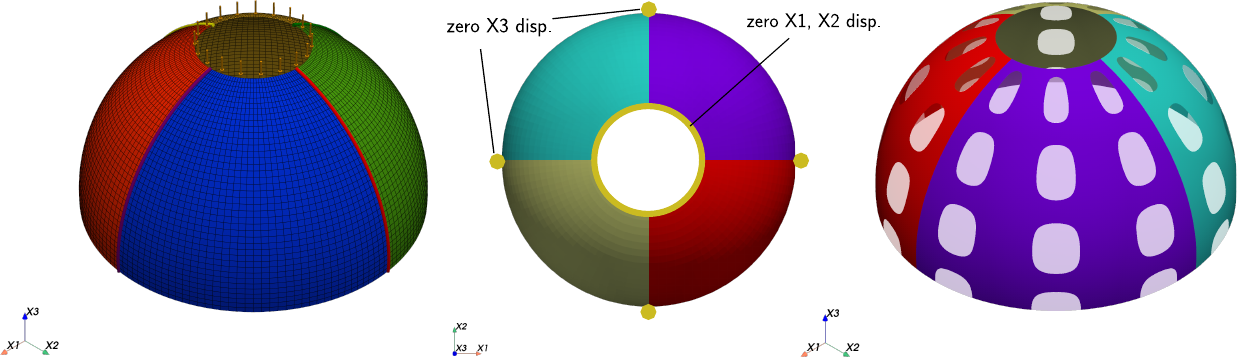}
        \caption{
        Multi-patch hemisphere benchmark. Geometry: Half-sphere of radius $R = 10\,m$ with a central hole of $18^\circ$, discretized into 4 NURBS patches (cubic basis, $40 \times 40$ elements per patch). Material: Young’s modulus $E = 200$, Poisson’s ratio $\nu = 0.3$, thickness $d = 0.02$. Distributed load $\boldsymbol{g} = -6 \times 10^4 \mathbf{e}_z$ (red arrows). Boundary conditions: Clamped at the hole (blue), free elsewhere.
}
        \label{fig:hemisphere_det_ini}
    \end{figure}
%
%
%
    \begin{figure}[H]
        \centering
        \includegraphics[width=.8\linewidth]{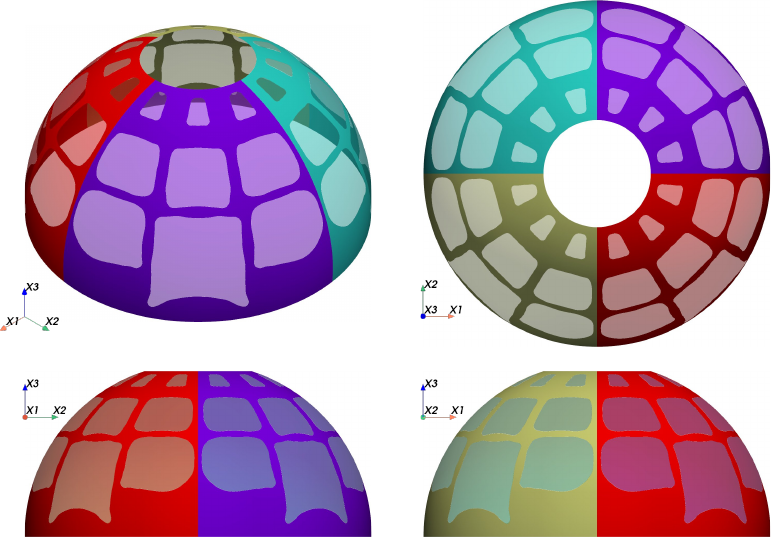}
        \caption{
Optimized topology for the multi-patch hemisphere 150 iterations with Rib-like supports (blue) connecting the clamped region to the loaded area. One can observe the symmetric reinforcement structure without explicit symmetry constraints.
For comparison, we achieved similar performance to \cite{dede_isogeometric_2012} ($128 \times 128 \times 1$ solid IGA elements) with $40 \times 40$ cubic NURBS KL elements.
}
        \label{fig:hemisphere_det_res}
    \end{figure}
%
%
%
     \begin{figure}[H]
        \centering
        \includegraphics[width=0.7\linewidth]{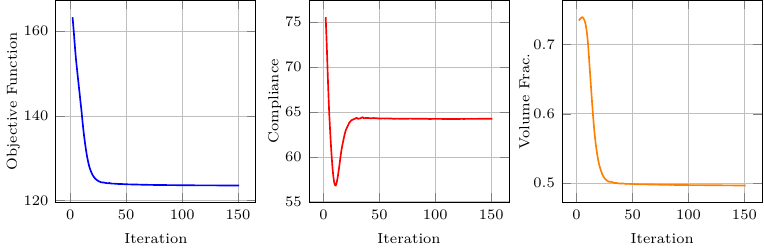}
        \caption{     
Convergence history for the multi-patch hemisphere optimization:
Objective function $\mathcal{J}(\Omega)$ (blue), compliance $\mathcal{J}_C$ (red), and volume fraction $V_0$ (yellow) vs. iteration number.
Settings: Initial timestep $\tau = 0.1$, line search parameters $\alpha = 1.2$,$\gamma = 0.5$. The objective function stabilizes after $\approx 30$ iterations.
}        
\label{fig:hemisphere_det_convergence}
    \end{figure}
The optimal shapes obtained are shown in \autoref{fig:hemisphere_det_res}. We remark that the method results in a relatively symmetric result despite not employing symmetry conditions. 
Similar shapes are obtained in \cite{dede_isogeometric_2012} on a slightly different geometry with solid elements with $128\times 128\times 1$ elements and degree 2, instead of $40\times 40$ elements  used here with cubic NURBS. The convergence of this process is displayed in \autoref{fig:hemisphere_det_convergence}.
\end{example}
\FloatBarrier
\begin{example}[Single-Patch Double Parabola]
    \label{ex:doublePara}
    \subex[Deterministic Point Loading]
    The next example demonstrates the effectiveness of the method on curved surfaces by optimizing a double parabola geometry. 
    A similar study of this geometry is conducted in \cite{cai_integrated_2023}, analyzing the effect of thickness and curvature on the optimal shapes. We study now the effect of the direction of the load $\ff$ on the optimal shape. Here, we also use cubic NURBS, but significantly fewer elements than in \cite{cai_integrated_2023} ($40^2$, instead of $64^2$) and we do not use additional integration points on the boundaries. 
    
    The geometry is described by a double parabola of side length of $5$ and height of $0.6$. The material has a Young's modulus $2\cdot 10^{11}$ and Poisson ratio $0.3$ is considered. Three load cases are investigated, with three distinct concentrated loads applied to the central position of the double parabola. The corresponding force vectors and Lagrange multipliers are (a) Vertical load: $\ff=[0,0,-5\cdot10^4]$ and $\Lambda=2$ (b) Horizontal load: $\ff=[5\cdot10^4,0,0]$ and $\Lambda=5$ (c) Diagonal load:  $\ff=[-\frac{10^4}{\sqrt{2}},-\frac{10^4}{\sqrt{2}},0]$  and $\Lambda = 5.0$\footnote{The volume Lagrange multiplier is chosen differently for the vertical load and the horizontal loads due to the difference in compliance for bending and membrane dominated load cases. The volume fractions would otherwise be too different for a reasonable comparison. }. 
    The Lagrange multipliers and forces differ due to the different stiffnesses in the different directions and planes. 
    The geometry, load, boundary conditions and initialization are shown in \autoref{fig:doublePara_geom}. 
    The optimal shapes obtained via the optimization are shown in Figure~\ref{fig:dp_res_loadCasesDet}.
%
%
%
\begin{figure}[H]
    \centering
\includegraphics[width=\linewidth]{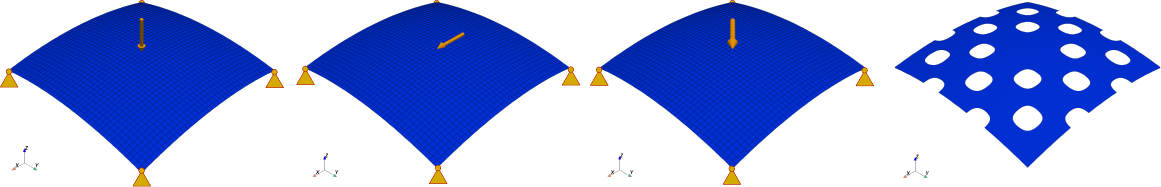}
    \begin{minipage}[t]{0.24\linewidth}
        \centering
        \vspace{2pt}
        \footnotesize Vertical Loading\\ $\boldsymbol{f}=[0,0,-5\cdot10^4]$.
    \end{minipage}\hfill 
    \begin{minipage}[t]{0.24\linewidth}
        \centering
        \vspace{2pt}
        \footnotesize Horizontal Loading\\ $\boldsymbol{f}=[5\cdot10^4, 0,0]$.
    \end{minipage}\hfill 
    \begin{minipage}[t]{0.24\linewidth}
        \centering
        \vspace{2pt}
        \footnotesize Diagonal Loading\\ $\boldsymbol{f}=[-\frac{10^4}{\sqrt{2}},-\frac{10^4}{\sqrt{2}},0]$.
    \end{minipage}\hfill 
    \begin{minipage}[t]{0.24\linewidth}
        \centering
        \vspace{2pt}
        \footnotesize Initialization
    \end{minipage}
    \caption{
Double parabola geometry and deterministic load cases. Geometry: Side length $L = 5$, height $H = 0.6$, cubic NURBS ($40 \times 40$ elements). Material: $E = 200 $,$\nu = 0.3$.
Load cases:
(a) vertical, 
(b) horizontal,
(c) diagonal and (d) initial configuration.
}
    \label{fig:doublePara_geom}
    \end{figure}
   
    \begin{figure}[H]
        \centering
        \begin{minipage}[t]{0.32\linewidth}
            \centering
\includegraphics[width=\linewidth]{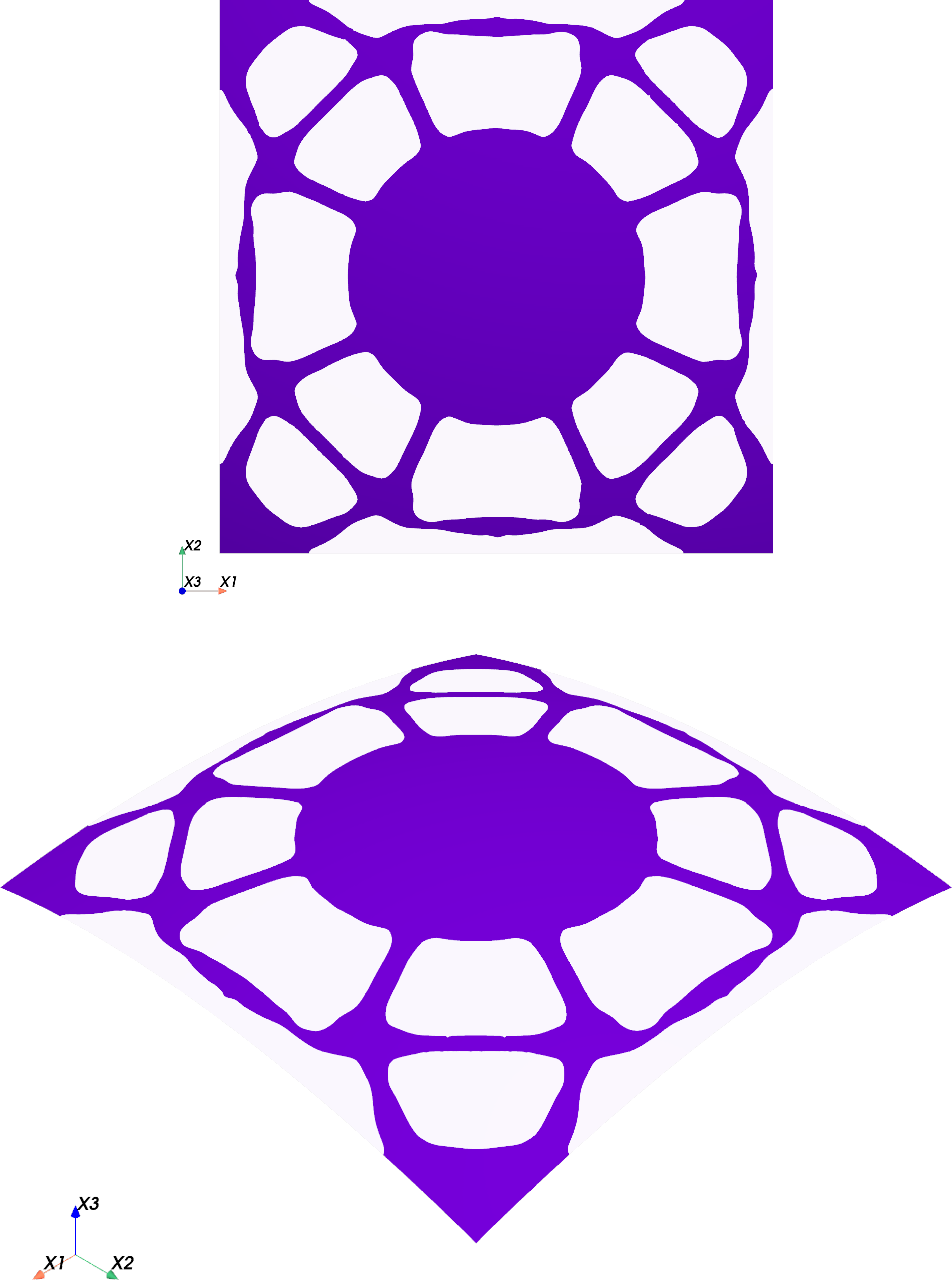}
            \vspace{2pt}
            \footnotesize Vertical Loading\\ $\boldsymbol{f}=[0,0,-5\cdot10^4]$.
        \end{minipage}\hfill 
        \begin{minipage}[t]{0.32\linewidth}
            \centering
\includegraphics[width=\linewidth]{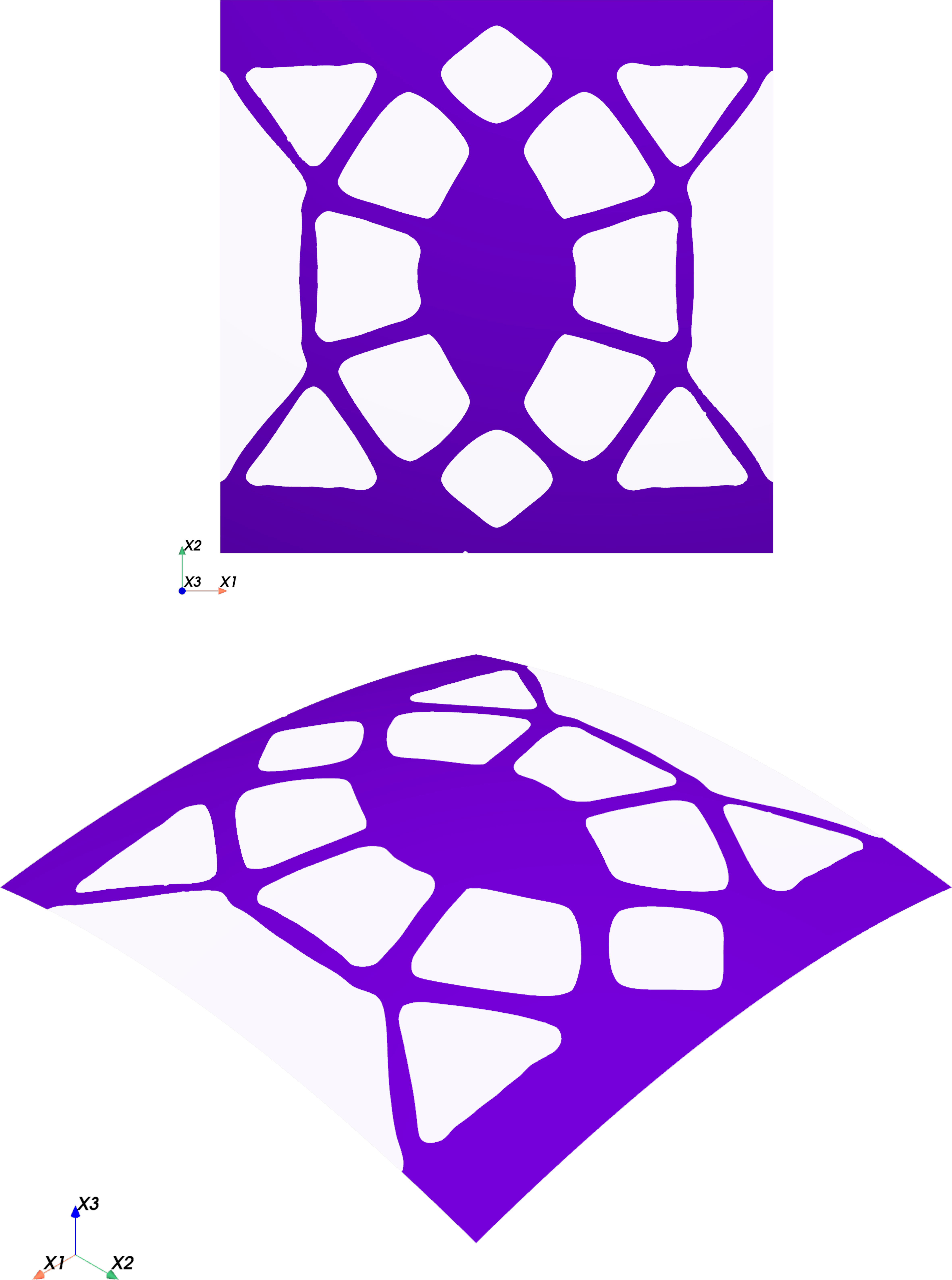}
            \vspace{2pt}
            \footnotesize Horizontal Loading\\ $\boldsymbol{f}=[5\cdot10^4, 0,0]$.
        \end{minipage}\hfill 
        \begin{minipage}[t]{0.32\linewidth}
            \centering
\includegraphics[width=\linewidth]{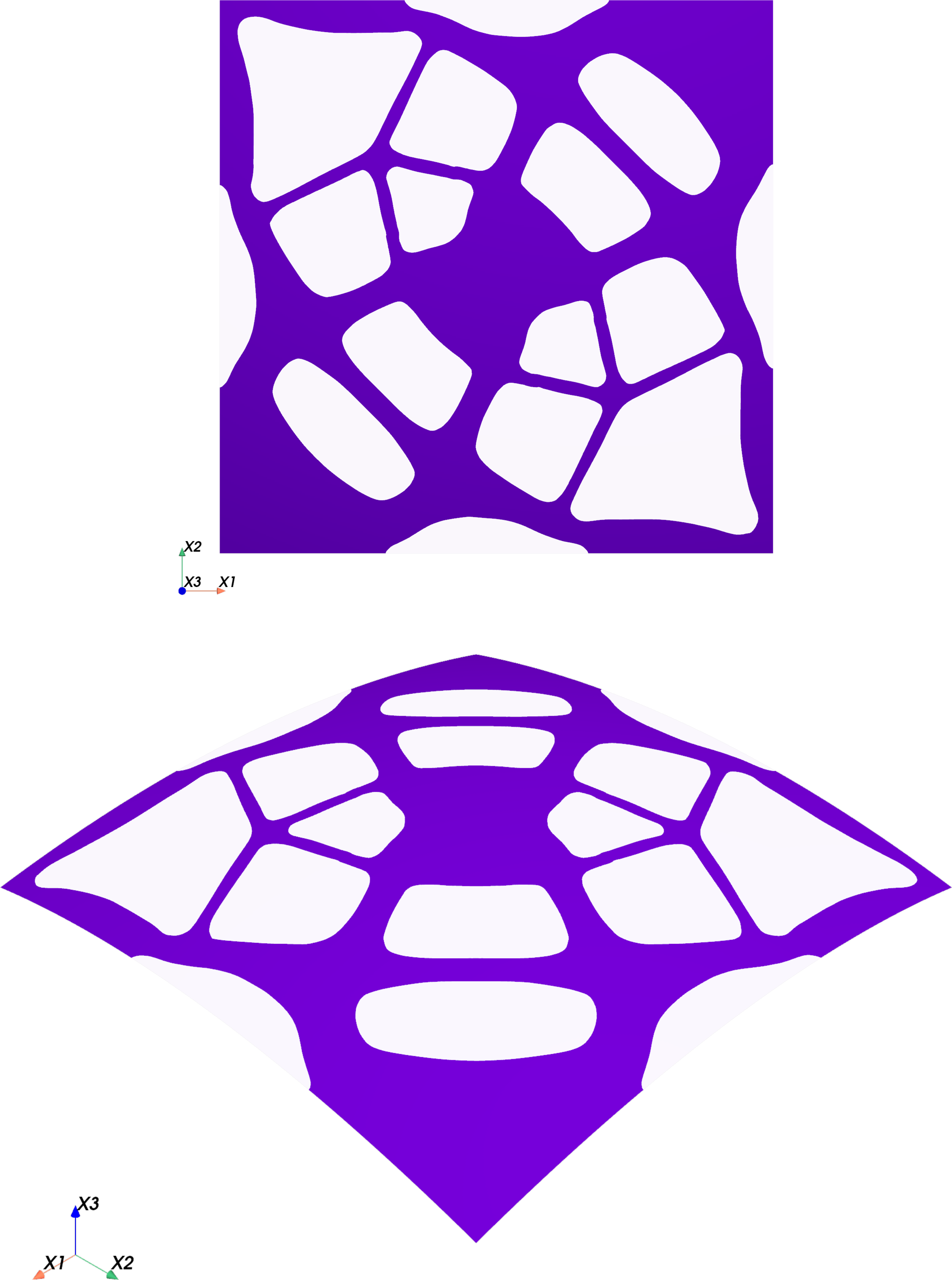}
            \vspace{2pt}
            \footnotesize Diagonal Loading\\ $\boldsymbol{f}=[-\frac{10^4}{\sqrt{2}},-\frac{10^4}{\sqrt{2}},0]$.
        \end{minipage}
        \caption{
Optimized topologies for the double parabola under deterministic point loads. Remarks: 
(a) Vertical load: reinforcement along z-axis with cross-shaped supports.
(b) Horizontal load: asymmetric struts aligned with the load direction.
(c) Diagonal load: diagonal trusses connecting the load point to clamped boundaries.
}   
        \label{fig:dp_res_loadCasesDet}
\end{figure}

\subex[Uncertainty in Loading Direction]
    \label{ex:robDoubleParaRobPara}
    While the previous example shows the applicability of the KL shell based topology optimization for deterministic point loads, the extension to uncertainty in loading direction is studied.
    With the developed optimization algorithms for robust cases, the optimization can handle uncertainty in loading. 
    The problem is based on the examples from \cite{martinez-frutos_robust_2016} and uses the same Gauss Hermite collocation model.  
    However, we treat a shell example and a curved geometry. \\
    We model the uncertainty by deviations in the loading directions, which is normally distributed with a standard deviation of $10^\circ$. 
    The expected loading direction corresponds to the diagonal loading from the previous example. 
    The penalty of the compliance standard deviation $\alpha$ is varied from $0$ to $5$ and the optimization aims for a target volume fraction of $50\%$.

    \begin{figure}[H]
        \centering
\includegraphics[width=0.8\linewidth]{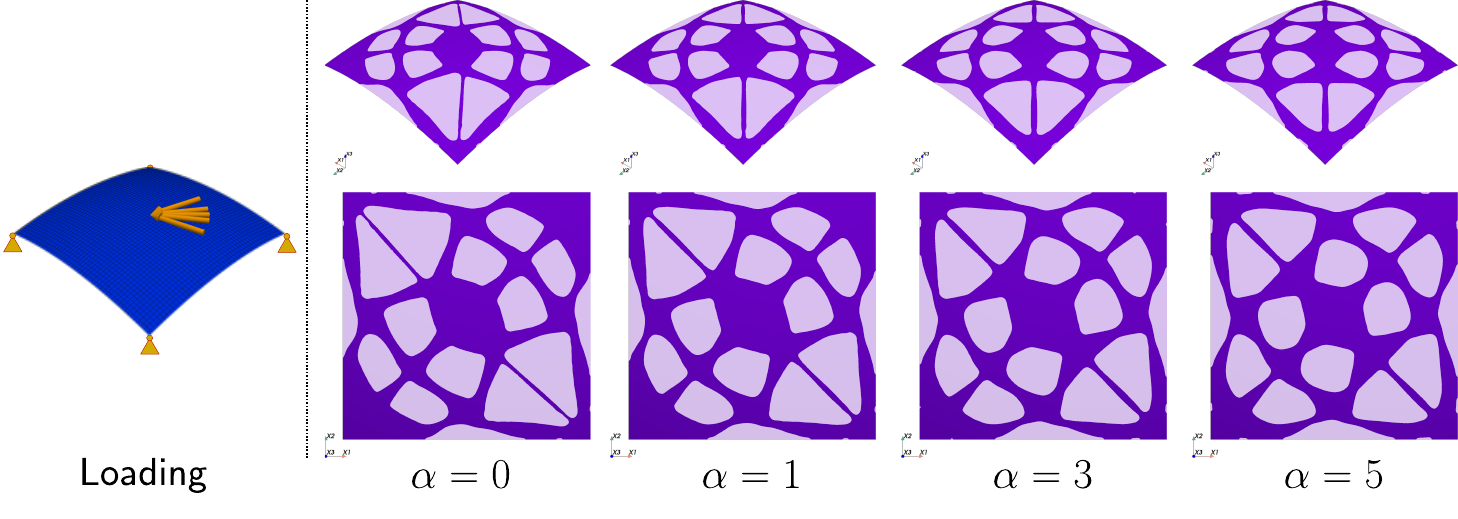}
        \caption{
Double parabola with uncertain loading direction, i.e. load orientation varies with standard deviation $10^\circ$ (Gaussian distribution). Robustness parameter: $\alpha = 0$ (mean-compliance), $1, 3, 5$ (left to right).
Observations: $\alpha = 0$ sparse diagonal supports optimized for mean load, 
$\alpha > 0$ induces thicker diagonal reinforcements to reduce compliance variance.
Higher values of $\alpha$ penalize variance, leading to more conservative designs.
}
\label{fig:res_rob_dp_alpha}
\end{figure}
    \begin{figure}[H]
        \centering
        \begin{minipage}[c]{0.4\textwidth}
            \centering
            \begin{tabular}{c|ccc}
                \toprule
                $\alpha$ & $\mathcal{J}_C^{R*}(\Omega)$ & $\mathbb{E}[{\mathcal{J}}_C^*]$ & $\sqrt{\mathbb{V}[{\mathcal{J}}_C^*]}$  \\ 
                \midrule
                0 & 27.61 & 27.61 & 1.11 \\
                1 & 28.60 & 27.84 & 0.76 \\
                3 & 29.86 & 28.45 & 0.47 \\
                5 & 30.81 & 29.31 & 0.30 \\
                \bottomrule
            \end{tabular}
            \captionof{table}{Comparison between robust compliance, the expected compliance, and its standard deviation for different robustness parameters $\alpha$.}
            \label{tab:compliance_data}
        \end{minipage}
        \hfill 
        \begin{minipage}[c]{0.55\textwidth}
            \centering
    
    
            \includegraphics{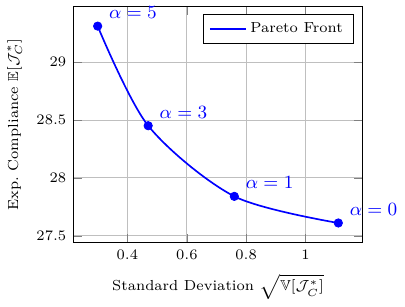}
            \caption{
            Expected Compliance $\mathbb{E}[\mathcal{J}_C]$ (y-axis) vs. compliance standard deviation $\sqrt{\mathbb{V}[\mathcal{J}_C]}$ (x-axis).
            The plot illustrates the trade-off: increasing $\alpha$ reduces variance but increases expected compliance, illustrating the risk-averse design compromise.            }
            \label{fig:pareto_front}
        \end{minipage}
    \end{figure}     
    \autoref{fig:res_rob_dp_alpha} shows the results of this optimization, which show that for a larger parameter $\alpha$, the resistance with respect to less probable force orientations (normal to the loading diagonals) is more optimized for. 
    This translates to greater reinforcement along the diagonal when the parameter $\alpha$ is higher.
    Figure~\ref{fig:pareto_front} shows the trade off between an increase in expected compliance and the variance of the compliance. 
    While an increase in $\alpha$ effectively reduces the variance, it leads to a higher expected compliance. 
      
\begin{figure}[h]
        \centering
        \includegraphics[width=0.7\linewidth]{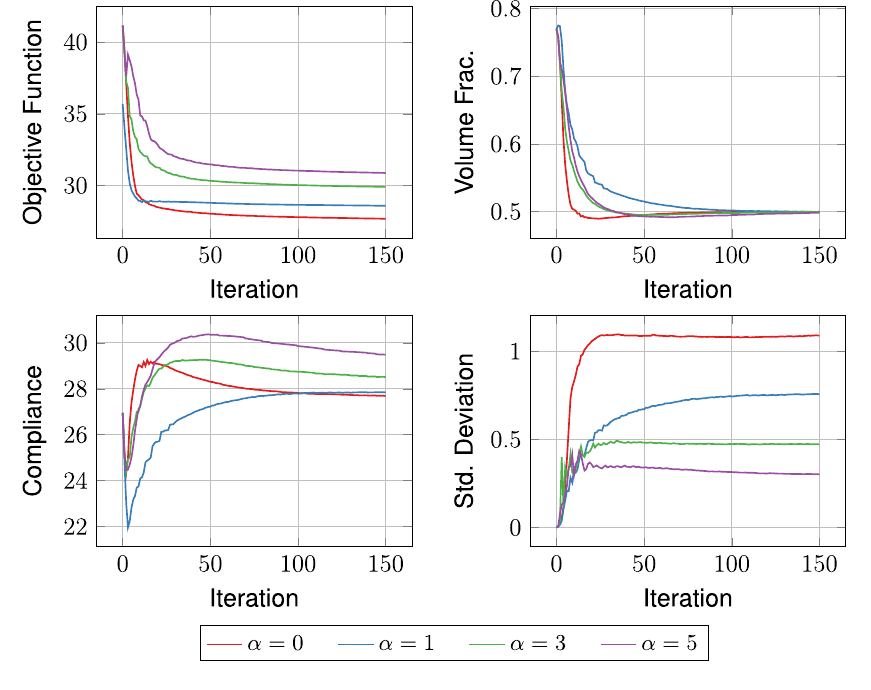}
        \caption{
%
Convergence history for the double parabola with uncertain loading direction for various values of $\alpha = 0,1,3,5$. Objective function, volume fraction, compliance, standard deviation.
        }        
        \label{fig:robDP_conv}
    \end{figure}
    
    \subex[Deterministic Surface Loading]
    \label{ex:dp_det_surf}
    In addition to point load, various engineering applications require the optimization with respect to distributed surface loads, which can for instance correspond to loading cases of building subject to wind. 
    In the following example the ersatz parameter of the weak phase must be adapted to ensure well-posedness of the problem and is chosen as $
    \epsilon_{er}=0.05$. 
    The loading is applied in $\boldsymbol{e}_z$ with magnitude $12500$. 
    The geometry and initialization are equivalent to the previous examples. 

    \begin{figure}[h]
        \includegraphics[width=\linewidth]{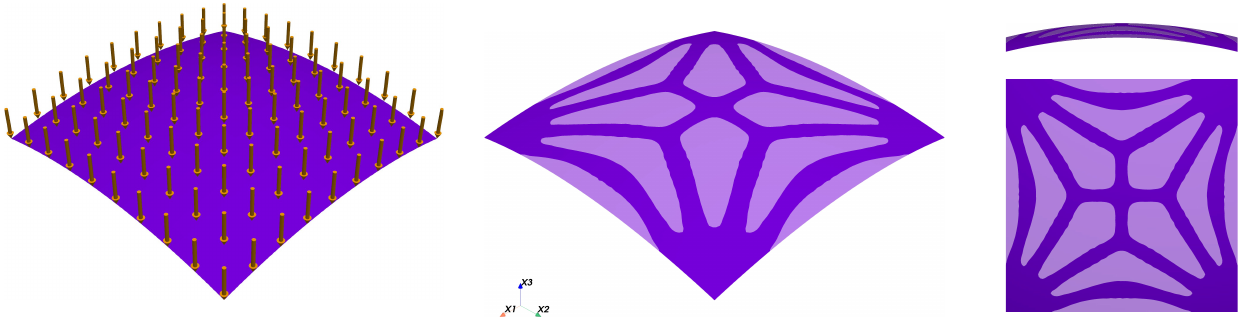}
        \caption{Double Parabola with distributed surface loading. Loading and optimized shape: 3D representation and projection on the axial planes.}
        \label{fig:dp_resDetSurfLoad}
    \end{figure}
    The resulting optimized geometries are shown in \autoref{fig:dp_resDetSurfLoad}. 
    Comparing the result to the result for a point load in vertical direction, the results are very different (see \autoref{fig:dp_res_loadCasesDet}). 
    The distributed load yields a geometry with a much more even distribution of material to ensure the rigidity of the whole surface. 
    The point load on the other hand yields a geometry with a large material concentration in the center to distribute the point as a lack of rigidity does not greatly change the cost function. 
    \FloatBarrier
    \subex[Distributed Surface Loading by GRF]
    \label{ex:robustDP_distLoad}
    For accurate modeling of engineering phenomena, the distributed load must be extended to account for spatial deviations, which could be due to uneven distributions in wind (for example due to vortices). 
    The problem setting remains the same, however, the loading is applied by a Gaussian random field with Matérn covariance and correlation length $\rho = 0.25$, the load is chosen as:
    \begin{equation}
        f(\fx, \omega)= f_0(1+\sigma q(\fx, \omega)), 
    \end{equation}
    where $f_0=12500$, $\sigma = 0.4$ and $q(\fx, \omega)$ corresponds to a sample of a GRF with the desired properties. 
    
%
%
%
%
    \begin{figure}[h]
        \centering
\includegraphics[width= \linewidth]{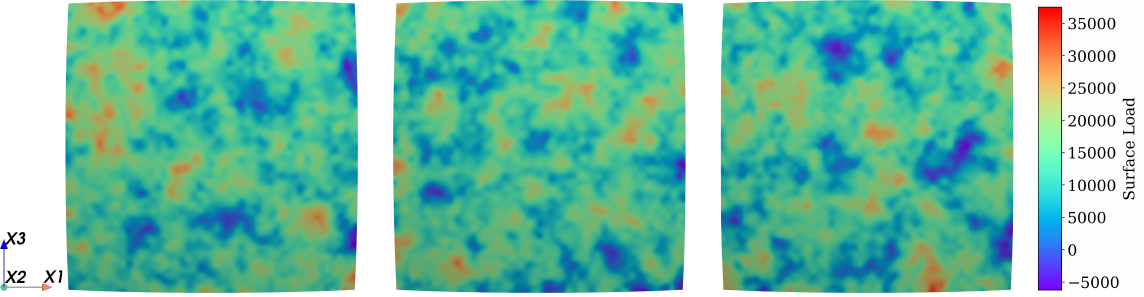}
        \caption{Double parabola with uncertain, GRF surface load. Three samples of the amplitudes of the applied surface loadings}
        \label{fig:loadDpSurf}
    \end{figure}
    \begin{figure}[h]
        \centering
\includegraphics[width= 0.8\linewidth]{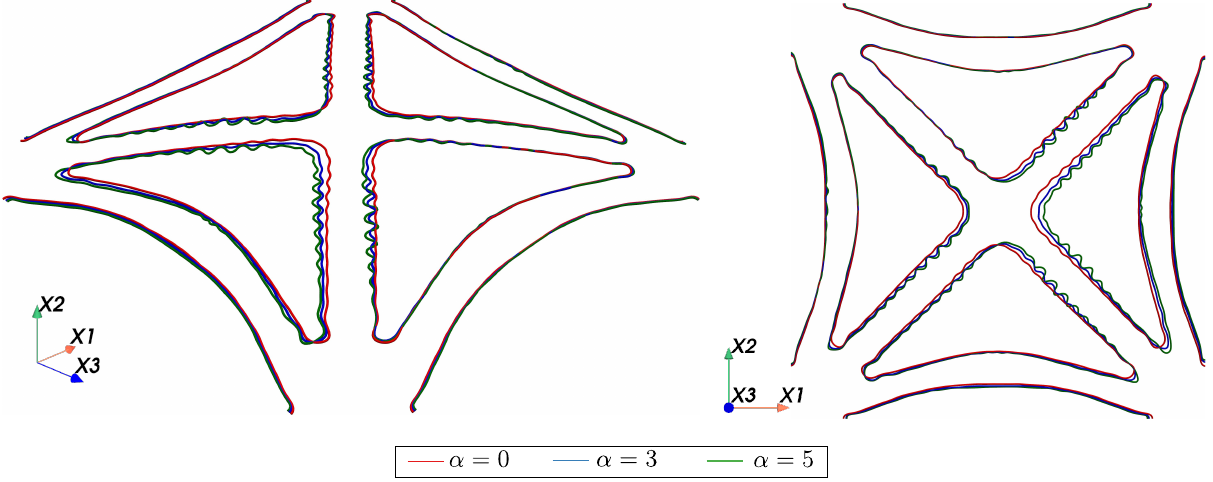}
        \caption{Double parabola with uncertain, GRF surface load. Comparison of the optimized shapes for different robustness values: $\alpha=0,\ 3,\ 5$.}
        \label{fig:resDpSurf}
    \end{figure}
    In \autoref{fig:resDpSurf}, the resulting geometries under uncertainty are shown. 
    At higher robustness levels $\alpha$, the diagonal geometries show a wavy contour. 
    This can be attributed to the developing reinforcements. 
    Previous papers (see for example \cite{duswald_finite_2024,diaz-avalos_formopt_2026}) have shown very fine details for distributed loads. 
    Due to the regularization and the relatively coarse geometries used, these details seem to not develop with the current approach. 
    However, from a manufacturing point of view, very small localized connections, can pose significant problems, thus, geometries without these details can (depending on the application) be preferable. 
    
%
%
%
    \begin{figure}[h]
        \centering
        \includegraphics[width=0.8\linewidth]{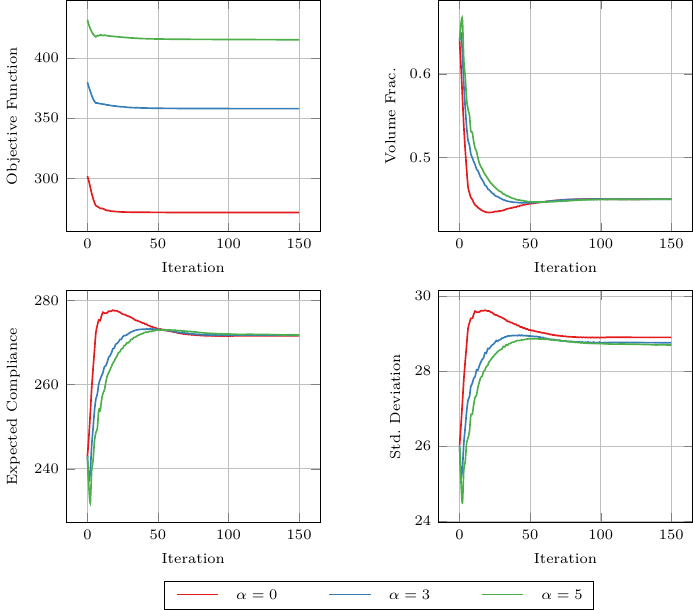}
        \caption{Double parabola with GRF surface load. Development of the compliance, volume fraction, standard deviation and cost function throughout the iterations.}
        \label{fig:dpSurf_costDel}
    \end{figure}
    In \autoref{fig:dpSurf_costDel}, the development of compliance, volume fraction, standard deviation and cost function is shown. 
    The cost function for the robustness level $\alpha=0$ is significantly lower than the cost function at higher robustness level which is due to the penalty of the variance. 
    \subex[Double Parabola with Lognormal Young's Modulus]
    As variations in the Young's modulus distribution have an effect on the resulting stress state, the effect on the resulting optimized geometry should be studied. 
    Similar applications for 2D elasticity problems or full field 3D problems have been studied in \cite{chen_level_2010,martinez-frutos_robust_2016,zheng_levelset_2019}. \\
    As a test example, the double parabola geometry with horizontal loading is studied with the diagonal loading from \autoref{fig:doublePara_geom} and the same initialization and geometry.\\
    Due to energy bounds, at every point inside the material domain $E>0$ and thus, the Young's modulus cannot be modeled by a Normal distribution and is chosen to be lognormally distributed (see \cite{chen_level_2010,martinez-frutos_robust_2016,zheng_levelset_2019}). 

    For the Gaussian random field $\rho=0.5$ and $\gamma=1$. For the MC sampling, the sample size was chosen as 40. 
%
%
\begin{figure}[h]
        \centering
\includegraphics[width= \linewidth]{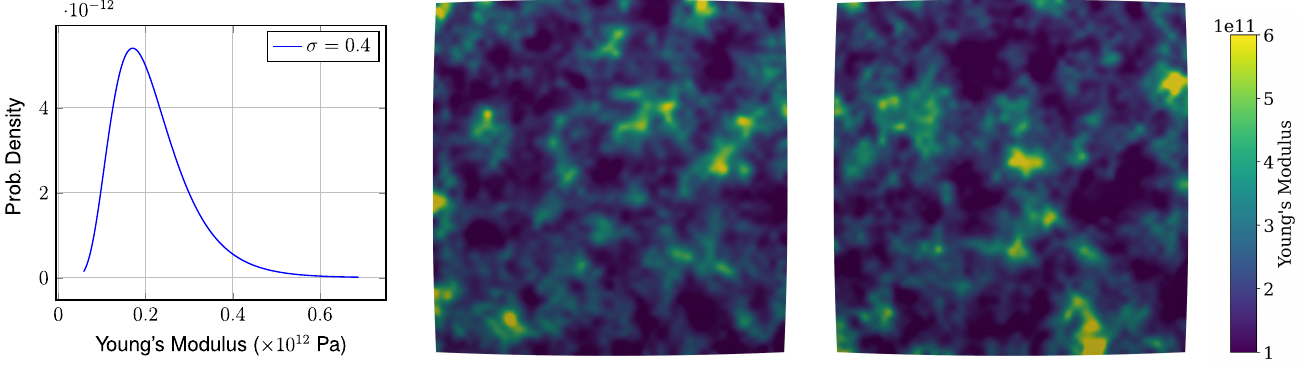}
        \caption{
Lognormal Young’s modulus field for the double parabola to ensure $E > 0$ everywhere with GRF parameters: Matérn covariance, correlation length $\rho = 0.5$, smoothness $\nu = 1$. Plot of the Lognormal distribution and two independent realizations of $E(\boldsymbol{x}, \omega)$. 
}
        \label{fig:robDP_matSamps}
\end{figure}

\begin{figure}[h]
        \centering
\includegraphics[width= 0.7\linewidth]{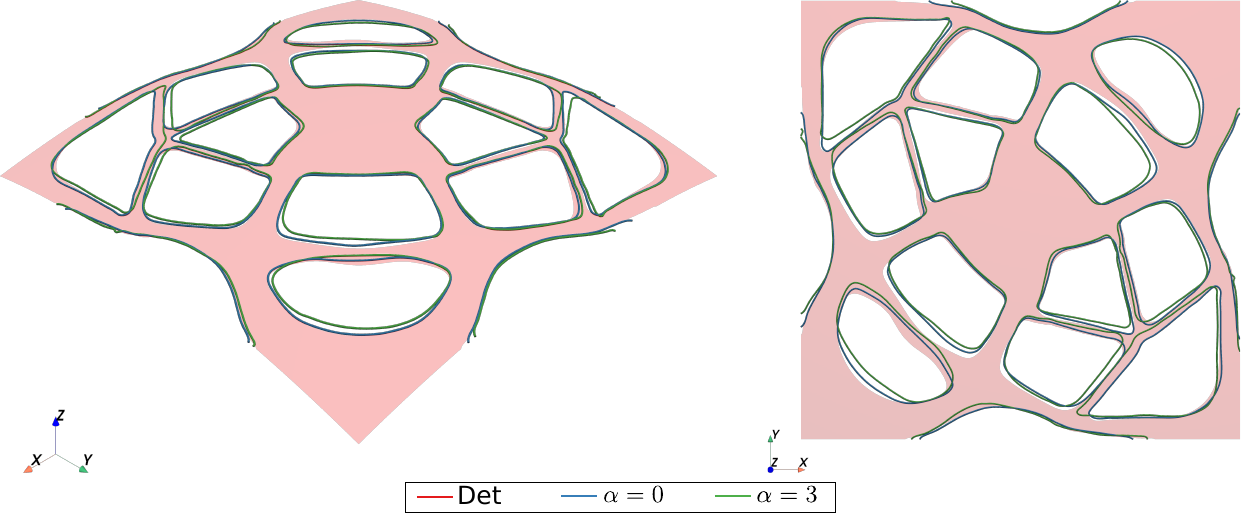}
        \caption{
Optimized topologies for the double parabola with lognormal Young’s modulus for  robustness levels: $\alpha = 0, 3$  and the deterministic result.
Observations, moderate variations in $E$ (due to long $\rho$) leads to minor topology changes and diagonal arms are retained.
}
        \label{fig:resRobDP_lorgnorm}
    \end{figure}
    
    The resulting geometry shown in \autoref{fig:resRobDP_lorgnorm} shares qualitatively most similarities to the previous result from \autoref{fig:dp_res_loadCasesDet}. 
    This is attributed to the moderate variations in Young's modulus and the long correlation length. 
    Furthermore, this demonstrates the robustness of the geometries the deterministic optimization algorithm finds. 

    \subex[Double Parabola with Geometric Perturbation by GRF]
    To model the effect of manufacturing uncertainties on the optimized geometry, an approach similar to the work in \cite{khristenko_statistical_2020,khristenko_statistically_2022} is studied. 
    Here, uncertainty is introduced through a perturbation of the geometry's level set by a Gaussian random field.  
    This can be used to model uncertainty due to fabrication and depending on the fabrication method, the GRF and the associated noise level $\gamma$ can be adapted. 
    The chosen perturbation method can be written as:
    \begin{align}
    \phi_{\text{p}}(\fx, \omega) = (1-\gamma) \phi_{\text{id}}(\fx)+ \gamma\phi_{\text{GRF}}(\fx, \omega), 
    \label{eq:lsMixPerturbation}
    \end{align}
    where $\phi_{\text{id}}$ corresponds to the ideal, unperturbed level set function and $\phi_{\text{GRF}}$ is the noise field depending on the sample indicator $\omega$. 
    Depending on the noise level $\gamma$, the effect of this perturbation field can be tuned. 
    
    \goodbreak
%
%
%

    In \autoref{fig:perturbationDP} this process is visualized. 
    In (a), the unperturbed level set is plotted and in (b), the perturbation field is shown. 
    The resulting sum is shown in (c) with the cut of the perturbation field shown in black and the unperturbed geometry is shown in gray. 
   \begin{figure}[!htbp]
        \centering
        \begin{minipage}[t]{0.2\linewidth}
            \centering
            \includegraphics[width=\linewidth]{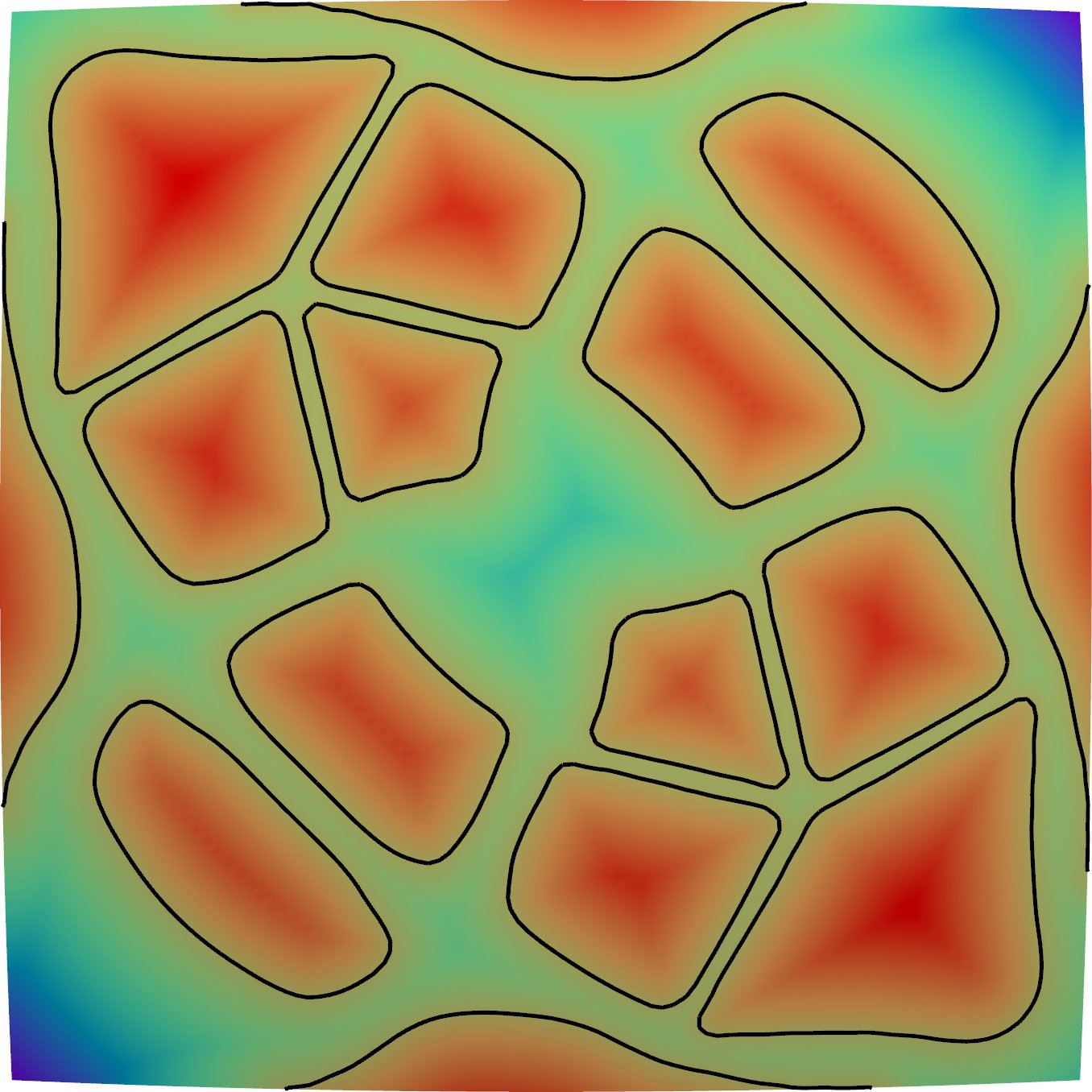}
            \vspace{2pt}
            \footnotesize $\phi_{\text{id}}(\fx)$
        \end{minipage}
        \begin{minipage}[t]{0.2\linewidth}
            \centering
            \includegraphics[width=\linewidth]{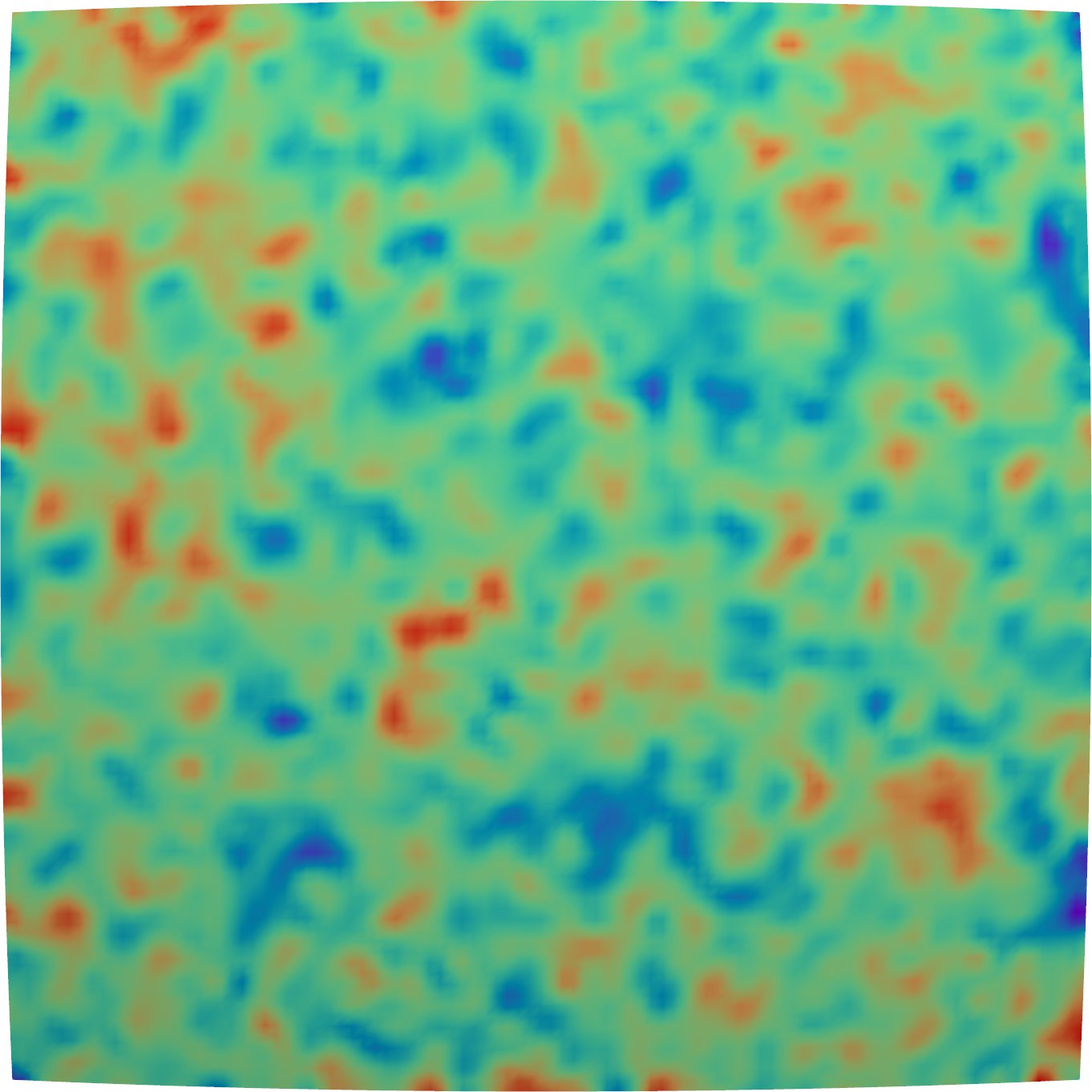}
            \vspace{2pt}
            \footnotesize $\phi_{\text{GRF}}(\fx, \omega)$
        \end{minipage}
        \begin{minipage}[t]{0.2\linewidth}
            \centering
            \includegraphics[width=\linewidth]{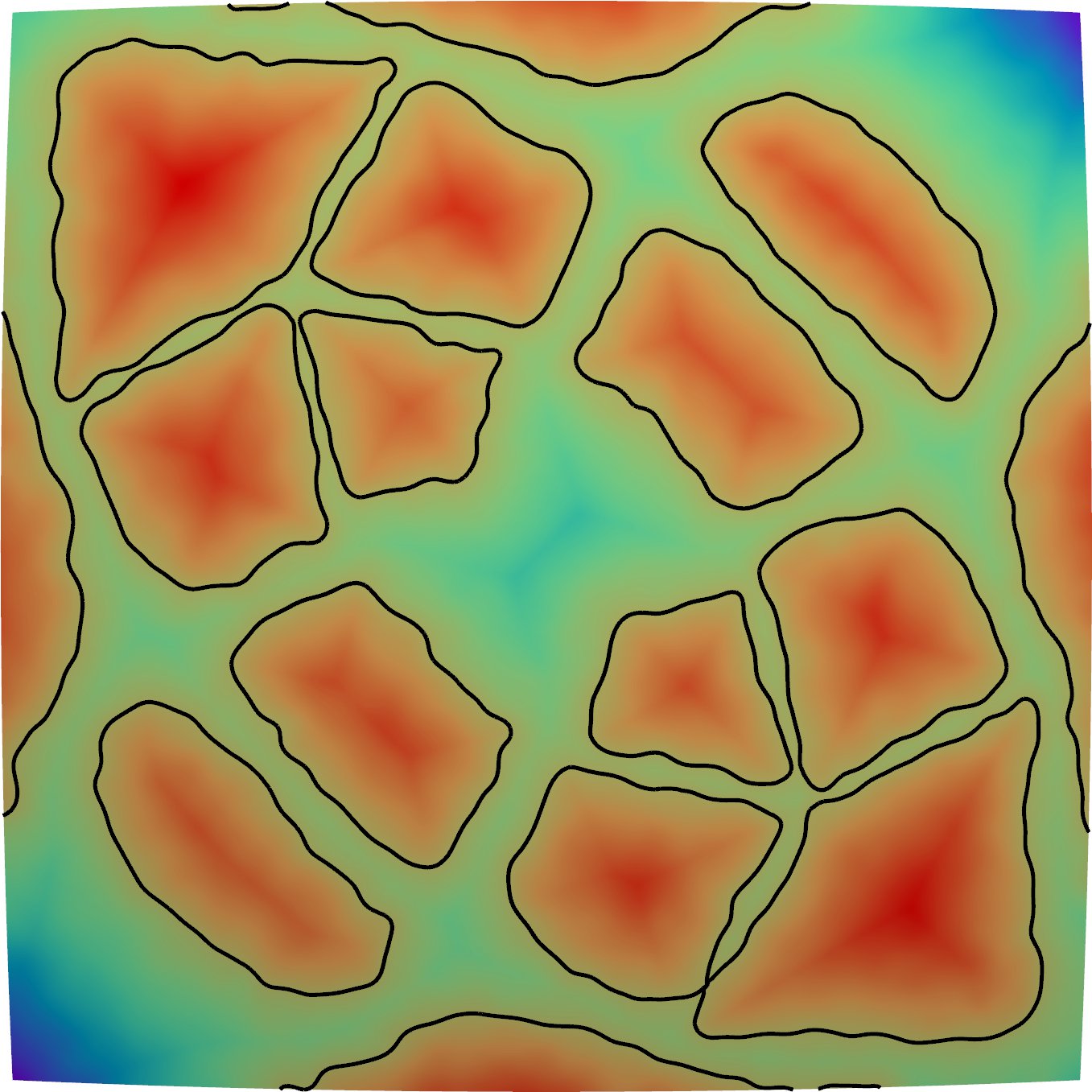}
            \vspace{2pt}
            \footnotesize $(1-\gamma) \phi_{\text{id}}(\fx)+ \gamma\phi_{\text{GRF}}(\fx, \omega)$
        \end{minipage}
        \caption{
Geometric perturbation model for manufacturing uncertainty.
Unperturbed optimized level set  $\phi_{\text{id}}(\boldmath{x})$  (left),  sample from perturbation field with Matérn covariance (center),  sample from perturbed level set with $\gamma = 0.02$ (right).
}
        \label{fig:perturbationDP}
    \end{figure}
     The following examples are based on the deterministic optimization results and the aim of the robust topology optimization is to fine tune the obtained results for a given uncertainty level that might be motivated by the fabrication process and must be adapted in accordance to the geometry of the problem. 
        The noise level set $\phi_{\text{GRF}}$ is normally distributed with unit standard deviation. The topological support, however, follows signed distance properties. 
        Thus, the noise level needs to be chosen depending on the size of the domain.
    To test the effect of the developed perturbation model on the resulting geometry, in a first step, the double parabola geometry given a GRF based perturbation is optimized. 
    Here, the optimization was conducted based on the deterministic optimization result from the deterministic point loading and the fine-tuning was conducted by optimizing with respect to the robust compliance under perturbations of the geometry. \\
    Here, the loading is chosen as diagonal loading with magnitude $F = 4.5\cdot 10^4$ and a target volume fraction of $\varphi^* =0.45$ (see \autoref{fig:doublePara_geom} for the loading).\\

    Comparing the result from \autoref{fig:dp_res_loadCasesDet} to the results under uncertainty shown in \autoref{fig:resDpPerturb}, the change in geometry for the fine arms is visible, which is due to the perturbation field. 
    As the diagonal reinforcements are thinner with respect to the other reinforcements, their sensitivity with respect to perturbations is higher and hence they were removed. 
    \begin{figure}[!htbp]
        \centering
\includegraphics[width= 0.8\linewidth]{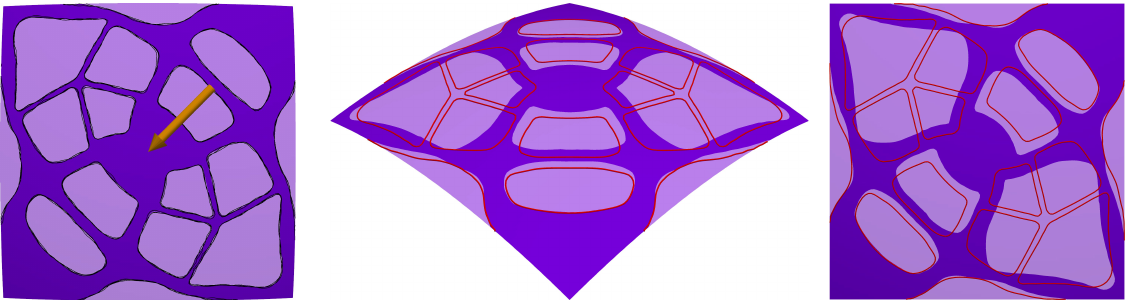}
        \caption{Double parabola with perturbation. Loading with initialization (left) and geometry of topology optimization given perturbation of GRF (center and right). The red contour corresponds to the deterministic optimization results while the blue surface corresponds to the optimization under uncertainty.}
        \label{fig:resDpPerturb}
    \end{figure}
\end{example}

\begin{example}[Multi-Patch Cylinder]
    \label{ex:mpCylinder}
    To demonstrate the effectiveness of the approach for distributed loads, the geometry is now loaded by a line distributed load on the right edge of the cylinder. 
    The load is applied as a torsion-type load with a magnitude of $255$. 
    Furthermore, a target volume fraction of $\phi^*=0.5$ is applied and included in the cost function. This volume fraction enforcement is based on an \emph{augmented Lagrangian} approach, as used in \cite{allaire_optimizing_2018}. See \cite{nocedal_numerical_2006,birgin_practical_2014} for more details on augmented Lagrangian methods and their theoretical foundations. This extends the formulation written in Problem~\ref{eq:optProb}.\\
    The initialization is adapted to ensure the loading edge is fully part of the material phase (see \autoref{fig:ini_cylinderLineLoad}). 

%
%
%

    \begin{figure}[H]
        \centering
\includegraphics[width=0.6\linewidth]{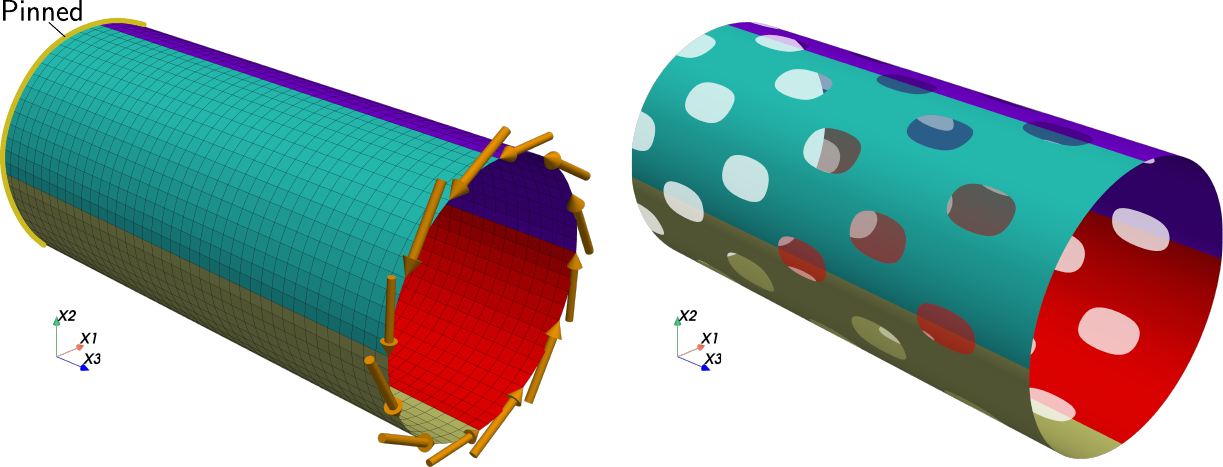}
        \caption{Multi-patch cylinder under torsion: IGA mesh ($40 \times 15 $ elements per patch) and boundary conditions pinned and distributed load on the two extremal edges (left) and initial shape of  the optimization process.}
        \label{fig:ini_cylinderLineLoad}
    \end{figure}

    \subex[Deterministic Torsion of Cylinder]
    \label{ex:cyl_detTors}

    \begin{figure}[H]
        \centering
\includegraphics[width=.8\linewidth]{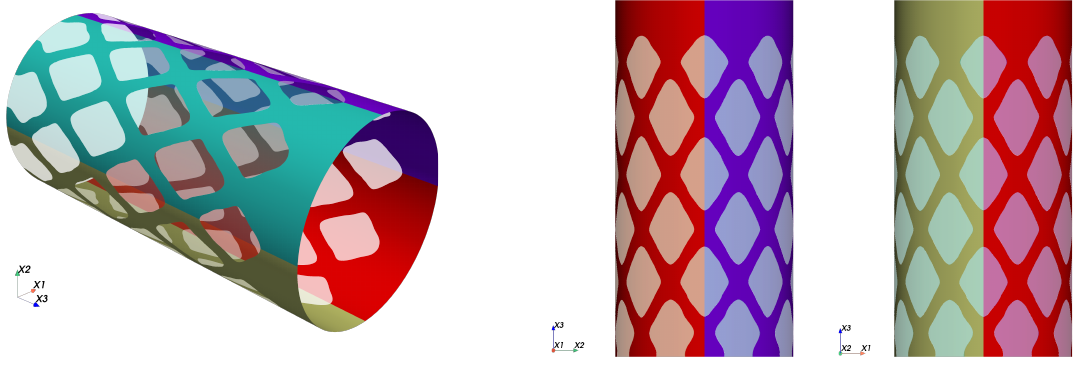}
        \caption{Multi-patch cylinder under torsion. Optimal design obtained at iteration 100. }
        \label{fig:results_cylinderLineLoad}
    \end{figure}
%
%
%
     
    \begin{figure}[h]
        \centering
        \includegraphics[width=0.8\linewidth]{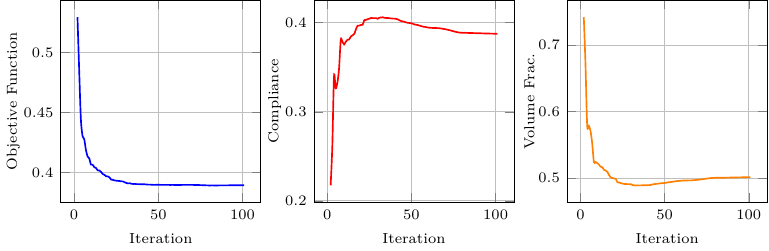}
        \caption{Multi-patch cylinder under torsion with pinned and distributed torsional load on extremal edges: histories of the cost-functional, compliance and volume fraction. }
        
    \end{figure}
%
In \autoref{fig:results_cylinderLineLoad}, the final geometry at iteration 100 is visualized. 
We remark, that the resulting geometry is visually very similar to the geometry obtained after 30 iterations, which can also be seen in the cost function development having the largest decrease over the first 30 iterations. 

 \subex[Torsion with axial loading by GRF]
    \label{ex:robustCylinderRobComp}
    While the previous example computes optimized geometries for pure axial loadings, axial components might be present in real world scenarios.
    To extend the methodology for this, a GRF with Matérn covariance is used to superimpose a random axial loading component. \\
   On the left of Figure~\ref{fig:smoothnessEffect} the GRF based loading component is shown for a correlation length of $\rho = 0.2$ and a smoothness of $\gamma=[1, 3]$, corresponding to a GRF generated by a harmonic and biharmonic SPDE.

     \begin{figure}[htb]
        \centering
        \begin{minipage}{0.25\textwidth}
            \centering
            \includegraphics[width=\linewidth]{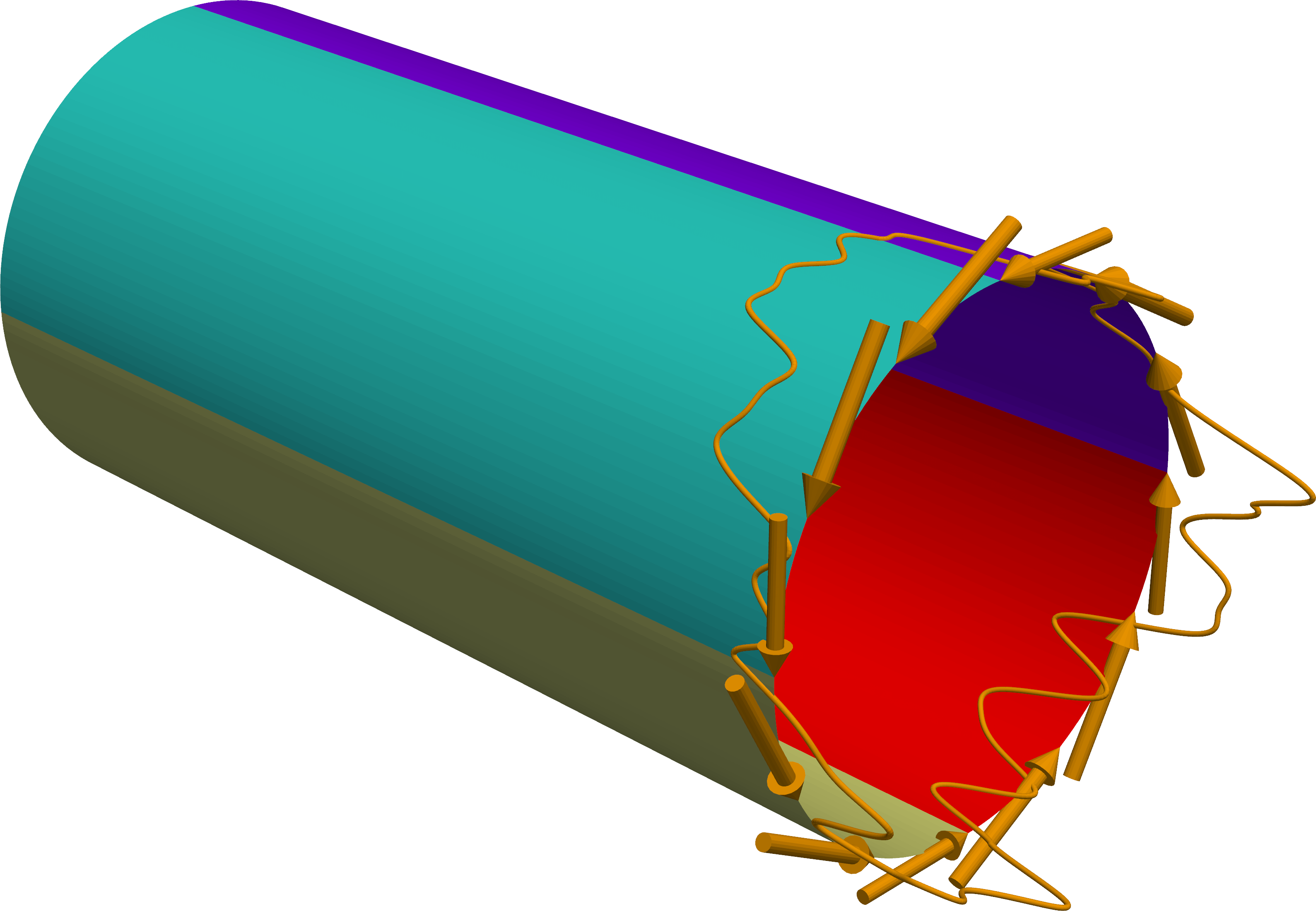}\\
            \vspace{2pt}
            \footnotesize $\rho = 0.2$, $\gamma = 1$
        \end{minipage}\hfill
        \begin{minipage}{0.35\textwidth}
            \centering
            \includegraphics[width=\linewidth]{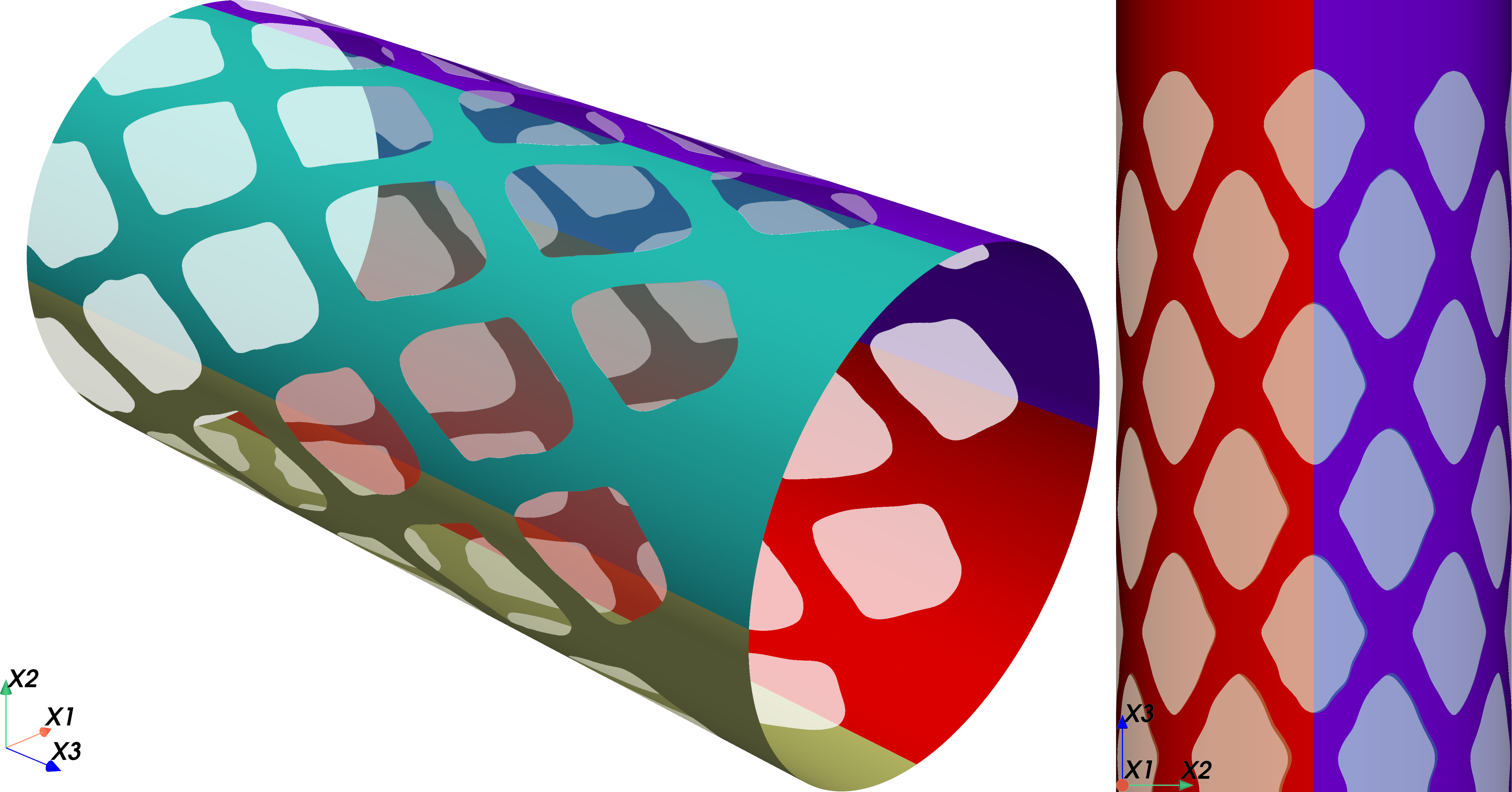}
            \vspace{2pt}
            \footnotesize $\alpha=0$, $\gamma=1$, $\rho=0.2$
        \end{minipage}\hfill
        \begin{minipage}{0.35\textwidth}
            \centering
            \includegraphics[width=\linewidth]{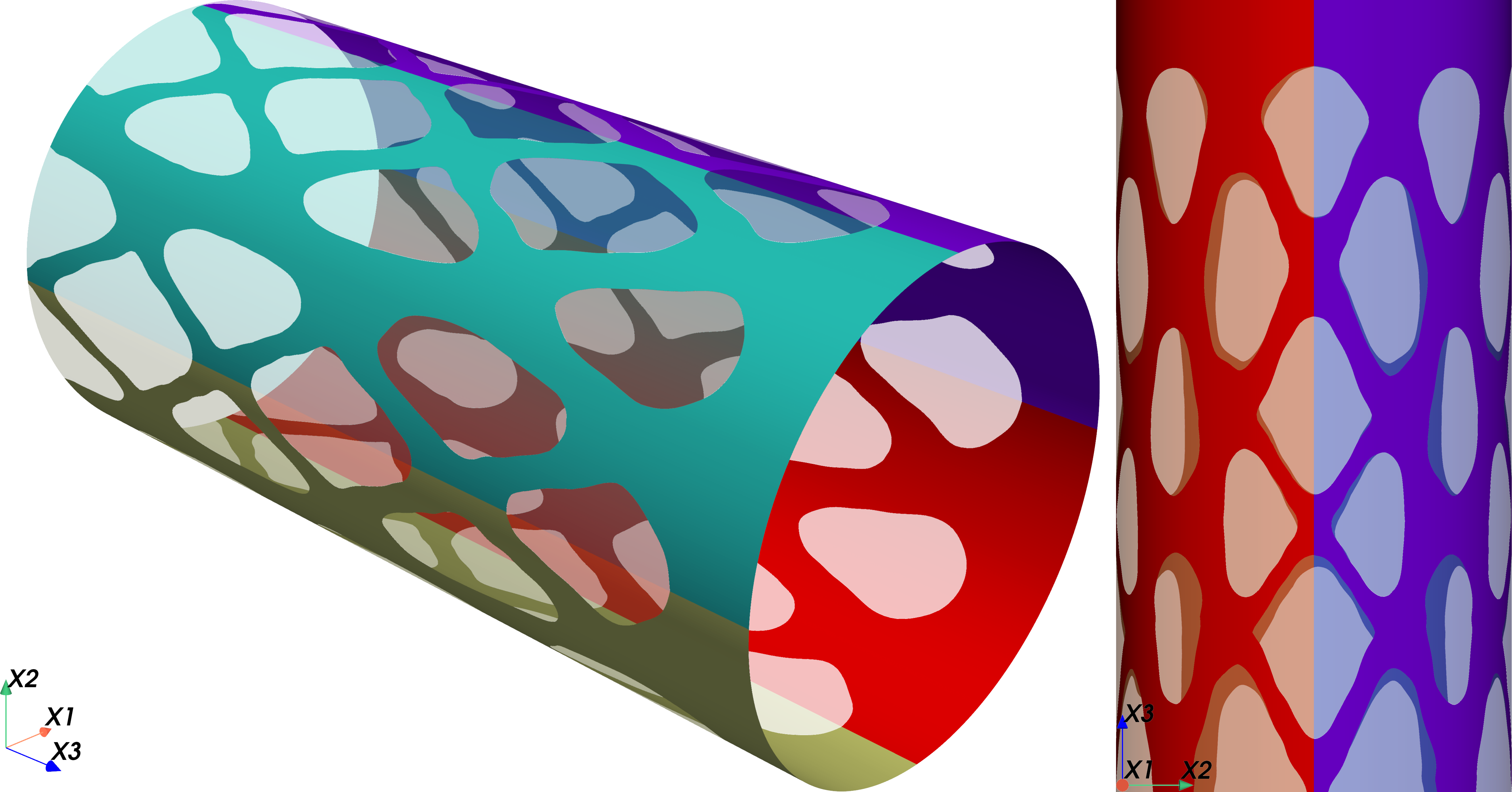}
            \vspace{2pt}
            \footnotesize $\alpha=10$, $\gamma=1$, $\rho=0.2$
        \end{minipage}
        \caption{
Multi-patch cylinder under torsion with pinned and distributed torsional load (with stochastic axial component) on the extremal edges respectively:   Optimized shapes for the deterministic case, $\alpha = 0,\ 10$ and fixed correlation length $\rho$ and smoothness $\gamma$. }
        \label{fig:smoothnessEffect}
    \end{figure}
\end{example}



\section{Conclusion}
This work presented a CAD-compatible topology optimization framework for
Kirchhoff--Love shells, accounting for non-matching multi-patch geometries. The
approach relies on distributed, volume-based shape derivatives, aligning with the
isogeometric paradigm of coarse, high-order discretizations. This formulation
avoids Dirac-type boundary approximations and eliminates the need for trimming,
local refinement, or selective integration. Because the bending strains involve
second derivatives of the velocity field, the shape derivative is regularized
through a biharmonic problem posed over the surface, which is solved with the same
multi-patch discretization as the state problem and is therefore compatible with
standard IGA solvers.

The method reproduces literature results using significantly fewer elements than
comparable FEM or refined IGA techniques: for the multi-patch hemisphere,
$40\times 40$ cubic NURBS elements per patch recover the reinforcement pattern
obtained in \cite{dede_isogeometric_2012} with $128\times128\times1$ solid
elements. The framework has been validated on single- and multi-patch geometries
under point, line and surface loading, and under uncertainty in loading direction,
loading amplitude, material properties and geometry.

Our approach models uncertainty in loading direction and amplitude for multi-patch
geometries at the cost of one linear solve per sample. For the double parabola
under uncertain loading direction, increasing the robustness parameter from
$\alpha=0$ to $\alpha=5$ reduces the standard deviation of the compliance by
roughly $73\%$ at the price of a $6\%$ increase in its expectation, and produces
thicker reinforcements aligned with the diagonal, i.e.\ material is redistributed
towards the less probable load orientations. In the limit of a vanishing input
dispersion the mean-compliance design $\alpha=0$ recovers the deterministic
solution; we emphasise that $\alpha=0$ is not itself the deterministic case, since
it still minimises the expected compliance over the sampled loads.

Matérn-type random fields provide realistic models of spatially correlated loading
and material data; future work could calibrate these fields to specific
applications. The advantage of the robust compliance lies in its interpretability,
as both terms carry the units of the compliance itself. 
Alternative cost functionals should be explored,
particularly peak von Mises stress to reduce the risk of plastic deformation, and
Conditional Value-at-Risk for managing extreme events; both are accessible within
the present framework, since the shape derivative enters only through a
realization-dependent scalar weight.

The optimized geometries obtained under material uncertainty show minimal
qualitative differences from their deterministic counterparts. The topology is
largely unchanged and the diagonal arms are retained, with the boundaries adapting
only locally as $\alpha$ increases. We attribute this to the moderate dispersion
and the long correlation length $\rho=0.5$ relative to the structural members: the
sampled fields vary slowly over the scale of the reinforcements, so each
realization sees an almost uniform effective stiffness. This is consistent with the
loading results, where the topology changes appreciably only when the uncertainty
alters the character of the load, as when a stochastic axial component is
superimposed on the torsion of the cylinder. It also suggests that the designs
found by the deterministic algorithm are already robust to this class of
perturbation. By contrast, the geometric perturbation model acts at the scale of
the members themselves and does modify the topology, removing the thinnest diagonal
reinforcements, whose sensitivity to boundary perturbations is highest.

Comparable studies report finer geometric detail under loading and material
uncertainty than we obtain here
\cite{martinez-frutos_risk-averse_2021,duswald_finite_2024}. Our coarse
discretization and the biharmonic regularization together impose a length scale
below which such features cannot develop. From a manufacturing standpoint this is
not necessarily a drawback, as very small localized connections are difficult to
fabricate; resolving them would nonetheless require finer discretizations, a
weaker regularization, or trimming-based integration.

The principal limitation of the approach is the cost of Monte Carlo sampling, which
scales the cost of each iteration by the sample size. The samples are independent
and the implementation is amenable to parallelization, but the Python-based Bézier
extraction dominates the per-iteration overhead in the present code. Variance
reduction or multilevel sampling strategies would reduce this cost. The theoretical
framework is independent of the implementation and transfers directly to more
efficient environments. Future work will address trimming-based integration for
improved accuracy near the interface, and stress-based and risk-averse objectives.

The numerical implementation is open-source and available at
\url{https://github.com/figaroIGA/figaro}. 

\section*{Acknowledgements}
This project was provided with computer and storage resources by GENCI at IDRIS thanks to the grant 2025-AD010917305 on the supercomputer Jean Zay's CSL partition and was granted access to the HPC resources of IDCS support unit from Ecole polytechnique.
This work has partially been supported by the \textsc{Design and modeling of innovative materials} Chair at Ecole polytechnique, sustained by Saint-Gobain and PE would like to acknowledge financial support through a PhD fellowship by Institut Polytechnique de Paris. 

\newpage

\begin{appendices}

\section{Derivation of result (\ref{eqn:shapederivvol})}
\label{sec:demSDvol}
By definition, the shape derivative of $\mathcal{J}_V$ is 

\begin{equation}
    \mathcal{J}_V'(\Omega)(\boldsymbol{\theta}) = \lim_{t \searrow 0} \frac{1}{t}\left( \int_{\Omega_t} \diff \mathbf{X}_t -  \int_{\Omega} \diff \mathbf{X} \right)
\end{equation}
Mapping to the parametric space, this equation yields

\begin{equation}
    \mathcal{J}_V'(\Omega)(\boldsymbol{\theta}) = \lim_{t \searrow 0} \frac{1}{t}\left( \int_{\hat{\Omega}_t} \hat{g} (\boldsymbol{\xi}_t)\diff \boldsymbol{\xi}_t -  \int_{\hat{\Omega}} \hat{g}(\boldsymbol{\xi})\diff \boldsymbol{\xi} \right) =  \int_{\hat{\Omega}} \frac{\partial}{\partial t } \left( (\hat{g} \circ T_t) \psi(t)\right)\bigg|_{t=0} \diff \boldsymbol{\xi} 
\end{equation}

Using the product rule, we get

\begin{equation}
    \mathcal{J}_V'(\Omega)(\boldsymbol{\theta}) =  \int_{\hat{\Omega}} ( \hat{g} \circ T_t )\frac{\partial}{\partial t }\left( \psi(t)\right)\bigg|_{t=0} +\frac{\partial}{\partial t } \left( \hat{g}\circ T_t \right)\bigg|_{t=0} \psi(0)\diff \boldsymbol{\xi}
\end{equation}

The term $\frac{\partial}{\partial t } \left( \hat{g}\circ T_t \right)$ is developed with the chain rule as
\begin{equation}
    \frac{\partial}{\partial t } \left( \hat{g}\circ T_t \right)\Bigg|_{t=0} = \left( \frac{\partial (\hat{g}\circ T_t)}{\partial \xi^\alpha }\frac{\partial \xi^\alpha}{\partial (\xi_t)^\beta } \frac{\partial (\xi_t)^\beta}{\partial t}\right)\Bigg|_{t=0} = \hat{\nabla}\hat{g} \cdot \hat{\boldsymbol{\vartheta}}
\end{equation}

With Lemma 2.31 from \cite{sokolowski_introduction_1992} divergence on the surface $\operatorname{div}$ the desired result is obtained by pushing forward to the physical domain:

\begin{equation}
\begin{split}
    \mathcal{J}_V'(\Omega)(\boldsymbol{\theta}) &=  \int_{\hat{\Omega}}  \hat{g} \widehat{\text{ div}} \hat{\boldsymbol{\vartheta}}  +\hat{\nabla}\hat{g} \cdot \hat{\boldsymbol{\vartheta}} \diff \boldsymbol{\xi} = \int_{\hat{\Omega}}  \frac{1}{\hat{g}} \widehat{\text{ div}} \left( \hat{g} \hat{\boldsymbol{\vartheta}} \right) \hat{g} \diff \boldsymbol{\xi}  \\ &= \int_{\hat{\Omega}}  \frac{1}{\hat{g}} \widehat{\text{ div}} \left( \hat{g} \hat{G}^{-1}\hat{F}^T \hat{F} \hat{\boldsymbol{\vartheta}} \right) \hat{g} \diff \boldsymbol{\xi}= \int_{\Omega} \text{ div} \left(F \boldsymbol{\vartheta} \right) \diff \boldsymbol{X}.
\end{split}
\label{eqn:divergence}
\end{equation}

\section{Derivation of results (\ref{eqn:dotepsilon}), (\ref{eqn:dotkappa}) and (\ref{eqn:cdot})}
\label{sec:technicalLemma}

To obtain the results stated in Equations (\ref{eqn:dotepsilon}), (\ref{eqn:dotkappa}) and (\ref{eqn:cdot}), the membrane and bending strains are written over a pseudo time $t$.

\textbf{Rewriting the parametric gradient operator.} The function composition $\hat{\boldsymbol{u}}^t(\boldsymbol{\xi}) = (\hat{\boldsymbol{u}}_t \circ T_t)(\boldsymbol{\xi}) = \hat{\boldsymbol{u}}_t(T_t(\boldsymbol{\xi}))$ ties a relation between the gradients of $\hat{\boldsymbol{u}}^t$ and $\hat{\boldsymbol{u}}_t$ by the chain rule and the Jacobian of $T_t$ with

\begin{align}
\begin{split}
    \frac{\partial \hat{\boldsymbol{u}}^t}{\partial \xi^\alpha} = \frac{\partial}{\partial \xi^\alpha} \left(\hat{ \boldsymbol{u}}_t \circ T_t \right) =\left( \frac{\partial \hat{\boldsymbol{u}} _t}{\partial (\xi_t)^\mu} \circ T_t\right) \frac{\partial (\xi_t)^\mu}{\partial\xi^\alpha}=\left( \frac{\partial \hat{\boldsymbol{u}}_t}{\partial (\xi_t)^\mu} \circ T_t\right) \hat{J}^\mu_{\cdot \alpha}.
    \label{eqn:changeOfCoordsGradient}
\end{split}
\end{align}

\textbf{Rewriting the parametric Hessian operator.} 
The same logic is applicable to the second derivative of $\hat{\boldsymbol{u}}^t$:

\begin{align}
\begin{split}
    \frac{\partial^2 \hat{\boldsymbol{u}}^t}{\partial \xi^\beta \partial \xi^\alpha} &= \frac{\partial}{\partial \xi^\beta} \left( \left( \frac{\partial \hat{\boldsymbol{u}}_t}{\partial (\xi_t)^\mu} \circ T_t\right) \frac{\partial (\xi_t)^\mu}{\partial\xi^\alpha}\right)\\
    &=   \left( \frac{\partial^2 \hat{\boldsymbol{u}}_t}{\partial (\xi_t)^\lambda (\xi_t)^\mu} \circ T_t\right) \frac{\partial (\xi_t)^\lambda}{\partial\xi^\beta} \frac{\partial (\xi_t)^\mu}{\partial\xi^\alpha} +\left( \frac{\partial \hat{\boldsymbol{u}}_t}{\partial (\xi_t)^\mu} \circ T_t\right) \frac{\partial^2 (\xi_t)^\mu}{\partial\xi^\beta \partial\xi^\alpha}  \\
\end{split}
\end{align}

On the term on the right, we recognize the quantity $\frac{\partial \hat{\boldsymbol{u}}_t}{\partial (\xi_t)^\mu} \circ T_t$ also present parametric gradient transformation in \autoref{eqn:changeOfCoordsGradient}. The trick here is to artificially add the Jacobian by means of the identity 

\begin{equation}
\delta^\alpha_{ \beta} = \hat{J}^\alpha_{\cdot \mu} (\hat{J}^{-1})^\mu_{\cdot\beta} . 
\label{eqn:jacobianIdentity}
\end{equation}

This allows us to write to further simplify the second derivative to:

\begin{align}
\begin{split}
    \frac{\partial^2 \hat{\boldsymbol{u}}^t}{\partial \xi^\beta \partial \xi^\alpha} =   \left( \frac{\partial^2 \hat{\boldsymbol{u}}_t}{\partial (\xi_t)^\lambda (\xi_t)^\mu} \circ T_t\right) \hat{J}^\lambda_{\cdot \beta} \hat{J}^\mu_{\cdot \alpha} + \frac{\partial \hat{\boldsymbol{u}}^t}{\partial \xi^\lambda} (\hat{J}^{-1})^\lambda_{\cdot\mu} \hat{H}^\mu_{\alpha\beta} 
    \label{eqn:changeOfCoordsSecondGradient}
\end{split}
\end{align}

\textbf{Rewriting the membrane strain. }  In a similar fashion, the membrane strain in \autoref{eqn:membraneStrainSmallDispl} is rewritten over $\hat{\Omega}$ via the definition

\begin{equation}
    (\hat{\varepsilon}^t)_{\alpha\beta}(\hat{\boldsymbol{u}}^t):= (\hat{\varepsilon}_t)_{\alpha\beta}(\hat{\boldsymbol{u}}_t)\circ T_t.
\end{equation}

More precisely, 
\begin{equation}
    (\hat{\varepsilon}^t)_{\alpha \beta}(\hat{\boldsymbol{u}}^t) = \frac{1}{2} \left(((\hat{\mathbf{A}}_t)_\alpha \circ T_t ) \cdot \left(\frac{\partial \hat{\boldsymbol{u}}_t}{\partial (\xi_t) ^\beta} \circ T_t \right)  + \left( (\hat{\mathbf{A}}_t)_\beta \circ T_t \right) \cdot \left( \frac{\partial \hat{\boldsymbol{u}}_t}{\partial (\xi_t) ^\alpha} \circ T_t\right) \right).
\end{equation}
One can recognize the term $\frac{\partial \hat{\boldsymbol{u}}_t}{\partial (\xi_t) ^\alpha} \circ T_t$ also present in the parametric gradient transformation in \autoref{eqn:changeOfCoordsGradient}. By inserting the identity $
\delta^\alpha_{ \beta} = \hat{J}^\alpha_{\cdot \mu} (\hat{J}^{-1})^\mu_{\cdot\beta} $ into the expression, the membrane strain $(\varepsilon^t)_{\alpha \beta} $ is fully described by quantities defined over $\hat{\Omega}$:

\begin{equation}
    (\hat{\varepsilon}^t)_{\alpha \beta}(\hat{\boldsymbol{u}}^t) = \frac{1}{2} \left((\hat{\mathbf{A}}^t)_{\alpha} \cdot \left(\frac{\partial \hat{\boldsymbol{u}}^t}{\partial \xi ^\mu} (\hat{J}^{-1})^\mu_{\cdot \beta} \right)  + (\hat{\mathbf{A}}^t)_{\beta} \cdot \left(\frac{\partial \hat{\boldsymbol{u}}^t}{\partial \xi ^\mu} (\hat{J}^{-1})^\mu_{\cdot \alpha} \right) \right).
\end{equation}

\textbf{Rewriting the bending strain. }  Let us consider writing the expression of curvature tensor $(\hat{\kappa}^t)_{\alpha\beta}$ with the following definition:

\begin{equation}
    (\hat{\kappa}^t)_{\alpha\beta}(\hat{\boldsymbol{u}}^t) := (\hat{\kappa}_t)_{\alpha\beta} (\hat{\boldsymbol{u}}_t)\circ T_t 
\end{equation}

Replacing $(\kappa_t)_{\alpha\beta}(\hat{\boldsymbol{u}}_t)$ by its linearized expression in \autoref{eqn:bendingStrainSmallDispl}, we obtain the following

\begin{equation}
    (\hat{\kappa}^t)_{\alpha\beta}(\hat{\boldsymbol{u}}^t)=-((\hat{\mathbf{A}}_t)_3\circ T_t)
    \cdot \left(\frac{\partial ^2  \hat{\boldsymbol{u}}_t}{\partial (\xi_t)^\alpha\partial(\xi_t)^\beta} \circ T_t- \left((\hat{\Gamma}_t)^\lambda_{\alpha\beta} \circ T_t\right) \frac{\partial \hat{\boldsymbol{u}}_t}{\partial (\xi_t)^\lambda} \circ T_t\right)
\end{equation}

Once again, we apply the identity $
\delta^\alpha_{ \beta} = \hat{J}^\alpha_{\cdot \mu} (\hat{J}^{-1})^\mu_{\cdot\beta}$ to rewrite the terms $\frac{\partial \hat{\boldsymbol{u}}_t}{\partial (\xi_t)^\alpha } \circ T_t$ and $ \frac{\partial ^2  \hat{\boldsymbol{u}}_t}{\partial (\xi_t)^\alpha\partial(\xi_t)^\beta}\circ T_t$ with the equations (\ref{eqn:changeOfCoordsGradient}) and (\ref{eqn:changeOfCoordsSecondGradient}) of the transformed parametric gradient and Hessian, respectively. This yields the following function over $\hat{\Omega}$:

\begin{equation}
\begin{aligned}
    (\hat{\kappa}^t)_{\alpha\beta}(\hat{\boldsymbol{u}}^t)
     &=  && -(\hat{\mathbf{A}}^t)_3\cdot \left( \frac{\partial ^2  \hat{\boldsymbol{u}}^t}{\partial \xi^\gamma \partial \xi^\eta} - \frac{\partial \hat{\boldsymbol{u}}^t}{\partial \xi^\mu} (\hat{J}^{-1})^\mu_{\cdot \lambda} \hat{H}^\lambda_{\gamma \eta} \right) (\hat{J}^{-1})^\gamma_{\cdot \beta}(\hat{J}^{-1})^\eta_{\cdot \alpha} \\& &&+(\hat{\mathbf{A}}^t)_3\cdot \left( (\Gamma^t)^\lambda_{\alpha\beta} 
     \frac{\partial \hat{\boldsymbol{u}}^t}{\partial \xi^\mu} (\hat{J}^{-1})^\mu_{\cdot \lambda}\right).
\end{aligned}
\end{equation}

\textbf{Membrane strain time derivative} Consider the expression of the membrane strain as written in terms of the original parametric domain $\hat{\Omega}$ and a displacement field $\hat{\boldsymbol{u}}^t\in H^2(\hat{\Omega})$:

\begin{equation}
      (\hat{\varepsilon}^t)_{\alpha \beta}(t, \hat{\boldsymbol{u}}^t)=  \frac{1}{2} \left((\hat{\mathbf{A}}^t)_{\alpha} \cdot \left(\frac{\partial \hat{\boldsymbol{u}}^t}{\partial \xi ^\mu} (\hat{J}^{-1})^\mu_{\cdot \beta} \right)  + (\hat{\mathbf{A}}^t)_{\beta} \cdot \left(\frac{\partial \hat{\boldsymbol{u}}^t}{\partial \xi ^\mu} (\hat{J}^{-1})^\mu_{\cdot \alpha} \right) \right)
\end{equation}

By differentiating over time, we get 

\begin{equation}
\begin{aligned}
    \left(\frac{\partial}{\partial t}(\hat{\varepsilon}^t)_{\alpha\beta}\right)\Bigg|_{t=0} = \frac{1}{2} \Bigg(& \frac{\partial}{\partial t}(\hat{\mathbf{A}}^t)_{\alpha} \cdot \left(\frac{\partial \hat{\boldsymbol{u}}^t}{\partial \xi ^\mu} (\hat{J}^{-1})^\mu_{\cdot \beta} \right)  + \frac{\partial}{\partial t}(\hat{\mathbf{A}}^t)_{\beta} \cdot \left(\frac{\partial \hat{\boldsymbol{u}}^t}{\partial \xi ^\mu} (\hat{J}^{-1})^\mu_{\cdot \alpha} \right) \\
    +& (\hat{\mathbf{A}}^t)_{\alpha} \cdot \left(\frac{\partial \hat{\boldsymbol{u}}^t}{\partial \xi ^\mu} \frac{\partial}{\partial t}(\hat{J}^{-1})^\mu_{\cdot \beta} \right)  + (\hat{\mathbf{A}}^t)_{\beta} \cdot \left(\frac{\partial \hat{\boldsymbol{u}}^t}{\partial \xi ^\mu} \frac{\partial}{\partial t}(\hat{J}^{-1})^\mu_{\cdot \alpha} \right) \Bigg)\Bigg|_{t=0}.
\end{aligned}
\end{equation}

The derived result is obtained by using $\frac{\partial}{\partial t}(\hat{J}^{-1})^\mu_{\cdot \alpha} = -\hat{\vartheta}^\mu_{,\alpha}$ and noting that
\begin{equation}
     \frac{\partial}{\partial t}(\hat{\mathbf{A}}^t)_{\alpha}\Bigg|_{t=0}=  \frac{\partial}{\partial \xi^\lambda} (\hat{\mathbf{A}}^t)_{\alpha}\frac{\partial \xi^\lambda}{\partial(\xi_t)^\mu} \frac{\partial (\xi_t)^\mu}{\partial t}\Bigg|_{t=0} = \mathbf{A}_{\alpha, \mu} \hat{\vartheta}^\mu
\end{equation}

\textbf{Bending strain time derivative}

Moreover, the time derivation of $(\hat{\kappa}^t)_{\alpha\beta} $ yields:
{\footnotesize
\begin{equation}
\begin{aligned}
    \left(\frac{\partial}{\partial t}(\hat{\kappa}^t)_{\alpha\beta}\right)\Bigg|_{t=0} 
     &=  && -(\hat{\mathbf{a}}^t)_3\cdot \frac{\partial ^2  \hat{\boldsymbol{u}}^t}{\partial \xi^\gamma \partial \xi^\eta} \left( -\frac{\partial \hat{\vartheta}^\gamma}{\partial \xi^\beta} \delta^\eta_\alpha- \delta^\gamma_\beta\frac{\partial \hat{\vartheta}^\eta}{\partial \xi^\alpha} \right) \\ & &&+(\hat{\mathbf{A}}^t)_3\cdot \frac{\partial \hat{\boldsymbol{u}}^t}{\partial \xi^\mu} \left( \delta^\mu_\lambda \delta^\gamma_\beta \delta^\eta_\alpha \frac{\partial^2 \hat{\vartheta}^\lambda}{\partial \xi^\gamma \partial \xi^\eta}  + (\ldots) \hat{H}^\lambda_{\gamma\eta}\bigg|_{t=0}\right) \\& && + \left( (\hat{\mathbf{A}}^t)_1 \times \left(\frac{\partial \hat{\boldsymbol{u}}^t}{\partial \xi^\mu } \frac{\partial \hat{\vartheta}^\mu}{\partial \xi^2} \right)  - (\hat{\mathbf{A}}^t)_2 \times \left(\frac{\partial \hat{\boldsymbol{u}}^t}{\partial \xi^\mu } \frac{\partial \hat{\vartheta}^\mu}{\partial \xi^1}\right)\right) \cdot \frac{(\hat{\mathbf{A}}^t)_{\alpha,\beta}}{|(\hat{\mathbf{A}}^t)_1 \times (\hat{\mathbf{A}}^t)_2|},
\end{aligned}
\end{equation}
}
which is reduced to:
{\footnotesize
\begin{equation}
\begin{aligned}
    \left(\frac{\partial}{\partial t}(\hat{\kappa}^t)_{\alpha\beta}\right)\Bigg|_{t=0} 
     &=  && (\hat{\mathbf{A}}^t)_3\cdot  \left( \frac{\partial ^2  \hat{\boldsymbol{u}}^t}{\partial \xi^\gamma \partial \xi^\alpha}\frac{\partial \hat{\vartheta}^\gamma}{\partial \xi^\beta} +  \frac{\partial ^2  \hat{\boldsymbol{u}}^t}{\partial \xi^\gamma \partial \xi^\beta}\frac{\partial \hat{\vartheta}^\gamma}{\partial \xi^\alpha} +  \frac{\partial \hat{\boldsymbol{u}}^t}{ \partial \xi^\gamma}\frac{\partial^2 \hat{\vartheta}^\gamma}{\partial \xi^\alpha \partial\xi^\beta } \right) \\& && + \left( (\hat{\mathbf{A}}^t)_1 \times \left(\frac{\partial \hat{\boldsymbol{u}}^t}{\partial \xi^\mu } \frac{\partial \hat{\vartheta}^\mu}{\partial \xi^2} \right)  - (\hat{\mathbf{A}}^t)_2 \times \left(\frac{\partial \hat{\boldsymbol{u}}^t}{\partial \xi^\mu } \frac{\partial \hat{\vartheta}^\mu}{\partial \xi^1}\right)\right) \cdot \frac{(\hat{\mathbf{A}}^t)_{\alpha,\beta}}{|(\hat{\mathbf{A}}^t)_1 \times (\hat{\mathbf{A}}^t)_2|}.
\end{aligned}
\end{equation}
}



\section{Generation of GRFs by SPDE Sampling}
\label{sec:grfWithSPDE}
It has been shown in \cite{whittle_stochastic_1963, lindgren_explicit_2011,lindgren_spde_2022} that a Gaussian random field can be computed by sampling from a stochastic partial differential equation, denoted as Whittle SPDE. 
This approach is studied and applied in various publications, see for example \cite{lindgren_spde_2022,khristenko_analysis_2019,chu_stochastic_2021,duswald_finite_2024, harbrecht2026simulation}. 
The resulting random field $u: \mathcal{D}\to \mathbb{R}$ is found as the solution to the fractional SPDE, called Whittle SPDE: 
\begin{equation}
(1-\kappa^{-2}\Delta)^k u(\boldsymbol{x}; \omega)=\eta \mathcal{W}(\boldsymbol{x}; \omega), 
\label{eqn:whittle}
\end{equation}
with fractional parameter $k$, where $\omega\in\Omega$ is a realization from an abstract probability space $(\Omega, \mathcal{F}, \mathbb{P})$ (see \cite[Section 2.1]{torquato_random_2002}), and  $\mathcal{W} \sim \mathcal{N}(0,1)$ denotes a realization of Gaussian white noise.  
This equation provides a field with smoothness defined through the exponents $k$ with $\gamma=2k-d/2$. Above, the parameter $\eta > 0$ is a normalizing constant given by
\begin{equation}
    \eta^2 = \sigma^2 \frac{ (4\pi)^{d/2} \Gamma(\gamma+d/2)}{\kappa^d \Gamma(\gamma)}.
\end{equation}
The fractional exponent $k$ defines the smoothness of the resulting random field. 

\begin{remark}
    \textit{Whittle has shown in \cite{whittle_stationary_1954} that when solving Equation~\ref{eqn:whittle} on an unbounded domain, the resulting field will have Matérn covariance. 
    In numerical practice, however, the equation is solved on bounded domains and appropriate boundary conditions need to be imposed. 
    We remark that this approximation will lead to errors with respect to the Matérn covariance close to the boundary. 
    For more details on errors arising from the choice of boundary conditions and the SPDE approximation, we refer to \cite{khristenko_analysis_2019}.}
\end{remark}

For most engineering applications, random fields with general smoothness are not required and we thus limit ourselves to the integer case $k\in\mathbb{N}$, however, we refer to \cite{duswald_finite_2024} and the references therein for a method to compute GRFs with general smoothness. 
The SPDE from Equation~(\ref{eqn:whittle}) can be discretized numerically through different approaches, in the later sections, we detail the generation based on a \emph{Finite Element} discretization and \emph{multi-patch IGA}. 

As different shell-based geometries are studied, Equation~\ref{eqn:whittle} is solved in the case $d=2$, to later study $\mathcal{D}$ as a 2-manifold. For integer values of $k$, the following approaches are distinguished. 

\textit{Harmonic equation ($k=1$, $\gamma=1$).}  Widely studied in the context of standard FEM (see for instance \cite{lindgren_explicit_2011,koh_stochastic_2023,chu_stochastic_2021,duswald_finite_2024}), a first-order polynomial functional space suffices to weakly solve the following Poisson-type equation

\begin{equation}
u(\boldsymbol{x};\omega)-\kappa^{-2}\Delta u(\boldsymbol{x}; \omega)=\eta \mathcal{W}(\boldsymbol{x};\omega). 
\label{eqn:harmonicSPDE}
\end{equation}

\textit{Biharmonic equation ($k=2$, $\gamma=3$).} Leveraging the higher regularity of NURBS basis functions, higher smoothness can be achieved by considering the biharmonic equation
    
\begin{equation}
u(\boldsymbol{x};\omega)-2\kappa^{-2}\Delta u(\boldsymbol{x};\omega) + \kappa^{-4}\Delta^2 u(\boldsymbol{x};\omega)=\eta \mathcal{W}(\boldsymbol{x};\omega). 
\label{eqn:biharmonicSPDE}
\end{equation}

\paragraph*{Generation with IGA}
The two SPDEs \ref{eqn:harmonicSPDE} and \ref{eqn:biharmonicSPDE} can directly be discretized with multi-patch IGA by using the SIPG coupling as applied for the multi patch regularization problem.
Applying this for a simple multi-patch geometry, the resulting fields for the harmonic and biharmonic SPDEs with associated smoothnesses of $\gamma=1$ and $\gamma=3$ are compared.

\begin{figure}[H]
    \centering
    \includegraphics[width=0.87\linewidth]{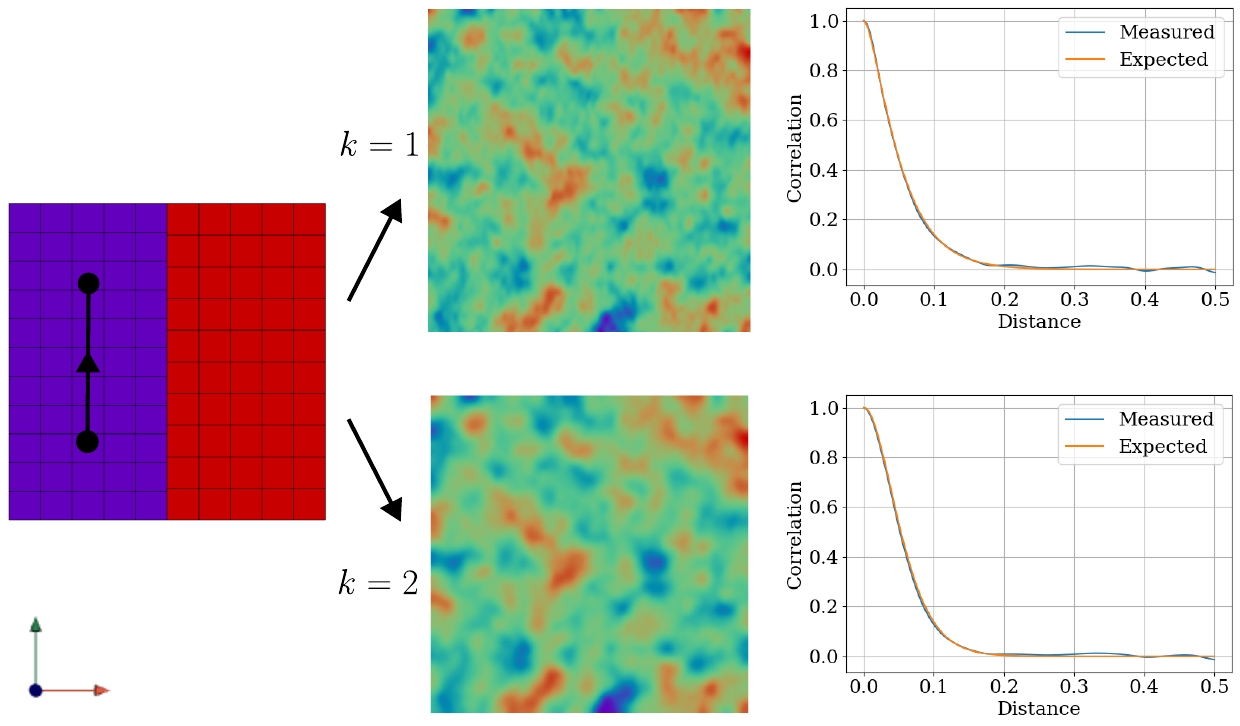}
    \caption{Left: Rectangular NURBS patches with the path used for computing the correlation function. Middle: generated sample of Gaussian random fields for harmonic (top) and biharmonic (bottom). Right: averaged covariance function observed for a sample size of $4\cdot10^4$ in the respective equations with the ground truth Matérn covariance. }
    \label{fig:example01}
\end{figure}

\paragraph*{Comparison of harmonic and biharmonic GRF generation} This example presents a comparison between the Gaussian random fields generated with a harmonic and a biharmonic equation, respectively, on a non-matching two patch domain. The geometry is a unit square, constructed with two cubic NURBS surfaces representing a $0.5\times1$ rectangle. We recall here that the smoothness of the Matérn covariance function is predetermined by $k$, with the relation $\gamma = 2k-d/2$, i.e. $\gamma=1$ for harmonic problem and $\gamma=3$ for the biharmonic, and the correlation length is set to $\rho=0.1$. In Figure \ref{fig:example01}, the two numerically  generated GRFs are displayed. On the left, the patches used for the simulation are shown, along with a path used to validate the observed covariance function. 

\end{appendices}
\newpage
\bibliographystyle{elsarticle-num}
\bibliography{bibliography_zotero}

@book{abraham_manifolds_1988,
  title = {Manifolds, {{Tensor Analysis}}, and {{Applications}}},
  author = {Abraham, Ralph and Marsden, Jerrold E. and Ratiu, Tudor},
  editor = {Marsden, J. E. and Sirovich, L. and John, F.},
  year = 1988,
  series = {Applied {{Mathematical Sciences}}},
  volume = {75},
  publisher = {Springer},
  address = {New York, NY},
  doi = {10.1007/978-1-4612-1029-0},
  copyright = {http://www.springer.com/tdm},
  isbn = {978-1-4612-6990-8 978-1-4612-1029-0}
}

@article{agnelli_design_2022,
  title = {Design of Thin Micro-Architectured Panels with Extension--Bending Coupling Effects Using Topology Optimization},
  author = {Agnelli, Filippo and Nika, Grigor and Constantinescu, Andrei},
  year = 2022,
  month = mar,
  journal = {Comput. Methods Appl. Mech. Eng.},
  volume = {391},
  pages = {114496},
  doi = {10.1016/j.cma.2021.114496}
}

@book{allaire_conception_2007,
  title = {Conception {{Optimale}} de {{Structures}}},
  author = {Allaire, Gr{\'e}goire},
  year = 2007,
  month = jan,
  publisher = {Springer Berlin Heidelberg},
  address = {Berlin Heidelberg},
  doi = {10.1007/978-3-540-36856-4}
}

@article{allaire_level-set_2005,
  title = {A Level-Set Method for Vibration and Multiple Loads Structural Optimization},
  author = {Allaire, Gr{\'e}goire and Jouve, Fran{\c c}ois},
  year = 2005,
  journal = {Comput. Methods Appl. Mech. Eng.},
  volume = {194},
  number = {30},
  pages = {3269--3290},
  doi = {10.1016/j.cma.2004.12.018}
}

@article{allaire_linearized_2014,
  title = {A Linearized Approach to Worst-Case Design in Parametric and Geometric Shape Optimization},
  author = {Allaire, Gr{\'e}goire and Dapogny, Charles},
  year = 2014,
  month = oct,
  journal = {Math. Models Methods Appl. Sci.},
  volume = {24},
  number = {11},
  pages = {2199--2257},
  doi = {10.1142/S0218202514500195},
  langid = {english}
}

@article{allaire_optimizing_2018,
  title = {Optimizing Supports for Additive Manufacturing},
  author = {Allaire, Gr{\'e}goire and Bogosel, Beniamin},
  year = 2018,
  month = dec,
  journal = {Struct Multidisc Optim},
  volume = {58},
  number = {6},
  pages = {2493--2515},
  doi = {10.1007/s00158-018-2125-x},
  langid = {english}
}

@article{ambati_isogeometric_2018,
  title = {Isogeometric {{Kirchhoff}}--{{Love}} Shell Formulation for Elasto-Plasticity},
  author = {Ambati, Marreddy and Kiendl, Josef and De Lorenzis, Laura},
  year = 2018,
  month = oct,
  journal = {Computer Methods in Applied Mechanics and Engineering},
  volume = {340},
  pages = {320--339},
  doi = {10.1016/j.cma.2018.05.023},
  langid = {english}
}

@article{amstutz_topological_2010,
  title = {Topological Optimization of Structures Subject to {{Von Mises}} Stress Constraints},
  author = {Amstutz, Samuel and Novotny, A.A.},
  year = 2010,
  month = apr,
  journal = {Struct. Multidiscip. Optim.},
  volume = {41},
  pages = {407--420},
  doi = {10.1007/s00158-009-0425-x}
}

@misc{antil_risk-averse_2026,
  title = {Risk-Averse Optimization under Distributional Uncertainty with {{Rockafellian}} Relaxation},
  author = {Antil, Harbir and Bustos, Alonso J. and Carney, Sean P. and Venegas, Benjam{\'i}n},
  year = 2026,
  month = mar,
  number = {arXiv:2604.00226},
  eprint = {2604.00226},
  primaryclass = {math},
  publisher = {arXiv},
  doi = {10.48550/arXiv.2604.00226},
  archiveprefix = {arXiv}
}

@article{bathe_geometric_1980,
  title = {A Geometric and Material Nonlinear Plate and Shell Element},
  author = {Bathe, Klaus-J{\"u}rgen and Bolourchi, Sa{\"i}d},
  year = 1980,
  month = feb,
  journal = {Computers \& Structures},
  volume = {11},
  number = {1-2},
  pages = {23--48},
  doi = {10.1016/0045-7949(80)90144-3},
  copyright = {https://www.elsevier.com/tdm/userlicense/1.0/},
  langid = {english}
}

@article{bazilevs_isogeometric_2006,
  title = {Isogeometric Analysis: Approximation, Stability and Error Estimates for h-Refined Meshes},
  author = {Bazilevs, Yuri and {Beirao da Veiga}, L and Cottrell, J Austin and Hughes, Thomas JR and Sangalli, Giancarlo},
  year = 2006,
  journal = {Math. Models Methods Appl. Sci.},
  volume = {16},
  number = {07},
  pages = {1031--1090}
}

@book{birgin_practical_2014,
  title = {Practical Augmented Lagrangian Methods for Constrained Optimization},
  author = {Birgin, Ernesto G. and Mart{\'i}nez, Jos{\'e} Mario},
  year = 2014,
  series = {Fundamentals of Algorithms},
  number = {FA10},
  publisher = {SIAM},
  address = {Philadelphia},
  isbn = {978-1-61197-335-8},
  langid = {english},
  lccn = {519.6}
}

@article{borden_isogeometric_2011,
  title = {Isogeometric Finite Element Data Structures Based on {{B\'ezier}} Extraction of {{NURBS}}},
  author = {Borden, Michael J. and Scott, Michael A. and Evans, John A. and Hughes, Thomas J. R.},
  year = 2011,
  month = jul,
  journal = {Numerical Meth Engineering},
  volume = {87},
  number = {1-5},
  pages = {15--47},
  doi = {10.1002/nme.2968},
  copyright = {http://onlinelibrary.wiley.com/termsAndConditions\#vor},
  langid = {english}
}

@article{bouclier_efficient_2013,
  title = {Efficient Isogeometric {{NURBS-based}} Solid-Shell Elements: {{Mixed}} Formulation and {{B}} - Method},
  shorttitle = {Efficient Isogeometric {{NURBS-based}} Solid-Shell Elements},
  author = {Bouclier, Robin and Elguedj, Thomas and Combescure, Alain},
  year = 2013,
  month = dec,
  journal = {Computer Methods in Applied Mechanics and Engineering},
  volume = {267},
  pages = {86--110},
  doi = {10.1016/j.cma.2013.08.002},
  langid = {english}
}

@article{cai_integrated_2023,
  title = {An Integrated Design Approach for Simultaneous Shape and Topology Optimization of Shell Structures},
  author = {Cai, Shouyu and Zhang, Hualin and Zhang, Weihong},
  year = 2023,
  month = oct,
  journal = {Computer Methods in Applied Mechanics and Engineering},
  volume = {415},
  pages = {116218},
  doi = {10.1016/j.cma.2023.116218},
  langid = {english}
}

@book{chapelle_finite_2010,
  title = {The Finite Element Analysis of Shells-Fundamentals},
  author = {Chapelle, Dominique and Bathe, Klaus-Jurgen},
  year = 2010,
  publisher = {Springer Science \& Business Media}
}

@article{chen_level_2010,
  title = {Level Set Based Robust Shape and Topology Optimization under Random Field Uncertainties},
  author = {Chen, Shikui and Chen, Wei and Lee, Sanghoon},
  year = 2010,
  month = apr,
  journal = {Struct Multidisc Optim},
  volume = {41},
  number = {4},
  pages = {507--524},
  doi = {10.1007/s00158-009-0449-2},
  copyright = {http://www.springer.com/tdm},
  langid = {english}
}

@article{chen_new_2011,
  title = {A New Level-Set Based Approach to Shape and Topology Optimization under Geometric Uncertainty},
  author = {Chen, Shikui and Chen, Wei},
  year = 2011,
  month = jul,
  journal = {Struct Multidisc Optim},
  volume = {44},
  number = {1},
  pages = {1--18},
  doi = {10.1007/s00158-011-0660-9},
  copyright = {http://www.springer.com/tdm},
  langid = {english}
}

@article{chu_stochastic_2021,
  title = {Stochastic {{Modeling}} and Identification of Material Parameters on Structures Produced by Additive Manufacturing},
  author = {Chu, Shanshan and Guilleminot, Johann and Kelly, Cambre and Abar, Bijan and Gall, Ken},
  year = 2021,
  journal = {Comput. Methods Appl. Mech. Eng.},
  volume = {387},
  pages = {114166},
  doi = {10.1016/j.cma.2021.114166}
}

@book{ciarlet_introduction_2006,
  title = {An {{Introduction}} to {{Differential Geometry}} with {{Applications}} to {{Elasticity}}},
  author = {Ciarlet, Philippe G.},
  year = 2006,
  edition = {1st ed},
  publisher = {Springer Netherlands},
  address = {Dordrecht},
  isbn = {978-1-4020-4247-8 978-1-4020-4248-5},
  langid = {english}
}

@incollection{ciarlet_introduction_2008,
  title = {An {{Introduction}} to {{Shell Theory}}},
  booktitle = {Differential {{Geometry}}: {{Theory}} and {{Applications}}},
  author = {Ciarlet, Philippe G. and Mardare, Cristinel},
  year = 2008,
  volume = {9},
  pages = {94--184},
  publisher = {CO-PUBLISHED WITH HIGHER EDUCATION PRESS},
  doi = {10.1142/9789812771476_0002}
}

@article{cirak_subdivision_2000,
  title = {Subdivision Surfaces: A New Paradigm for Thin-Shell Finite-Element Analysis},
  shorttitle = {Subdivision Surfaces},
  author = {Cirak, Fehmi and Ortiz, Michael and Schr{\"o}der, Peter},
  year = 2000,
  month = apr,
  journal = {Int. J. Numer. Meth. Engng.},
  volume = {47},
  number = {12},
  pages = {2039--2072},
  doi = {10.1002/(SICI)1097-0207(20000430)47:12<2039::AID-NME872>3.0.CO;2-1},
  copyright = {http://doi.wiley.com/10.1002/tdm\_license\_1.1},
  langid = {english}
}

@article{conti_risk_2011,
  title = {Risk {{Averse Shape Optimization}}},
  author = {Conti, Sergio and Held, Harald and Pach, Martin and Rumpf, Martin and Schultz, R{\"u}diger},
  year = 2011,
  month = jan,
  journal = {SIAM J. Control Optim.},
  volume = {49},
  number = {3},
  pages = {927--947},
  doi = {10.1137/090754315},
  langid = {english}
}

@article{conti_stochastic_2018,
  title = {Stochastic {{Dominance Constraints}} in {{Elastic Shape Optimization}}},
  author = {Conti, Sergio and Rumpf, Martin and Schultz, R{\"u}diger and T{\"o}lkes, Sascha},
  year = 2018,
  month = jan,
  journal = {SIAM J. Control Optim.},
  volume = {56},
  number = {4},
  pages = {3021--3034},
  doi = {10.1137/16M108313X},
  langid = {english}
}

@book{cottrell_isogeometric_2009,
  title = {Isogeometric {{Analysis}}: {{Toward Integration}} of {{CAD}} and {{FEA}}},
  shorttitle = {Isogeometric {{Analysis}}},
  author = {Cottrell, J. Austin and Hughes, Thomas J. R. and Bazilevs, Yuri},
  year = 2009,
  month = aug,
  edition = {1},
  publisher = {Wiley},
  doi = {10.1002/9780470749081},
  copyright = {http://doi.wiley.com/10.1002/tdm\_license\_1.1},
  isbn = {978-0-470-74873-2 978-0-470-74908-1},
  langid = {english}
}

@article{croci_efficient_2018,
  title = {Efficient {{White Noise Sampling}} and {{Coupling}} for {{Multilevel Monte Carlo}} with {{Nonnested Meshes}}},
  author = {Croci, M. and Giles, M. B. and Rognes, M. E. and Farrell, P. E.},
  year = 2018,
  month = jan,
  journal = {SIAM/ASA J. Uncertainty Quantification},
  volume = {6},
  number = {4},
  pages = {1630--1655},
  doi = {10.1137/18M1175239},
  langid = {english}
}

@misc{dapogny_distributionally_2025,
  title = {Distributionally {{Robust Shape}} and {{Topology Optimization}}},
  author = {Dapogny, Charles and Prando, Julien and Thibert, Boris},
  year = 2025,
  month = jul,
  number = {arXiv:2507.21574},
  eprint = {2507.21574},
  primaryclass = {math},
  publisher = {arXiv},
  doi = {10.48550/arXiv.2507.21574},
  archiveprefix = {arXiv}
}

@article{de_gournay_shape_2008,
  title = {Shape and Topology Optimization of the Robust Compliance {\emph{via}} the Level Set Method},
  author = {De Gournay, Fr{\'e}d{\'e}ric and Allaire, Gr{\'e}goire and Jouve, Fran{\c c}ois},
  year = 2008,
  month = jan,
  journal = {ESAIM: COCV},
  volume = {14},
  number = {1},
  pages = {43--70},
  doi = {10.1051/cocv:2007048}
}

@article{de_gournay_velocity_2006,
  title = {Velocity {{Extension}} for the {{Level-set Method}} and {{Multiple Eigenvalues}} in {{Shape Optimization}}},
  author = {De Gournay, Fr{\'e}d{\'e}ric},
  year = 2006,
  month = jan,
  journal = {SIAM J. Control Optim.},
  volume = {45},
  number = {1},
  pages = {343--367},
  doi = {10.1137/050624108},
  langid = {english}
}

@article{de_topology_2020,
  title = {Topology Optimization under Uncertainty Using a Stochastic Gradient-Based Approach},
  author = {De, Subhayan and Hampton, Jerrad and Maute, Kurt and Doostan, Alireza},
  year = 2020,
  month = nov,
  journal = {Struct Multidisc Optim},
  volume = {62},
  number = {5},
  pages = {2255--2278},
  doi = {10.1007/s00158-020-02599-z},
  langid = {english}
}

@article{dede_isogeometric_2012,
  title = {Isogeometric {{Analysis}} for {{Topology Optimization}} with a {{Phase Field Model}}},
  author = {Ded{\`e}, Luca and Borden, Micheal J. and Hughes, Thomas J. R.},
  year = 2012,
  month = sep,
  journal = {Arch. Comput. Methods Eng.},
  volume = {19},
  number = {3},
  pages = {427--465},
  doi = {10.1007/s11831-012-9075-z}
}

@book{delfour_shapes_2011,
  title = {Shapes and {{Geometries}}: {{Metrics}}, {{Analysis}}, {{Differential Calculus}}, and {{Optimization}}, {{Second Edition}}},
  shorttitle = {Shapes and {{Geometries}}},
  author = {Delfour, M. C. and Zol{\'e}sio, J. -P.},
  year = 2011,
  month = jan,
  edition = {Second},
  publisher = {{Society for Industrial and Applied Mathematics}},
  doi = {10.1137/1.9780898719826},
  isbn = {978-0-89871-936-9 978-0-89871-982-6},
  langid = {english}
}

@misc{diaz-avalos_formopt_2026,
  title = {{{FormOpt}}: {{A FEniCSx}} Toolbox for Level Set-Based Shape Optimization Supporting Parallel Computing},
  shorttitle = {{{FormOpt}}},
  author = {{D{\'i}az-Avalos}, Josu{\'e} D. and Laurain, Antoine},
  year = 2026,
  month = jan,
  number = {arXiv:2601.05709},
  eprint = {2601.05709},
  primaryclass = {math},
  publisher = {arXiv},
  doi = {10.48550/arXiv.2601.05709},
  archiveprefix = {arXiv}
}

@article{duswald_finite_2024,
  title = {Finite Elements for {{Mat\'ern-type}} Random Fields: {{Uncertainty}} in Computational Mechanics and Design Optimization},
  shorttitle = {Finite Elements for {{Mat\'ern-type}} Random Fields},
  author = {Duswald, Tobias and Keith, Brendan and Lazarov, Boyan and Petrides, Socratis and Wohlmuth, Barbara},
  year = 2024,
  month = sep,
  journal = {Computer Methods in Applied Mechanics and Engineering},
  volume = {429},
  pages = {117146},
  doi = {10.1016/j.cma.2024.117146},
  langid = {english}
}

@article{epshteyn_estimation_2007,
  title = {Estimation of Penalty Parameters for Symmetric Interior Penalty {{Galerkin}} Methods},
  author = {Epshteyn, Yekaterina and Rivière, Béatrice},
  year = {2007},
  journaltitle = {Journal of Computational and Applied Mathematics},
  journal = {Journal of Computational and Applied Mathematics},
  volume = {206},
  number = {2},
  pages = {843--872},
  doi = {10.1016/j.cam.2006.08.029},
  langid = {english}
}

@book{frankel_geometry_2004,
  title = {The {{Geometry}} of {{Physics}}: {{An Introduction}}},
  shorttitle = {The {{Geometry}} of {{Physics}}},
  author = {Frankel, Theodore},
  year = 2004,
  publisher = {Cambridge University Press},
  googlebooks = {2iq2EGNgX24C},
  isbn = {978-0-521-53927-2},
  langid = {english}
}

@article{guo_confidence_2009,
  title = {Confidence Structural Robust Design and Optimization under Stiffness and Load Uncertainties},
  author = {Guo, Xu and Bai, Wei and Zhang, Weisheng and Gao, Xixin},
  year = 2009,
  month = sep,
  journal = {Computer Methods in Applied Mechanics and Engineering},
  volume = {198},
  number = {41-44},
  pages = {3378--3399},
  doi = {10.1016/j.cma.2009.06.018},
  copyright = {https://www.elsevier.com/tdm/userlicense/1.0/},
  langid = {english}
}

@article{hennig_bezier_2016,
  title = {B\'ezier Extraction and Adaptive Refinement of Truncated Hierarchical {{NURBS}}},
  author = {Hennig, P. and M{\"u}ller, S. and K{\"a}stner, M.},
  year = 2016,
  month = jun,
  journal = {Computer Methods in Applied Mechanics and Engineering},
  volume = {305},
  pages = {316--339},
  doi = {10.1016/j.cma.2016.03.009},
  langid = {english}
}

@article{herrema_penalty_2019,
  title = {Penalty Coupling of Non-Matching Isogeometric {{Kirchhoff}}--{{Love}} Shell Patches with Application to Composite Wind Turbine Blades},
  author = {Herrema, Austin J. and Johnson, Emily L. and Proserpio, Davide and Wu, Michael C.H. and Kiendl, Josef and Hsu, Ming-Chen},
  year = 2019,
  month = apr,
  journal = {Computer Methods in Applied Mechanics and Engineering},
  volume = {346},
  pages = {810--840},
  doi = {10.1016/j.cma.2018.08.038},
  langid = {english}
}

@article{hirschler_embedded_2019,
  title = {The Embedded Isogeometric {{Kirchhoff}}--{{Love}} Shell: {{From}} Design to Shape Optimization of Non-Conforming Stiffened Multipatch Structures},
  shorttitle = {The Embedded Isogeometric {{Kirchhoff}}--{{Love}} Shell},
  author = {Hirschler, T. and Bouclier, R. and Duval, A. and Elguedj, T. and Morlier, J.},
  year = 2019,
  month = jun,
  journal = {Computer Methods in Applied Mechanics and Engineering},
  volume = {349},
  pages = {774--797},
  doi = {10.1016/j.cma.2019.02.042},
  langid = {english}
}

@article{hsu_high-performance_2011,
  title = {High-Performance Computing of Wind Turbine Aerodynamics Using Isogeometric Analysis},
  author = {Hsu, Ming-Chen and Akkerman, Ido and Bazilevs, Yuri},
  year = 2011,
  month = oct,
  journal = {Computers \& Fluids},
  volume = {49},
  number = {1},
  pages = {93--100},
  doi = {10.1016/j.compfluid.2011.05.002},
  copyright = {https://www.elsevier.com/tdm/userlicense/1.0/},
  langid = {english}
}

@article{hubner_scherer_stress-based_2024,
  title = {Stress-Based Topological Shape Optimization for Thick Shells Using the Level Set Method and Trimmed Non-Conforming Multi-Patch Isogeometric Analysis},
  author = {H{\"u}bner Scherer, Fernando and Zarroug, Malek and Naceur, Hakim and Constantinescu, Andrei},
  year = 2024,
  month = oct,
  journal = {Struct Multidisc Optim},
  volume = {67},
  number = {10},
  pages = {177},
  doi = {10.1007/s00158-024-03892-x},
  langid = {english}
}

@article{hubner_scherer_topology_2024,
  title = {Topology Optimization of Curved Thick Shells Using Level Set Method and Non-Conforming Multi-Patch Isogeometric Analysis},
  author = {H{\"u}bner Scherer, Fernando and Zarroug, Malek and Naceur, Hakim and Constantinescu, Andrei},
  year = 2024,
  month = oct,
  journal = {Computer Methods in Applied Mechanics and Engineering},
  volume = {430},
  pages = {117205},
  doi = {10.1016/j.cma.2024.117205},
  langid = {english}
}

@article{hughes_isogeometric_2005,
  title = {Isogeometric Analysis: {{CAD}}, Finite Elements, {{NURBS}}, Exact Geometry and Mesh Refinement},
  author = {Hughes, T. J. R. and Cottrell, J. A. and Bazilevs, Y.},
  year = 2005,
  journal = {Comput. Methods Appl. Mech. Eng.},
  volume = {194},
  number = {39},
  pages = {4135--4195},
  doi = {10.1016/j.cma.2004.10.008}
}

@article{jahangiry_isogeometrical_2017,
  title = {An Isogeometrical Approach to Structural Level Set Topology Optimization},
  author = {Jahangiry, Hassan A. and Tavakkoli, S. Mehdi},
  year = 2017,
  month = jun,
  journal = {Computer Methods in Applied Mechanics and Engineering},
  volume = {319},
  pages = {240--257},
  doi = {10.1016/j.cma.2017.02.005},
  langid = {english}
}

@article{kamensky_tigar_2019,
  title = {{{tIGAr}}: {{Automating}} Isogeometric Analysis with {{FEniCS}}},
  author = {Kamensky, David and Bazilevs, Yuri},
  year = 2019,
  journal = {Comput. Methods Appl. Mech. Eng.},
  volume = {344},
  pages = {477--498},
  doi = {10.1016/j.cma.2018.10.002}
}

@article{kharmanda_reliability-based_2004,
  title = {Reliability-Based Topology Optimization},
  author = {Kharmanda, G. and Olhoff, N. and Mohamed, A. and Lemaire, M.},
  year = 2004,
  month = mar,
  journal = {Structural and Multidisciplinary Optimization},
  volume = {26},
  number = {5},
  pages = {295--307},
  doi = {10.1007/s00158-003-0322-7},
  copyright = {http://www.springer.com/tdm}
}

@article{khristenko_analysis_2019,
  title = {Analysis of {{Boundary Effects}} on {{PDE-Based Sampling}} of {{Whittle--Mat\'ern Random Fields}}},
  author = {Khristenko, U. and Scarabosio, L. and Swierczynski, P. and Ullmann, E. and Wohlmuth, B.},
  year = 2019,
  month = jan,
  journal = {SIAM/ASA J. Uncertainty Quantification},
  volume = {7},
  number = {3},
  pages = {948--974},
  doi = {10.1137/18M1215700},
  langid = {english}
}

@article{khristenko_statistical_2020,
  title = {A {{Statistical Framework}} for {{Generating Microstructures}} of {{Two-Phase Random Materials}}: {{Application}} to {{Fatigue Analysis}}},
  author = {Khristenko, Ustim and Constantinescu, Andrei and Tallec, Patrick Le and Oden, J. Tinsley and Wohlmuth, Barbara},
  year = 2020,
  journal = {Multiscale Model. Simul.},
  volume = {18},
  number = {1},
  pages = {21--43},
  doi = {10.1137/19M1259286}
}

@article{khristenko_statistically_2022,
  title = {Statistically Equivalent Surrogate Material Models: {{Impact}} of Random Imperfections on the Elasto-Plastic Response},
  author = {Khristenko, Ustim and Constantinescu, Andrei and Tallec, Patrick Le and Wohlmuth, Barbara},
  year = 2022,
  journal = {Comput. Methods Appl. Mech. Eng.},
  volume = {402},
  pages = {115278},
  doi = {10.1016/j.cma.2022.115278}
}

@article{kiendl_isogeometric_2009,
  title = {Isogeometric Shell Analysis with {{Kirchhoff}}--{{Love}} Elements},
  author = {Kiendl, J. and Bletzinger, K.-U. and Linhard, J. and W{\"u}chner, R.},
  year = 2009,
  journal = {Comput. Methods Appl. Mech. Eng.},
  volume = {198},
  number = {49},
  pages = {3902--3914},
  doi = {10.1016/j.cma.2009.08.013}
}

@article{kiendl_isogeometric_2014,
  title = {Isogeometric Shape Optimization of Shells Using Semi-Analytical Sensitivity Analysis and Sensitivity Weighting},
  author = {Kiendl, J. and Schmidt, R. and W{\"u}chner, R. and Bletzinger, K.-U.},
  year = 2014,
  month = jun,
  journal = {Computer Methods in Applied Mechanics and Engineering},
  volume = {274},
  pages = {148--167},
  doi = {10.1016/j.cma.2014.02.001},
  langid = {english}
}

@article{kiendl_phase-field_2016,
  title = {Phase-Field Description of Brittle Fracture in Plates and Shells},
  author = {Kiendl, Josef and Ambati, Marreddy and De Lorenzis, Laura and Gomez, Hector and Reali, Alessandro},
  year = 2016,
  month = dec,
  journal = {Computer Methods in Applied Mechanics and Engineering},
  volume = {312},
  pages = {374--394},
  doi = {10.1016/j.cma.2016.09.011},
  langid = {english}
}

@article{kodakkal_risk-averse_2022,
  title = {Risk-Averse Design of Tall Buildings for Uncertain Wind Conditions},
  author = {Kodakkal, Anoop and Keith, Brendan and Khristenko, Ustim and Apostolatos, Andreas and Bletzinger, Kai-Uwe and Wohlmuth, Barbara and W{\"u}chner, Roland},
  year = 2022,
  month = dec,
  journal = {Computer Methods in Applied Mechanics and Engineering},
  volume = {402},
  pages = {115371},
  doi = {10.1016/j.cma.2022.115371},
  langid = {english}
}

@article{koh_stochastic_2023,
  title = {Stochastic {{PDE}} Representation of Random Fields for Large-Scale {{Gaussian}} Process Regression and Statistical Finite Element Analysis},
  author = {Koh, Kim Jie and Cirak, Fehmi},
  year = 2023,
  journal = {Comput. Methods Appl. Mech. Eng.},
  volume = {417},
  pages = {116358},
  doi = {10.1016/j.cma.2023.116358}
}

@article{laurain_distributed_2015,
  title = {Distributed Shape Derivative via Averaged Adjoint Method and Applications},
  author = {Laurain, Antoine and Sturm, Kevin},
  year = 2015,
  month = sep,
  journal = {ESAIM Math. Model. Numer. Anal.},
  volume = {50},
  doi = {10.1051/m2an/2015075}
}

@article{laurain_level_2017,
  title = {A Level Set-Based Structural Optimization Code Using {{FEniCS}}},
  author = {Laurain, Antoine},
  year = 2018,
  month = may,
  journal = {Struct. Multidiscip. Optim.},
  volume = {58},
  doi = {10.1007/s00158-018-1950-2}
}

@article{lindgren_explicit_2011,
  title = {An Explicit Link between {{Gaussian}} Fields and {{Gaussian Markov}} Random Fields: The Stochastic Partial Differential Equation Approach},
  author = {Lindgren, Finn and Rue, H{\aa}vard and Lindstr{\"o}m, Johan},
  year = 2011,
  journal = {J. R. Stat. Soc. Ser. B Stat. Methodol.},
  volume = {73},
  number = {4},
  pages = {423--498},
  doi = {10.1111/j.1467-9868.2011.00777.x}
}

@article{lindgren_spde_2022,
  title = {The {{SPDE}} Approach for {{Gaussian}} and Non-{{Gaussian}} Fields: 10 Years and Still Running},
  author = {Lindgren, Finn and Bolin, David and Rue, H{\aa}vard},
  year = 2022,
  journal = {Spat. Stat.},
  volume = {50},
  pages = {100599},
  doi = {10.1016/j.spasta.2022.100599}
}

@book{logg_automated_2012,
  title = {Automated {{Solution}} of {{Differential Equations}} by the {{Finite Element Method}}: {{The FEniCS Book}}},
  shorttitle = {Automated {{Solution}} of {{Differential Equations}} by the {{Finite Element Method}}},
  editor = {Logg, Anders and Mardal, Kent-Andre and Wells, Garth},
  year = 2012,
  series = {Lecture {{Notes}} in {{Computational Science}} and {{Engineering}}},
  volume = {84},
  publisher = {Springer Berlin Heidelberg},
  address = {Berlin, Heidelberg},
  doi = {10.1007/978-3-642-23099-8},
  copyright = {https://www.springernature.com/gp/researchers/text-and-data-mining},
  isbn = {978-3-642-23098-1 978-3-642-23099-8},
  langid = {english}
}

@article{macneal_proposed_1985,
  title = {A Proposed Standard Set of Problems to Test Finite Element Accuracy},
  author = {Macneal, Richard H. and Harder, Robert L.},
  year = 1985,
  journal = {Finite Elem. Anal. Des.},
  volume = {1},
  number = {1},
  pages = {3--20},
  doi = {10.1016/0168-874X(85)90003-4}
}

@article{martinez-frutos_risk-averse_2021,
  title = {Risk-Averse Approach for Topology Optimization of Fail-Safe Structures Using the Level-Set Method},
  author = {{Mart{\'i}nez-Frutos}, J. and Ortigosa, R.},
  year = 2021,
  month = nov,
  journal = {Comput Mech},
  volume = {68},
  number = {5},
  pages = {1039--1061},
  doi = {10.1007/s00466-021-02058-6},
  langid = {english}
}

@article{martinez-frutos_robust_2016,
  title = {Robust Shape Optimization of Continuous Structures via the Level Set Method},
  author = {{Mart{\'i}nez-Frutos}, Jes{\'u}s and {Herrero-P{\'e}rez}, David and Kessler, Mathieu and Periago, Francisco},
  year = 2016,
  journal = {Comput. Methods Appl. Mech. Eng.},
  volume = {305},
  pages = {271--291},
  doi = {10.1016/j.cma.2016.03.003}
}

@article{moore_discontinuous_2018,
  title = {Discontinuous {{Galerkin Isogeometric Analysis}} for the Biharmonic Equation},
  author = {Moore, Stephen Edward},
  year = 2018,
  journal = {Comput. Math. Appl.},
  volume = {76},
  number = {4},
  pages = {673--685},
  doi = {10.1016/j.camwa.2018.05.001}
}

@misc{moore_multipatch_2020,
  title = {Multipatch {{Discontinuous Galerkin IGA}} for the {{Biharmonic Problem On Surfaces}}},
  author = {Moore, Stephen E.},
  year = 2020,
  publisher = {arXiv},
  doi = {10.48550/ARXIV.2012.03425},
  copyright = {Creative Commons Attribution 4.0 International}
}

@article{nika_design_2019,
  title = {Design of Multi-Layer Materials Using Inverse Homogenization and a Level Set Method},
  author = {Nika, Grigor and Constantinescu, Andrei},
  year = 2019,
  month = apr,
  journal = {Comput. Methods Appl. Mech. Eng.},
  volume = {346},
  pages = {388--409},
  doi = {10.1016/j.cma.2018.11.029}
}

@book{nocedal_numerical_2006,
  title = {Numerical {{Optimization}}},
  author = {Nocedal, Jorge and Wright, Stephen},
  year = 2006,
  series = {Springer {{Series}} in {{Operations Research}} and {{Financial Engineering}}},
  publisher = {Springer New York},
  doi = {10.1007/978-0-387-40065-5},
  copyright = {http://www.springer.com/tdm},
  isbn = {978-0-387-30303-1},
  langid = {english}
}

@article{proserpio_penalty_2022,
  title = {Penalty Coupling of Trimmed Isogeometric {{Kirchhoff}}--{{Love}} Shell Patches},
  author = {Proserpio, Davide and Kiendl, Josef},
  year = 2022,
  month = mar,
  journal = {J. Mech.},
  volume = {38},
  pages = {156--165},
  doi = {10.1093/jom/ufac008},
  copyright = {https://creativecommons.org/licenses/by/4.0/},
  langid = {english}
}

@article{qian_full_2010,
  title = {Full Analytical Sensitivities in {{NURBS}} Based Isogeometric Shape Optimization},
  author = {Qian, Xiaoping},
  year = 2010,
  month = jun,
  journal = {Computer Methods in Applied Mechanics and Engineering},
  volume = {199},
  number = {29-32},
  pages = {2059--2071},
  doi = {10.1016/j.cma.2010.03.005},
  copyright = {https://www.elsevier.com/tdm/userlicense/1.0/},
  langid = {english}
}

@article{scott_isogeometric_2011,
  title = {Isogeometric Finite Element Data Structures Based on {{B\'ezier}} Extraction of {{T}}-splines},
  author = {Scott, Michael A. and Borden, Michael J. and Verhoosel, Clemens V. and Sederberg, Thomas W. and Hughes, Thomas J. R.},
  year = 2011,
  month = oct,
  journal = {Numerical Meth Engineering},
  volume = {88},
  number = {2},
  pages = {126--156},
  doi = {10.1002/nme.3167},
  copyright = {http://onlinelibrary.wiley.com/termsAndConditions\#vor},
  langid = {english}
}

@inproceedings{sederberg_free-form_1986,
  title = {Free-Form Deformation of Solid Geometric Models},
  booktitle = {Proc. 13th {{Annu}}. {{Conf}}. {{Comput}}. {{Graph}}. {{Interact}}. {{Tech}}.},
  author = {Sederberg, Thomas W and Parry, Scott R},
  year = 1986,
  pages = {151--160}
}

@article{sigmund_manufacturing_2009,
  title = {Manufacturing Tolerant Topology Optimization},
  author = {Sigmund, Ole},
  year = 2009,
  month = apr,
  journal = {Acta Mech Sin},
  volume = {25},
  number = {2},
  pages = {227--239},
  doi = {10.1007/s10409-009-0240-z},
  copyright = {http://www.springer.com/tdm},
  langid = {english}
}

@book{sokolowski_introduction_1992,
  title = {Introduction to {{Shape Optimization}}: {{Shape Sensitivity Analysis}}},
  author = {Sokolowski, Jan and Zolesio, Jean-Paul},
  year = 1992,
  series = {Springer {{Series}} in {{Computational Mathematics}}},
  edition = {1},
  publisher = {Springer Berlin, Heidelberg},
  doi = {10.1007/978-3-642-58106-9},
  isbn = {978-3-540-54177-6}
}

@book{spivak_comprehensive_1979,
  title = {A Comprehensive Introduction to Differential Geometry},
  author = {Spivak, Michael David},
  year = 1979,
  edition = {2nd ed},
  publisher = {Publish or Perish, Inc.}
}

@incollection{stein_models_2004,
  title = {Models and {{Finite Elements}} for {{Thin}}-{{Walled Structures}}},
  booktitle = {Encyclopedia of {{Computational Mechanics}}},
  author = {Bischoff, M. and Bletzinger, K.-U. and Wall, W. A. and Ramm, E.},
  editor = {Stein, Erwin and Borst, Ren{\'e} and Hughes, Thomas J. R.},
  year = 2004,
  month = aug,
  edition = {1},
  publisher = {Wiley},
  doi = {10.1002/0470091355.ecm026},
  isbn = {978-0-470-84699-5 978-0-470-09135-7},
  langid = {english}
}

@misc{teixeira_isogeometric_2025,
  title = {Isogeometric {{Topology Optimization Based}} on {{Topological Derivatives}}},
  author = {Teixeira, Guilherme Henrique and Krenn, Nepomuk and Gangl, Peter and Marussig, Benjamin},
  year = 2025,
  month = oct,
  number = {arXiv:2509.09236},
  eprint = {2509.09236},
  primaryclass = {math},
  publisher = {arXiv},
  doi = {10.48550/arXiv.2509.09236},
  archiveprefix = {arXiv}
}

@article{tepole_isogeometric_2015,
  title = {Isogeometric {{Kirchhoff}}--{{Love}} Shell Formulations for Biological Membranes},
  author = {Tepole, Adri{\'a}n Buganza and Kabaria, Hardik and Bletzinger, Kai-Uwe and Kuhl, Ellen},
  year = 2015,
  month = aug,
  journal = {Computer Methods in Applied Mechanics and Engineering},
  volume = {293},
  pages = {328--347},
  doi = {10.1016/j.cma.2015.05.006},
  langid = {english}
}

@book{torquato_random_2002,
  title = {Random {{Heterogeneous Materials}}: {{Microstructure}} and {{Macroscopic Properties}}},
  author = {Torquato, Salvatore},
  year = 2002,
  series = {Interdisciplinary {{Applied Mathematics}}},
  edition = {1},
  volume = {16},
  publisher = {Springer},
  address = {New York, NY},
  doi = {10.1007/978-1-4757-6355-3},
  isbn = {978-0-387-95167-6}
}

@article{wall_isogeometric_2008,
  title = {Isogeometric Structural Shape Optimization},
  author = {Wall, Wolfgang A. and Frenzel, Moritz A. and Cyron, Christian},
  year = 2008,
  month = jun,
  journal = {Computer Methods in Applied Mechanics and Engineering},
  volume = {197},
  number = {33-40},
  pages = {2976--2988},
  doi = {10.1016/j.cma.2008.01.025},
  copyright = {https://www.elsevier.com/tdm/userlicense/1.0/},
  langid = {english}
}

@article{wang_isogeometric_2016,
  title = {Isogeometric Analysis for Parameterized {{LSM-based}} Structural Topology Optimization},
  author = {Wang, Yingjun and Benson, David},
  year = 2016,
  month = jan,
  journal = {Comput. Mech.},
  volume = {57},
  doi = {10.1007/s00466-015-1219-1}
}

@article{whittle_stationary_1954,
  title = {On Stationary Processes in the Plane},
  author = {Whittle, Peter},
  year = 1954,
  journal = {Biometrika},
  pages = {434--449}
}

@article{whittle_stochastic_1963,
  title = {Stochastic {{Processes}} in {{Several Dimensions}}},
  author = {Whittle, Peter},
  year = 1963,
  journal = {Bull. Int. Stat. Inst.},
  volume = {40},
  pages = {974--994}
}

@article{xing_isogeometric_2026,
  title = {Isogeometric Topology Optimization for Smooth Shells Using a Conformal Mapping-Enhanced Level Set Method},
  author = {Xing, Zhi and Wang, Yingjun},
  year = 2026,
  month = mar,
  journal = {Applied Mathematical Modelling},
  volume = {151},
  pages = {116504},
  doi = {10.1016/j.apm.2025.116504},
  langid = {english}
}

@article{xu_level_2019,
  title = {Level Set-Based Isogeometric Topology Optimization for Maximizing Fundamental Eigenfrequency},
  author = {Xu, Manman and Wang, Shuting and Xie, Xianda},
  year = 2019,
  month = feb,
  journal = {Front. Mech. Eng.},
  volume = {14},
  doi = {10.1007/s11465-019-0534-1}
}

@article{yu_three-dimensional_2019,
  title = {Three-Dimensional Elastoplastic Solids Simulation by an Effective {{IGA}} Based on {{B\'ezier}} Extraction of {{NURBS}}},
  author = {Yu, Tiantang and Lai, Wenjiang and Bui, Tinh Quoc},
  year = 2019,
  month = mar,
  journal = {Int J Mech Mater Des},
  volume = {15},
  number = {1},
  pages = {175--197},
  doi = {10.1007/s10999-018-9405-x},
  langid = {english}
}

@article{zhao_automated_2024,
  title = {Automated Shape and Thickness Optimization for Non-Matching Isogeometric Shells Using Free-Form Deformation},
  author = {Zhao, Han and Kamensky, David and Hwang, John T. and Chen, Jiun-Shyan},
  year = 2024,
  month = dec,
  journal = {Engineering with Computers},
  volume = {40},
  number = {6},
  pages = {3495--3518},
  doi = {10.1007/s00366-024-01947-7},
  langid = {english}
}

@article{zhao_open-source_2022,
  title = {An Open-Source Framework for Coupling Non-Matching Isogeometric Shells with Application to Aerospace Structures},
  author = {Zhao, Han and Liu, Xiangbei and Fletcher, Andrew H. and Xiang, Ru and Hwang, John T. and Kamensky, David},
  year = 2022,
  journal = {Comput. Math. Appl.},
  volume = {111},
  pages = {109--123},
  doi = {10.1016/j.camwa.2022.02.007}
}

@article{zheng_levelset_2019,
  title = {Level-set Topology Optimization for Robust Design of Structures under Hybrid Uncertainties},
  author = {Zheng, Jing and Luo, Zhen and Jiang, Chao and Wu, Jinglai},
  year = 2019,
  month = feb,
  journal = {Numerical Meth Engineering},
  volume = {117},
  number = {5},
  pages = {523--542},
  doi = {10.1002/nme.5966},
  langid = {english}
}

@article{baumgarten2026risk,
  title={Risk-averse Optimization in Random Materials: Algorithmic Advances and HPC Acceleration},
  author={Baumgarten, Niklas and Koch, Marcel and Schneiderhan, David and Schrader, Tim and Wei{\ss}, Robin},
  journal={arXiv preprint arXiv:2608.25967},
  year={2026}
}

@article{harbrecht2026simulation,
  title={Simulation of Gaussian random fields on surfaces using the isogeometric finite element method},
  author={Harbrecht, Helmut and Sonderegger, Florian and von Rickenbach, Remo},
  journal={arXiv preprint arXiv:2605.24882},
  year={2026}
}

\end{document}